\documentclass[a4paper,11pt,final]{amsart}

\usepackage{amsmath,amssymb,amsthm,amscd}

\numberwithin{equation}{section}
\allowdisplaybreaks

\theoremstyle{definition}
 \newtheorem{thm}{Theorem}[section]
 \newtheorem*{thm*}{Theorem}
 \newtheorem{prop}[thm]{Proposition}
 \newtheorem{lem}[thm]{Lemma}
 \newtheorem{cor}[thm]{Corollary}
  \newtheorem{con}[thm]{Conjecture}

\theoremstyle{definition}
 
 \newtheorem{dfn}[thm]{Definition}
 \newtheorem*{ack}{Acknowledgements}
 
 \newtheorem{rmk}[thm]{Remark}

\newcommand{\ep}{\epsilon}
\newcommand{\ve}{\varepsilon}
\newcommand{\bbC}{\mathbb{C}}
\newcommand{\bbF}{\mathbb{F}}
\newcommand{\bbN}{\mathbb{N}}
\newcommand{\bbP}{\mathbb{P}}
\newcommand{\bbQ}{\mathbb{Q}}
\newcommand{\bbZ}{\mathbb{Z}}
\newcommand{\calA}{\mathcal{A}}

\newcommand{\calF}{\mathcal{F}}

\newcommand{\calU}{\mathcal{U}}
\newcommand{\calW}{\mathcal{W}}

\newcommand{\calM}{\mathcal{M}}
\newcommand{\calO}{\mathcal{O}}
\newcommand{\calP}{\mathcal{P}}
\newcommand{\frakh}{\mathfrak{h}}
\newcommand{\frakS}{\mathfrak{S}}
\newcommand{\whcalA}{\widehat{\mathcal{A}}}
\newcommand{\wtcalA}{\widetilde{\mathcal{A}}}
\newcommand{\wbcalA}{\overline{\mathcal{A}}}

\newcommand{\seteq}{\mathbin{:=}}
\newcommand{\Ind}{\operatorname{Ind}}
\newcommand{\Sym}{\operatorname{Sym}}
\newcommand{\sgn}{\operatorname{sgn}}

\newcommand{\qtpr}[1]{\left\langle #1\right\rangle_{q,t}}

\title[A commutative algebra and Macdonald polynomials]{A commutative algebra on degenerate $\bbC\bbP^1$\\ and Macdonald polynomials}

\date{April 15,\ 2009}

\author{B. Feigin, K. Hashizume, A. Hoshino, J. Shiraishi and S. Yanagida}
\address{BF: Landau Institute for Theoretical Physics,
Russia, Chernogolovka, 142432, prosp. Akademika Semenova, 1a,   \newline
Higher School of Economics, Russia, Moscow, 101000,  Myasnitskaya ul., 20 and
\newline
Independent University of Moscow, Russia, Moscow, 119002,
Bol'shoi Vlas'evski per., 11}
\email{bfeigin@gmail.com}
\address{KH, SY: Kobe University, Department of Mathematics, Rokko, Kobe 658-8501, Japan}
\email{hashizume@math.kobe-u.ac.jp, yanagida@math.kobe-u.ac.jp}
\address{AH. Department of Mathematics, Sophia University, Kioicyo, Tokyo, 102-8554, Japan}
\email{ayumu-h@mm.sophia.ac.jp}
\address{JS: Graduate School of Mathematical Sciences, University of Tokyo, Komaba, Tokyo
153-8914, Japan}
\email{shiraish@ms.u-tokyo.ac.jp}

\begin{document}

\begin{abstract}
We introduce a unital associative algebra $\calA$ over 
degenerate $\bbC\bbP^1$. We show that $\calA$ is a commutative algebra
and whose Poincar\'e series is given by the number of partitions. 
Thereby we can regard $\calA$ as a smooth degeneration limit of 
the elliptic algebra introduced by one of the authors and Odesskii \cite{FO:1997}.
Then we study the commutative family of the Macdonald difference operators 
acting on the space of symmetric functions. 
A canonical basis is proposed for this family by
using $\calA$ and 
the Heisenberg representation of the commutative family studied by one of the authors \cite{S:2006}.
It is found that the Ding-Iohara algebra \cite{DI:1997} 
provides us with an algebraic framework for 
the free filed construction. An elliptic deformation
of our construction
is discussed, showing connections with 
the Drinfeld quasi-Hopf twisting \cite{QHA} a la Babelon Bernard Billey \cite{BBB}, 
the Ruijsenaars difference operator \cite {Ruijsenaars:1987}
and the operator ${\sf M}(q,t_1,t_2)$ of Okounkov-Pandharipande \cite{OkounkovPandharipande:2004}.
\end{abstract}

\maketitle

\section{Introduction}\label{sect:intro}
The aim of this paper is to construct a graded algebra $\calA=\calA(q_1,q_2,q_3)$ having three shift parameters $q_1$, $q_2$ and $q_3$. We prove the basic properties of the algebra $\calA$, including the commutativity and the Poincar\'e series, by introducing a certain filtration associated to the dominance ordering among partitions. This stratification will be called the {\it Gordon filtration}, and it  is characterized by a sequence of null conditions associated to the partitions and to the shift parameter $q_i$. Our algebra can be regarded as a smooth limit of the elliptic algebra \cite{FO:1997}. Precisely speaking, the original algebra is constructed over an elliptic curve and our algebra $\calA$ is constructed over a degenerate  $\bbC\bbP^1$, i.e. a rational curve obtained from an elliptic curve by pinching one cycle.

\subsection{Commutative algebra with three shift parameters 
on degenerate  $\bbC {\bbP}^1$}

Let $\bbC {\bbP}^1=\{(\xi_1:\xi_2) \mid (\xi_1,\xi_2)\in\bbC^2\setminus \{(0,0)\}\}$ and $x=\xi_1/\xi_2$ and $1/x=\xi_2/\xi_1$ be the local coordinates. Let $x=0,\infty$ be the two marked points on $\bbC {\bbP}^1$, being fixed throughout the text. We denote the set of non-negative integers by $\bbN\seteq\{0,1,2,\ldots\}$ and the set of positive integers by $\bbN_+\seteq\{1,2,\ldots\}$.
Let $q_1,q_2$ be two independent indeterminates.  We set $q_3=1/q_1q_2$ 
throughout the paper. We denote by $\bbF\seteq \bbQ(q_1,q_2)$ the field of rational functions in $q_1,q_2$.

\begin{dfn}[Operators $\partial^{(0,k)}$ and $\partial^{(\infty,k)}$]
\label{dfn:bibun}
For $n,k\in\bbN_+$, we define two operators $\partial^{(0,k)},\partial^{(\infty,k)}$ acting on the space of symmetric rational functions in $n$ variables $x_1,\ldots,x_n$ by
\begin{align*}
\begin{array}{l l l l l l l l l}
\partial^{(0,k)}&: &f &\mapsto 
 &\displaystyle 
   \dfrac{n!}{(n-k)!} \lim_{\xi \to 0} 
   f(x_1,\ldots,x_{n-k},\xi x_{n-k+1},\xi x_{n-k+2},\ldots,\xi x_n)
\\
\partial^{(\infty,k)}&: &f &\mapsto 
 &\displaystyle 
    \dfrac{n!}{(n-k)!} \lim_{\xi \to \infty}
    f(x_1,\ldots,x_{n-k},\xi x_{n-k+1},\xi x_{n-k+2},\ldots,\xi x_n)
\end{array}
\end{align*}
whenever the limit exists. We also set $\partial^{(0,k)} c=0, \partial^{(\infty,k)} c=0$ for $c\in \bbF$. Finally we define $\partial^{(0,0)}$ and $\partial^{(\infty,0)}$ to be the identity operator.
\end{dfn}

\begin{dfn}[Space $\calA$]
\label{dfn:calA}
For $n\in\bbN$, the vector space $\calA_n=\calA_n(q_1,q_2,q_3)$ is defined by the following conditions (i), (ii), (iii) and (iv).
\\
(i)
$\calA_0 \seteq \bbF$. For $n\in\bbN_+$, $f(x_1,\ldots,x_n)\in \calA_n$ is a rational function with coefficients in $\bbF$, and symmetric with respect to the $x_i$'s.
\\
(ii) For $n\in\bbN$, $0\leq k\leq n$ and $f \in \calA_n$, 
the limits $\partial^{(\infty,k)}f$ and $\partial^{(0,k)}f$ both exist and coincide: $\partial^{(\infty,k)}f=\partial^{(0,k)}f$ ({\it degenerate $\bbC \bbP^1$ condition}).
\\
(iii)
The poles of $f\in \calA_n$ are located only on the diagonal $\{(x_1,\ldots,x_n) \mid \exists (i,j), i\neq j ,x_i=x_j\}$, and the orders of the poles are at most two.
\\
(iv) For $n\geq 3$, $f\in \calA_n$ satisfies the {\it wheel conditions}
\begin{align}
\label{eq:wheel_condition}
f( x_1,q_1 x_1,q_1 q_2 x_1,x_4,\ldots)=0,\qquad
f( x_1,q_2 x_1,q_1 q_2 x_1,x_4,\ldots)=0.
\end{align}

Then set the graded vector space 
$\calA=\calA(q_1,q_2,q_3)\seteq\bigoplus_{n\geq 0}\calA_n$.
\end{dfn}

\begin{rmk}
Using $q_3=1/q_1q_2$, we can rewrite the wheel conditions (\ref{eq:wheel_condition}) in a symmetric way as
\begin{align*}
f( q_1 x_1,q_1 q_2 x_1,q_1 q_2 q_3 x_1,x_4,\ldots)=0,\qquad 
f( q_3 x_1,q_3 q_2 x_1,q_3 q_2 q_1 x_1,x_4,\ldots)=0,
\end{align*}
which indicates that the space $\calA$ is symmetric under the interchange of the $q_i$'s.
\end{rmk}

We introduce a bilinear operation $*$ called the {\it star product}.

\begin{dfn}[Star product $*$]
\label{dfn:star_prod}
For an $m$-variable symmetric rational function $f$ and an $n$-variable symmetric 
rational function $g$, we define an $(m+n)$-variable symmetric rational function $f*g$ by 
\begin{align}\label{eq:star}
 (f*g)(x_1,\ldots,x_{m+n})&=
  \Sym
  \bigg[
   f(x_1,\ldots,x_m) g(x_{m+1},\ldots,x_{m+n})
   \prod_{\substack{1\le\alpha\le m\\m+1\le\beta\le m+n}}
   \omega(x_\alpha,x_\beta)
  \bigg].
\end{align}
Here $\omega(x,y)$ is the rational function
\begin{align}\label{eq:omega}
\omega(x,y)=\omega(x,y;q_1,q_2,q_3)=\dfrac{(x-q_1 y)(x - q_2 y)(x-q_3 y)}{(x-y)^3},
\end{align}
and the symbol $\Sym$ denotes the symmetrizer
\begin{align*}
\Sym (f(x_1,\ldots,x_n))
\seteq\dfrac{1}{n!}\sum_{\sigma\in\frakS_{n}}\sigma(f(x_1,\ldots,x_n)),
\end{align*}
where $\frakS_n$ is the $n$-th symmetric group acting on the indices of $x_i$'s.
\end{dfn}

A priori, it is unclear in which space the star product is closing. Neither is obvious whether the resulting algebra is commutative. Moreover the study of the Poincar\'e series usually requires some technique. Nevertheless we may solve all of these problems in a combinatorial manner.

\begin{thm}\label{thm:1}
The vector space $\calA$ is closed with respect to the 
star product $*$, hence the pair $(\calA,*)$ defines a unital associative algebra. The algebra $(\calA,*)$ is commutative. The Poincar\'e series is $\sum_{n\ge 0} (\dim_\bbF \calA_n) z^n=\prod_{m\ge 1}(1-z^m)^{-1}$.
\end{thm}

A proof of Theorem \ref{thm:1} will be given in \S \ref{sect:basics}.

\subsection{Gordon filtration and intersection}
\label{subsec:gordon}

Next, we introduce a stratification in the space $\calA$ by using a set of specialization maps.

We recall some basic definitions concerning partitions. We basically follow the notation in \cite{M:1995:book}. A partition of $n\in \bbN$ is a sequence  $\lambda=(\lambda_1,\lambda_2,\ldots)$ of non-negative integers satisfying $\lambda_1\geq\lambda_2\geq\cdots$. The weight and the length of $\lambda$ are defined by $|\lambda|\seteq\lambda_1+\lambda_2+\cdots$ and $\ell(\lambda)\seteq\#\{i\mid\lambda_i\neq0\}$. We write $\lambda \vdash n$ if $\lambda$ is a partition of $n$. As usual, we denote the conjugate of a partition $\lambda$ by $\lambda'$, which is corresponding to the transpose of the diagram of  $\lambda$. We work with the dominance partial ordering defined as follows:
\begin{align*}
\lambda\geq \mu \overset{\rm def}{\iff}
|\lambda|=|\mu|,\  \lambda_1+\cdots+\lambda_i\geq \mu_1+\cdots+\mu_i
\mbox{ for all } i\geq 1.
\end{align*}

\begin{dfn}[Specialization map $\varphi$]\label{dfn:specialization_map}
Let $p\in \bbF=\bbQ(q_1,q_2)$. For a partition $\lambda=(\lambda_1,\ldots,\lambda_m)$ of $n$, we define a linear map
\begin{align}
\label{eq:specialization_map}
\begin{array}{l c c c}
\varphi_\lambda^{(p)} : 
&\calA_n          &\longrightarrow&\bbF(y_1,\ldots,y_m)\\
&f(x_1,\ldots,x_n)&\mapsto        &f(y_1,p y_1,\ldots,p^{\lambda_1-1}y_1,
y_2,p y_2,\ldots,p^{\lambda_2-1}y_2,\ldots\\
&                 &               &\hskip 1.5em \ldots,y_m,p y_m\ldots,p^{\lambda_m-1}y_m).
\end{array}
\end{align}
\end{dfn}

\begin{dfn}[Gordon filtration]
\label{dfn:gordon_filtr}
For $q_i$ ($i=1,2,3$) and $\lambda\vdash n$, define $\calA_{n,\lambda}^{(q_i)}\subset \calA_n$ by
\begin{align}
\label{eq:gordon_filtr}
\calA_{n,\lambda}^{(q_i)}\seteq\bigcap_{\mu\not\le\lambda}\ker\varphi_\mu^{(q_i)}\quad (\lambda<(n)),\qquad
\calA_{n,(n)}^{(q_i)}\seteq \calA_n. 
\end{align}
\end{dfn}

The next inclusion property of the strata $\calA_{n,\lambda}^{(q_i)}$ easily follows from the definition.

\begin{lem}
If $\lambda\ge\mu$, then 
$\calA_{n,\lambda}^{(q_i)}\cap\ker\varphi_\lambda^{(q_i)}
\supseteq \calA_{n,\mu}^{(q_i)}$.
In particular $\calA_{n,\lambda}^{(q_i)}\supseteq \calA_{n,\mu}^{(q_i)}$.
\end{lem}

\begin{rmk}
Hence we have a filtration of the space $\calA_n$ whose inclusion sequence corresponds to the dominance ordering of the partitions of $n$.
\end{rmk}

By using only one of the shift parameters, say $q_1$, a combinatorial study of the Gordon filters $\calA_{n,\lambda}^{(q_1)}$ works well enough to prove Theorem \ref{thm:1}. Moreover, one finds a quite interesting interplay between the two different shift parameters, say $q_1$ and $q_2$.

\begin{thm}\label{thm:2}
Let $n\in\bbN$ and $\lambda$ be a partition of $n$. 
Then the intersection of the subspaces $\calA_{n,\lambda}^{(q_1)}$ and $\calA_{n,{\lambda'}}^{(q_2)}$ is 1-dimensional:
\begin{align*}
\dim_\bbF
\bigl(\calA_{n,\lambda}^{(q_1)}\cap\calA_{n,{\lambda'}}^{(q_2)}\bigr)=1.
\end{align*}
\end{thm}

A proof of Theorem \ref{thm:2} will be given in \S \ref{sect:mcd}. Our proof essentially uses the free field construction of the Macdonald difference operators \cite{AMOS:1995,S:2003:book,S:2006}.

\begin{rmk}
It is an interesting open problem to find a proof of Theorem \ref{thm:2} purely inside the algebra $\calA$ without referring to the free field construction.
\end{rmk}

\subsection{Heisenberg representation of the Macdonald difference operators}
\label{subsec:intro_ffr}
We recall the Heisenberg representation of the Macdonald difference operators,
which is initiated by Jing \cite{J:1994}. We basically follow the notation in \cite{S:2006}. Our aim is to introduce a canonical basis for the space of Macdonald difference operators in the sense of Theorem \ref{thm:3} below.

Let $q$ and $t$ be independent indeterminates, which are the
two parameters for the Macdonald theory in the standard notation.
We identify the parameters as 
$q_1=q^{-1}$ and $q_2=t$ (the condition $q_1q_2q_3=1$ means $q_3=t^{-1}q$),
and denotes our base field by the same symbol $\bbF=\bbQ(q,t)=\bbQ(q_1,q_2)$.

Consider the Heisenberg Lie algebra $\frakh$ over $\bbF$ with the generators $a_n$ ($n\in\bbZ$) and the relations
\begin{align}\label{eq:boson_macdonald}
[a_m,a_n] = m\dfrac{1-q^{|m|}}{1-t^{|m|}}\delta_{m+n,0}\,a_0.
\end{align}
Let $\frakh^{\ge0}$ (resp. $\frakh^{<0}$) be the subalgebra generated by $a_n$ for $n\ge 0$ (resp. $n< 0$) . Consider the Fock representation $\calF\seteq\Ind_{\frakh^{\ge0}}^\frakh\bbF$ of $\frakh$. We also use the notation $a_{-\lambda}\seteq a_{-\lambda_1} \cdots a_{-\lambda_l}$ for a partition $\lambda=(\lambda_1,\ldots,\lambda_l)$. 

Let $x=(x_1,x_2,\ldots)$ be a set of indeterminates and $\Lambda$ be the ring of symmetric functions in $x$ over $\bbZ$. As a $\bbF$-vector spaces, $\calF$ is isomorphic to $\Lambda_\bbF\seteq\Lambda\otimes_\bbZ\bbF$ via $
\iota:\calF\to\Lambda_\bbF$ defined by $a_{-\lambda}\cdot 1 \mapsto p_\lambda$. Here  $1\in\calF$ is the highest vector, $p_n=p_n(x)\seteq \sum_{i\geq 1}x_i^n$ is the power sum function and $p_\lambda\seteq p_{\lambda_1} \cdots p_{\lambda_l}$. 
For $n>0$ and $v\in \calF$, we have $a_{-n} v=p_n v$,
$a_{n} v=n{((1-q^n)/( 1-t^n))}\partial v/\partial p_n$, and $a_0v=v$.
In what follows we identify $\calF\simeq \Lambda_\bbF$ via this isomorphism $\iota$.

Introduce the vertex operator
\begin{align}
\begin{split}
\label{eq:vertex_operator}
\eta(z)
&=\,
\exp\Bigg( \sum_{n>0} \dfrac{1-t^{-n}}{n}a_{-n} z^{n} \Bigg)
\exp\Bigg(-\sum_{n>0} \dfrac{1-t^{n} }{n}a_n    z^{-n}\Bigg)\\
&=
:\,\exp\bigl(-\sum_{n\neq0}\dfrac{1-t^{n}}{n}a_{n}z^{-n}\bigr)\,:,
\end{split}
\end{align}
where the symbol $:\,\, :$ denotes the normal ordering with respect to the decomposition $\frakh=\frakh^{<0} \oplus \frakh^{\ge0}$, {\it i.e.} all the negative  generators $a_{-n}$ are moved to the left of the positive generators $a_n$. Then the Fourier modes $\eta_n$ of $\eta(z)$, defined by $\eta(z)=\sum_{n\in \bbZ}\eta_n z^{-n}$, are well-defined operators acting on the Fock space $\calF$ or on $\Lambda_\bbF$. The operator $E$ of Macdonald \cite[(VI.4.2)]{M:1995:book} is realized as $\eta_0=(t-1)E+1$. We note that the plethystic operator $\Delta'$ in \cite[(2.10)]{H:1999} and $\eta_0$ are identical.

Recall the existence theorem of the Macdonald symmetric functions:
\begin{prop}[Macdonald]\label{prop:Macd}
The Macdonald symmetric functions $P_\lambda(x;q,t)$ are uniquely determined by the conditions:
\begin{align*}
\tag{a} 
&P_\lambda(x;q,t)=m_\lambda+\sum_{\mu<\lambda} 
 c_{\lambda\mu}^{m\rightarrow P}m_\mu
 \qquad (c_{\lambda\mu}^{m\rightarrow P}\in \bbF),\\
\tag{b} 
&E P_\lambda(x;q,t)=\sum_{i\geq 1} (q^{\lambda_i}-1)t^{-i} \cdot  P_\lambda(x;q,t).
\end{align*}
Here $m_\lambda$ denotes the monomial symmetric function.
\end{prop}

\begin{rmk}\label{rmk:spectral}
{} From $\eta_0=(t-1)E+1$, we may write (b) as
\begin{align*}
& \dfrac{t^{-1}}{1-t^{-1}}\eta_0  
P_\lambda(x;q,t)=e_1(s^\lambda)  P_\lambda(x;q,t),
\end{align*}
where $e_1(x)=x_1+x_2+\cdots$ is the first elementary symmetric function, and $s^\lambda=(s_1^{\lambda},s_2^{\lambda},\ldots)=(t^{-1}q^{\lambda_1},t^{-2}q^{\lambda_2},\ldots)$. We note that the first factor in the LHS should be understood as the power series $t^{-1}/(1-t^{-1})=t^{-1}+t^{-2}+\cdots$, corresponding to $e_1(s^\lambda)\in \bbF[[t^{-1}]]$.
\end{rmk}

A detailed study of the vertex operator $\eta(z)$ 
derives the mutually commuting family of operators, 
which we call $\calM$, including $E$ as the first member.
We begin with the next calculation.

\begin{lem}
It holds that
\begin{align}\label{eq:eta_ope}
\eta(z)\eta(w)=\dfrac{(1-w/z)(1-q t^{-1}w/z)}{(1-q w/z)(1-t^{-1}w/z)}
 :\eta(z)\eta(w):.
\end{align}
Here and hereafter throughout the text 
the rational factors such as $1/(1-q w/z)$ and $1/(1-t^{-1} w/z)$
must be understood as power series in $w/z$.
\end{lem}

This gives the operator valued {\it symmetric Laurent series} 
\begin{align}\label{eq:eta_ep}
\begin{split}
&\dfrac{1}{\omega(z_1,z_2;q^{-1},t,q t^{-1})}
\eta(z_1)\eta(z_2)
=\dfrac{1}{\omega(z_2,z_1;q^{-1},t,q t^{-1})}
\eta(z_2)\eta(z_1)\\
=&\dfrac{q t^{-1}(1-z_2/z_1)^2 }
{(1-q z_2/z_1)(1-t^{-1}z_2/z_1)
}\dfrac{ (1-z_1/z_2)^2}
{(1-q z_1/z_2)(1-t^{-1}z_1/z_2)
}
:\eta(z_2)\eta(z_1):.
\end{split}
\end{align}

For $f\in\calA_n(q^{-1},t,q t^{-1})$, consider the operator product 
\begin{align*}
&f(z_1,\ldots,z_n) \prod_{1\leq i<j\leq n} \omega(z_i,z_j)^{-1}\cdot 
\eta(z_1)\cdots \eta(z_n)\\
&=
\biggl(f(z_1,\ldots,z_n)
\prod_{1\leq i<j\leq n} \dfrac{(z_i-z_j)(z_j-z_i)}{z_i z_j}
\biggr)\\
&\times
\biggl(\biggl[
 \prod_{1\leq i<j\leq n}
 \dfrac{q t^{-1}(1-z_j/z_i) }{(1-q z_j/z_i)(1-t^{-1}z_j/z_i)}
 \dfrac{ (1-z_i/z_j)}{(1-q z_i/z_j)(1-t^{-1}z_i/z_j)}
\biggr]:\eta(z_1)\cdots \eta(z_n):\biggr).
\end{align*}
Note that the first term in the RHS is a symmetric Laurent polynomial in $z_i$'s, hence the RHS in total is a symmetric Laurent series.

For a Laurent series in $n$ variables $f(z_1,\ldots,z_n)$ we denote by $[f(z_1,\ldots,z_n)]_1$ the constant term.

\begin{dfn}
\label{dfn:calO}
Let $f\in \calA_n(q^{-1},t,q t^{-1})$. Define a mapping 
$\calO(\cdot)=\calO(\cdot\, ;q,t):\calA\rightarrow {\rm End}_\bbF(\calF)$ 
by 
\begin{align*}
\calO(f)=\calO(f;q,t)\seteq\biggl[
 \dfrac{f(z_1,\ldots,z_n)}
       {\prod_{1\le i<j\le n}\omega(z_i,z_j;q^{-1},t,q t^{-1})}
 \eta(z_1)\cdots\eta(z_n)\biggr]_1,
\end{align*}
for $f\in \calA_n$ and extending it by linearity.
\end{dfn}

The star product $*$ and the operation $\calO(\cdot)$ are compatible in the following sense.

\begin{prop}\label{prop:calO}
For $f,g\in \calA$, we have
\begin{align*}
\calO(f*g)=\calO(f)\calO(g),
\end{align*}
which indicates that $[\calO(f),\calO(g)]=0$ from the commutativity of $\calA$.
\end{prop}

\begin{dfn}\label{dfn:calM}
Define the commutative ring $\calM$ of operators on $\Lambda_\bbF$ by $\calM\seteq\{\calO(f) \mid f\in \calA\}$.
\end{dfn}

\begin{prop}\label{prop:spectral_criterion}
The mapping $\calO(\cdot):\calA\to\calM$ gives an isomorphism of commutative rings.
\end{prop}

The proofs of these Propositions \ref{prop:calO} and \ref{prop:spectral_criterion} will be given in \S\S\,\ref{subsec:prf_prop_calO_spec}.

We are ready to state our interpretation of the intersection of the two Gordon filters given in Theorem \ref{thm:2}, in terms of the Macdonald symmetric functions.

\begin{thm}\label{thm:3}
Let $n\in \bbN$ and $\mu\vdash n$.  There exists a unique element $f_\mu \in \calA_n$ such that 
\begin{align*}
\calO(f_\mu)P_\lambda(x;q,t)=P_\mu(s^\lambda;q,t)P_\lambda(x;q,t)
\end{align*}
for any partition $\lambda$, where $s^\lambda$ is the same as in Remark \ref{rmk:spectral}. Then the intersection of the subspaces $\calA_{n,\mu}^{(q^{-1})}$ and $\calA_{n,{\mu'}}^{(t)}$ (which is one dimensional from Theorem \ref{thm:2}) is spanned by $f_\mu$.
\end{thm}

Theorem \ref{thm:3} will be proved in \S\S\,\ref{subsec:intersect}.

It is worthwhile to make a comment here. There is a direct connection between the Gordon filtrations and some particular bases of $\Lambda_\bbF$. Recall the elementary symmetric function $e_n(x)=P_{(1^n)}(x;q,t)$ and the symmetric function $g_n(x;q,t)\seteq\prod_{i=1}^n (1-t q^{i-1})/(1-q^i) \cdot P_{(n)}(x;q,t)$. Set $e_\lambda\seteq e_{\lambda_1}e_{\lambda_2}\cdots$ and $g_\lambda\seteq g_{\lambda_1}g_{\lambda_2}\cdots$. It is known that $(e_\lambda(x))$ and $(g_\lambda(x;q,t))$
are bases of $\Lambda_\bbF$.

\begin{prop}[{Haiman \cite[Proposition 2.6]{H:1999}}]\label{prop:Haiman}
Let $\lambda\vdash n$. Let $V_\lambda$ be the totality of the symmetric functions which 
can be expanded in two ways as
\begin{align*}
\sum_{\mu\ge\lambda}c_{\mu}g_\mu(x;q,t)
=\sum_{\nu\ge\lambda'}c'_{\nu}e_\nu(x)\qquad 
(c_{\mu},c'_{\nu}\in \bbF).
\end{align*}
Then $V_\lambda$ is one dimensional, and is spanned by $P_\lambda(x;q,t)$.
\end{prop}

It is clear that our Theorems \ref{thm:2} and \ref{thm:3} resemble Haiman's Proposition \ref{prop:Haiman}.

\subsection{Ding-Iohara algebra $\calU(q,t)$}
Recall that the Ding-Iohara algebra \cite{DI:1997} was introduced as a generalization of the quantum affine algebra, which respects the structure of the Drinfeld coproduct. 
We consider a particular case of the Ding-Iohara algebra having two parameters $q$ and $t$, which we call $\calU(q,t)$. 
As for the definition of $\calU(q,t)$, see Appendix \ref{app:DI}.
In the course of the study on the vertex operator $\eta(z)$, 
the authors were lead to a particular representation of 
$\calU(q,t)$. At present, we do not have a good understanding 
of the connection between the Macdonald theory and $\calU(q,t)$. 
Thus in this paper the authors decided to treat 
the connection as an unexpected byproduct (see \S\S\,\ref{subsec:DI}),
instead of starting from the Ding-Iohara algebra $\calU(q,t)$ and
deriving all the basic properties which we need for our free field representation.
However, 
the authors hope the present paper can also be 
read in the logical direction. Namely,
read Appendix \ref{app:DI} first where 
we collected the materials related to $\calU(q,t)$, and proceed to 
\S\S \ref{subsec:DI}. After that,
concerning the free field approach, every formulas can be 
understood from the viewpoint of the representation theories of 
the Ding-Iohara algebra.

As a related work, we mention the following.
In the paper \cite{FT}, one of the authors and 
Tsymbaliuk studied the equivariant $K$-theory by using 
the Ding-Iohara and shuffle algebras.

Our heuristic argument is explained as follows. 
We have introduced the vertex operator $\eta(z)$ acting on $\calF$. 
One can construct another one, which we call $\xi(z)=\sum_n\xi_n z^{-n}$.
See  (\ref{eq:xi}) for its definition. 
Then we find the following three.\\
(1) The permutation relation between $\xi(z)$ and $\xi(w)$ 
can be obtained by replacing $(\eta,q,t)$ to $(\xi,q^{-1},t^{-1}$), meaning that 
one can construct a family of commuting operators $\calM'$ in terms of $\calA$ and $\xi(z)$
with a suitable change of the parameters.\\
(2) Then we have the commutativity $[\eta_0,\xi_0]=0$. This indicates that
the operators $f\in\calM$ and $g\in M'$ commutes with each other.
See Proposition \ref{prop:M_M'}.\\
(3) When we prove $[\eta_0,\xi_0]=0$, 
two more vertex operators $\varphi^+(z)$ and $\varphi^-(z)$ inevitably come into the game.
See (\ref{eq:phi^+}) and (\ref{eq:phi^-}).

Hence we have four vertex operators $\eta(z),\xi(z),\varphi^\pm(z)$ in total. 
A simple calculation gives us a proof that these vertex operators 
give us a Fock representation of $\calU(q,t)$, which is 
our main statement concerning the Ding-Iohara algebra. 
See Proposition \ref{prop:DI_Fock_rep}.

\subsection{Elliptic deformation}

Our algebra $\calA(q_1,q_2,q_3)$ can be regarded as 
a degenerate limit of the elliptic algebra proposed in \cite{FO:1997},
which we denote by $\calA(p)=\calA(q_1,q_2,q_3,p)$. 
The basic properties of the elliptic algebra $\calA(p)$ 
are stated in Proposition \ref{prop:calA(p)}, 
showing that the deformation is {\it smooth} 
and the structure of $\calA$ remains the same while we deform it.

It has been recognized  that tensor representations of 
the so-called dynamical quantum groups, including 
elliptic algebras,
take their natural place when the Drinfeld quasi-Hopf twisting is applied \cite{QHA}. 
See \cite{Fel95,EF,Fron,Fron1,ABRR, JKOS1,JKOS2} for example and references therein.
We aim at 
making the representation theory of the Ding-Iohara algebra $\calU(q,t)$ 
compatible with the elliptic algebra $\calA(p)$, and apply the 
quasi-Hopf twisting prescribed by 
Babelon, Bernard and Billey \cite{BBB}. We call the resulting 
quasi-bialgebra  $\calU(q,t,p)$ for short. Namely, we
twist $\Delta$ (see Proposition \ref{prop:Hopf-alg})
into a $p$-depending one $\Delta_{p}$ (see Proposition \ref{prop:QBA}).
Then the dressed Drinfeld generators $\eta(z,p), \xi(z,p), 
\varphi^\pm(z,p)$, 
as in Definition \ref{dfn:dressed_currents}, 
obey the dressed Drinfeld coproduct
in Proposition \ref{prop:dressed_copro}.

Our main result, as regarding the elliptic deformation $\calM(p)$ of the 
commuting family $\calM$, is stated in 
Theorem \ref{thm:calM(p)}. 

Finally, 
we remark that 
the Ruijsenaars difference operator \cite{Ruijsenaars:1987}
and the operator ${\sf M}(q,t_1,t_2)$ of Okounkov and Pandharipande \cite{OkounkovPandharipande:2004} are appearing 
in our framework based on the quasi-Hopf algebra 
$\calU(q,t,p)$. 
As for the Ruijsenaars difference operator,
see Corollary \ref{cor:Ruij}.
Take the zero Fourier component of the dressed current $\eta(z,p q^{-1}t)$
and consider
the limit $q=e^\hbar$, $t=e^{\beta \hbar}$, $\hbar\rightarrow 0$,
while $p$ being fixed. 
Then one can recover ${\sf M}(q,t_1,t_2)$. 
See Proposition \ref{prop:O-P}.
We also make a comment that 
$\calU(q,t,p)$ provides us with a natural framework 
to understand the deformed $\calW_m$ algebra 
introduced in \cite{SKAO:1996,FF,AKOS:1996}
from the point of view of quasi-Hopf algebra. See Proposition \ref{prop:W_m}.

\subsection{Plan of the paper}
This paper is organized as follows. In \S\,\ref{sect:basics} we present our proof of Theorem \ref{thm:1}. \S\,\ref{sect:mcd} is devoted 
to the proof of Theorems \ref{thm:2} and \ref{thm:3}.
Since our proof uses the technique of the free field construction, 
materials related to the algebra $\calA(q^{-1},t,q t^{-1})$ and to the Ding-Iohara algebra $\calU(q,t)$ are dealt in a complementary manner. 
In \S\,\ref{sect:elliptic}, 
based on the elliptic algebras $\calA(q_1,q_2,q_3,p)$ and $\calU(q,t,p q^{-1}t)$, we consider the elliptic deformation $\calM(p)$ of the 
family of Macdonald difference operators $\calM$. 
Definition and basic properties of the Ding-Iohara algebra $\calU(q,t)$ are
summarized in Appendix \ref{app:DI}.

\tableofcontents

\section{Proof of Theorem \ref{thm:1}}
\label{sect:basics}
This section is devoted to the proof of Theorem \ref{thm:1}. We embed $\calA$ in an ambient space which is called $\wtcalA$. Then we gradually reduce our space  as $\wtcalA\supset\whcalA\supset\wbcalA\supset\calA$. 
The same statements listed in Theorem \ref{thm:1} 
will be proved for $\wbcalA$. After that we will show $\calA\supset \wbcalA$.

\subsection{Construction of the Algebra $\wbcalA$}
\label{subsec:wtcal}

First we recall the notation of \cite{M:1995:book} for the symmetric polynomials. Let $x_1,x_2,\ldots,x_n$ be a set of indeterminates and $\Lambda_n\seteq\bbZ[x_1,\ldots,x_n]^{\frakS_n}$ be the ring of symmetric polynomials. $\Lambda_n$ has a natural grading, $\Lambda_n=\bigoplus_{k\ge0}\Lambda_n^k$, where $\Lambda_n^k$ consists of the homogeneous symmetric polynomials of degree $k$, together with the zero polynomial. 
We denote the base changes by $\Lambda_{n,\bbF}\seteq\Lambda_n \otimes_\bbZ \bbF$ and $\Lambda_{n,\bbF}^k\seteq\Lambda_n^k\otimes_\bbZ \bbF$. We 
write the van der Monde determinant as 
$\Delta_n(x)=\Delta(x_1,\ldots,x_n)\seteq\prod_{1\le i<j\le n}(x_i-x_j)$.

\subsubsection{}

We start our consideration from the ambient space $\wtcalA$.

\begin{dfn}[Space $\wtcalA$]
For each $n\in\bbN$, we set
\begin{align*}
{\wtcalA}_n\seteq
\biggl\{\dfrac{p(x_1,\ldots,x_n)}{\Delta(x_1,\ldots,x_n)^2}
\ \bigg|\ p(x_1,\ldots,x_n)\in \Lambda_{n,\bbF}^{n(n-1)}\biggr\}\quad
(n\ge 1),\quad
\wtcalA_0\seteq \bbF.
\end{align*}
Define $\wtcalA \seteq \bigoplus_{n\ge0} \wtcalA_n$.
\end{dfn}

\begin{lem}
Let $m,n\in \bbN$. For $f \in \wtcalA_m$ and $g \in \wtcalA_n$ we have $f * g \in \wtcalA_{m+n}$.
\end{lem}

\begin{proof}
Set $f\seteq p(x_1,\ldots,x_m)/\Delta_m(x)^2$ and $g\seteq q(x_1,\ldots,x_n)/\Delta_n(x)^2$. Then we have
\begin{align*}
f * g
&=\dfrac{1}{(m+n)!}\sum_{\sigma\in\frakS_{m+n}}\sigma
\biggl(\dfrac{p(x_1,\ldots,x_m)}{\Delta_m(x)^2}
\dfrac{q(x_{m+1},\ldots,x_{m+n})}{\Delta(x_{m+1},\ldots,x_{m+n})^2}
\prod_{\substack{1\le\alpha\le m \\ m+1\le\beta\le m+n}}
\omega(x_\alpha,x_\beta)\biggl)\\
&=\dfrac{1}{\Delta_{m+n}(x)^3}
\sum_{\sigma\in\frakS_{m+n}}\sgn(\sigma)\sigma(\mbox{polynomial})
=\dfrac{1}{\Delta_{m+n}(x)^2}
\cdot\mbox{(symmetric polynomial)}.
\end{align*}
\end{proof}

\begin{lem}\label{lem:associative}
The star product $*$ is associative on $\wtcalA$.
\end{lem}

\begin{proof}
Take $f \in \wtcalA_l$, $g \in \wtcalA_m$ and $h \in \wtcalA_n$. Then
\begin{align*}
&(l+m+n)!\, (f*g)*h\\
&=\sum_{\sigma \in \frakS_{l+m+n}} \sigma
  \biggl(
   (f*g) \, h 
   \prod_{\substack{\ \ \ \ 1\le\alpha\le l+m\\l+m+1\le\beta\le l+m+n}}
    \omega(x_\alpha,x_\beta)
  \biggr)\\
&=\sum_{\sigma\in\frakS_{l+m+n}} \sigma
   \biggl[
    \dfrac{1}{(l+m)!}
    \sum_{\tau\in\frakS_{l+m}} \tau
     \biggl(
      f \, g 
      \prod_{\substack{1\le\gamma\le l\\l+1\le\delta\le l+m}}
       \omega(x_\gamma,x_\delta)
     \biggr)\, h
     \prod_{\substack{\ \ \ \ 1\le\alpha\le l+m\\l+m+1\le\beta\le l+m+n}}
      \omega(x_\alpha,x_\beta)
   \biggl]\\
&=\sum_{\sigma\in\frakS_{l+m+n}}
   \sigma
   \biggl[
    f \, g \, h \,
    \prod_{\substack{1\le\gamma\le l\\l+1\le\delta\le l+m}}
     \omega(x_\gamma,x_\delta) 
    \prod_{\substack{\ \ \ \ 1\le\alpha\le l+m\\l+m+1\le\beta\le l+m+n}}
     \omega(x_\alpha,x_\beta)
   \biggl].
\end{align*}
By the same argument we have
\begin{align*}
(l+m+n)!\, f*(g*h) =\sum_{\sigma\in\frakS_{l+m+n}} \sigma
   \biggl[
    f \, g \, h \,
    \prod_{\substack{\ l+1\le\gamma\le l+m\\l+m+1\le\delta\le l+m+n}}
     \omega(x_\gamma,x_\delta) \,
    \prod_{\substack{\ \ \ \ 1\le\alpha\le l\ \ \ \ \\ \ \ \ l+1\le\beta\le l+m+n}}
     \omega(x_\alpha,x_\beta)\biggl].
\end{align*}
\end{proof}

Hence, we have shown 
\begin{prop}\label{prop:tildeA}
$(\wtcalA,*)$ is a unital associative algebra.
\end{prop}

\subsubsection{}

Next, we consider the following subspace $\whcalA\subset \wtcalA$.

\begin{dfn}[Space $\whcalA$]
Set $\whcalA_0\seteq \bbF$. We define the subspace $\whcalA_n\subset \wtcalA_n$ by the condition: for any $f\in \whcalA_n$ and any $1\leq k\leq n$ the limit 
$\partial^{(\infty,k)} f$ (see Definition \ref{dfn:bibun}) exists. We set $\whcalA\seteq\bigoplus_{n\ge0}\whcalA_n$.
\end{dfn}

\begin{lem}\label{lem:hatA:multi}
If $f \in \whcalA_m$ and $g \in \whcalA_n$, then $f*g \in \whcalA_{m+n}$.
\end{lem}

\begin{proof}
Thanks to the definition of the star product (Definition \ref{dfn:star_prod}), the definition of $\omega(x,y)$ in (\ref{eq:omega}), and the definition of $\whcalA$, it is clearly seen that the required limits exist, hence $f*g\in \whcalA_{n+m}$.
\end{proof}

For our later purpose, we give an alternative characterization of the space $\whcalA$.
\begin{prop}\label{prop:hatA_2delta}
We have
\begin{align*}
\whcalA_n=
\bigg\{
 \dfrac{p(x_1,\ldots,x_n)}{\Delta_n(x)^{2}}\ \bigg|\ 
 p(x_1,\ldots,x_n)=
  \sum_{\lambda\le2\delta_n}
  c_\lambda m_\lambda(x_1,\ldots,x_n),\ c_\lambda\in\bbF
\bigg\},
\end{align*}
where $2\delta_n\seteq(2(n-1),2(n-2),\ldots,2,0)$ and  $m_\lambda$ denotes the monomial symmetric polynomial in $x_1,\ldots,x_n$ associated to the partition $\lambda$. 
\end{prop}

\begin{proof}
We show that for any $f\in\wtcalA_n$, the following conditions are equivalent:
\begin{align*}
f\in\whcalA_n
&\iff f=p(x_1,\ldots,x_n)/\Delta_n(x)^{2},\\
&\hskip 3em 
 \mbox{where }
 p\in \Lambda_{n,\bbF}^{n(n-1)}
 \mbox{\ is of the form\ }
 \sum_{\lambda\le2\delta_n}
  c_\lambda m_\lambda(x_1,\ldots,x_n),\ c_\lambda\in\bbF.
\end{align*}
Fix any $1\leq k\leq n$. Expand the numerator of $f\in\whcalA_n$, which is a symmetric polynomial in $x_{n-k+1},\ldots,x_n$, as
\begin{align*}
f(x_1,\ldots,x_n)=
 \dfrac{1}{\Delta_n(x)^2}
 \sum_{|\rho|\leq n(n-1)}
  h_\rho(x_1,\ldots,x_{n-k})m_\rho(x_{n-k+1},\ldots,x_{n}),
\end{align*}
where $h_\rho$ is a symmetric polynomial is $x_1,\ldots,x_{n-k}$. Shifting the variables as $x_{i}\rightarrow \xi x_{i}$ ($n-k+1\leq i\leq n$), we consider the limit
\begin{align*}
&\lim_{\xi\rightarrow \infty}
 f(x_1,\ldots,x_{n-k},\xi x_{n-k+1},\ldots,\xi x_n)\\
=
&\lim_{\xi\rightarrow \infty}
 \dfrac{1}{\Delta_{n-k}(x)^2 \xi^{k(2n-k-1)}+\mbox{lower.}}
 \sum_{|\rho|\leq n(n-1)}
 h_\rho(x_1,\ldots,x_{n-k})
 m_\rho(x_{n-k+1},\ldots,x_{n})\xi^{|\rho|},
\end{align*}
where the symbol `lower.' denotes the lower degree terms in $\xi$. To have this limit existing, we have the conditions $\sum_{|\rho|=j}h_\rho m_\rho=0$ for $j>k(2n-k-1)$. From the linear independence of $m_\lambda$ in $\Lambda_{k,\bbF[x_1,\ldots,x_{n-k}]}$, we have $h_\rho $ for all $\rho$ satisfying $|\rho|>k(2n-k-1)$.
The converse is trivial. 
\end{proof}

Summarizing the construction for $\whcalA$, we have  
\begin{prop}\label{prop:hatA}
$(\whcalA,*)$ is a unital associative algebra. 
\end{prop}

\subsubsection{}
We impose the wheel conditions (\ref{eq:wheel_condition}) on the space $\whcalA$.

\begin{dfn}[Space $\wbcalA$]\label{rmk:old_defn}
Set  $\wbcalA_0\seteq\bbF$, $\wbcalA_1\seteq\bbF$ and $\wbcalA_2\seteq\whcalA_2$. For $n\geq 3$, let $\wbcalA_n\subset \whcalA$ be the subspace specified by the wheel conditions
\begin{align*}
f(x_1,q_1x_1,q_1q_2x_1,x_4,\ldots,x_n)=0,\qquad 
f(x_1,q_2x_1,q_1q_2x_1,x_4,\ldots,x_n)=0.
\end{align*}
Set $\wbcalA\seteq \bigoplus_{n\geq 0}\wbcalA_n$.
\end{dfn}

\begin{prop}\label{prop:wbcalA_alg}
$(\wbcalA,*)$ is a unital associative algebra.
\end{prop}

\begin{proof}
It is enough to show that $f*g\in\wbcalA_{m+n}$ for $f\in\wbcalA_m$ and $g\in\wbcalA_n$. Recalling the star product $*$ in Definition \ref{dfn:star_prod}, one can write it as
\begin{align}
\label{eq:star_rewrite}
f*g=\dfrac{m!n!}{(m+n)!}\sum_{\substack{I\subseteq\{1,\ldots,m+n\}\\\#I=m}}
f(x_I)g(x_{I^c})\prod_{\alpha\in I,\, \beta\in I^c}\omega(x_\alpha,x_\beta),
\end{align}
where we have used the abbreviation $x_I=\{x_{\alpha}\mid\alpha\in I\}$ and $x_{I^c}=\{x_{\beta}\mid\beta\in I^c\}$. If $\{1,2,3\}\subseteq I$ or $\{1,2,3\}\subseteq I^c$, we have
\begin{align*}
f(x_I)g(x_{I^c})|_{x_2=q_1x_1,x_3=q_1q_2x_1}=0
\end{align*} 
{}from the wheel condition for $f$ and $g$. In the other cases we have
\begin{align*}
\prod_{\alpha\in I,\, \beta\in I^c}
 \omega(x_\alpha,x_\beta)|_{x_2=q_1x_1,x_3=q_1q_2x_1}=0
\end{align*}
by the definition of $\omega$.  The study of $(f*g)(x_1,q_2x_1,q_1q_2x_1,x_4,\ldots)$ goes exactly in the same way.
\end{proof}

\subsection{Gordon filtration}
\label{subsec:Gordon}
The goal of this subsection is to calculate the dimension of $\wbcalA_n$ (see Proposition \ref{prop:wbcalA_n:poincare} below). In this and subsequent subsections, we need to analyze zero conditions of $f\in \wbcalA_n$. 
Before we proceed, we recast the results obtained in the last subsection as follows.

\begin{dfn}
Let $\calP_n=\calP_n(q_1,q_2,q_3)$ be the vector space of polynomials $p(x)=p(x_1,\ldots,x_n)$ satisfying the following conditions.
\\
(i)
$p$ is a homogeneous symmetric polynomial of degree $n(n-1)$.
\\
(ii)
$p$ has an expansion 
$
p(x_1,\ldots,x_n)
=\sum_{ \lambda\le2\delta_n}
c_\lambda m_\lambda(x_1,\ldots,x_n),\ c_\lambda\in\bbF.
$
\\
(iii)
If $n\ge3$, then $p$ satisfies the wheel conditions
\begin{align*}
p(x_1,q_1x_1,q_1q_2x_1,x_4,\ldots,x_n)=0,\qquad
p(x_1,q_2x_1,q_1q_2x_1,x_4,\ldots,x_n)=0.
\end{align*}
\end{dfn}

Then we have $\wbcalA_n=\{ p(x_1,\ldots,x_n)\Delta_n(x)^{-2} \mid p\in \calP_n\}$. Applying the specialization map $\varphi$ (\ref{eq:specialization_map}), we define the Gordon filtration on the space $\wbcalA_n$.

\begin{dfn}[Gordon filtration]
\label{dfn:gordon_filtr_Ab_P}
For $q_i$ ($i=1,2,3$), and a partition $\lambda$ of $n$ satisfying $\lambda<(n)$, define $\wbcalA_{n,\lambda}^{(q_i)}\subset \wbcalA_n$ by $\wbcalA_{n,\lambda}^{(q_i)}\seteq\bigcap_{\mu\not\le\lambda}\ker\varphi_\mu^{(q_i)}$. We also set $\wbcalA_{n,(n)}^{(q_i)}\seteq \wbcalA_n$. We introduce a filtration on the space $\calP_n$ in the same manner, and denote the corresponding filter by $\calP_{n,\lambda}^{(q_i)}$. We have $\wbcalA_{n,\lambda}^{(q_i)}=\calP_{n,\lambda}^{(q_i)}\cdot \Delta_n(x)^{-2}$. 
\end{dfn}

\begin{rmk}
One may give a diagrammatic meaning to the specialization map $\varphi$: for a partition $\lambda=(\lambda_1,\ldots,\lambda_m)$, we assign a monomial to each box in the  Young diagram of $\lambda$ by the rule
\begin{align*}
\mbox{the box at the}\ i\mbox{-th row and}\ j\mbox{-th column}
\ \mapsto \  y_i q^{j-1}.
\end{align*}
In Figure \ref{fig:special}, we show an example (the case $\lambda=(4,4,2,1,1,1)$).
\begin{figure}[htbp]
\begin{center}
\unitlength 0.1in
\begin{picture}(15.0, 16.0)(2.0,-18.0)
\special{pn 8}%
\special{pa 200 200}%
\special{pa 500 200}%
\special{pa 500 500}%
\special{pa 200 500}%
\special{pa 200 200}%
\special{fp}%
\special{pn 8}%
\special{pa 200 500}%
\special{pa 500 500}%
\special{pa 500 800}%
\special{pa 200 800}%
\special{pa 200 500}%
\special{fp}%
\special{pn 8}%
\special{pa 200 800}%
\special{pa 500 800}%
\special{pa 500 1100}%
\special{pa 200 1100}%
\special{pa 200 800}%
\special{fp}%
\special{pn 8}%
\special{pa 200 1100}%
\special{pa 500 1100}%
\special{pa 500 1400}%
\special{pa 200 1400}%
\special{pa 200 1100}%
\special{fp}%
\special{pn 8}%
\special{pa 200 1400}%
\special{pa 500 1400}%
\special{pa 500 1700}%
\special{pa 200 1700}%
\special{pa 200 1400}%
\special{fp}%
\special{pn 8}%
\special{pa 200 1700}%
\special{pa 500 1700}%
\special{pa 500 2000}%
\special{pa 200 2000}%
\special{pa 200 1700}%
\special{fp}%
\special{pn 8}%
\special{pa 500 800}%
\special{pa 800 800}%
\special{pa 800 1100}%
\special{pa 500 1100}%
\special{pa 500 800}%
\special{fp}%
\special{pn 8}%
\special{pa 500 500}%
\special{pa 800 500}%
\special{pa 800 800}%
\special{pa 500 800}%
\special{pa 500 500}%
\special{fp}%
\special{pn 8}%
\special{pa 800 500}%
\special{pa 1100 500}%
\special{pa 1100 800}%
\special{pa 800 800}%
\special{pa 800 500}%
\special{fp}%
\special{pn 8}%
\special{pa 1100 500}%
\special{pa 1400 500}%
\special{pa 1400 800}%
\special{pa 1100 800}%
\special{pa 1100 500}%
\special{fp}%
\special{pn 8}%
\special{pa 500 200}%
\special{pa 800 200}%
\special{pa 800 500}%
\special{pa 500 500}%
\special{pa 500 200}%
\special{fp}%
\special{pn 8}%
\special{pa  800 200}%
\special{pa 1100 200}%
\special{pa 1100 500}%
\special{pa  800 500}%
\special{pa  800 200}%
\special{fp}%
\special{pn 8}%
\special{pa 1100 200}%
\special{pa 1400 200}%
\special{pa 1400 500}%
\special{pa 1100 500}%
\special{pa 1100 200}%
\special{fp}%
\put( 2.9,-19.2){\makebox(0,0)[lb]{$y_6$}}%
\put( 2.9,-16.2){\makebox(0,0)[lb]{$y_5$}}%
\put( 2.9,-13.2){\makebox(0,0)[lb]{$y_4$}}%
\put( 2.9,-10.2){\makebox(0,0)[lb]{$y_3$}}%
\put( 5.2,-10.2){\makebox(0,0)[lb]{$q_1y_3$}}%
\put( 2.9,- 7.2){\makebox(0,0)[lb]{$y_2$}}%
\put( 5.2, -7.2){\makebox(0,0)[lb]{$q_1y_2$}}%
\put( 8.2, -7.2){\makebox(0,0)[lb]{$q_1^2y_2$}}%
\put(11.2, -7.2){\makebox(0,0)[lb]{$q_1^3y_2$}}%
\put( 2.9, -4.2){\makebox(0,0)[lb]{$y_1$}}%
\put( 5.2, -4.2){\makebox(0,0)[lb]{$q_1 y_1$}}%
\put( 8.2, -4.2){\makebox(0,0)[lb]{$q_1^2y_1$}}%
\put(11.2, -4.2){\makebox(0,0)[lb]{$q_1^3y_1$}}%
\end{picture}%
\end{center}
\caption{Specialization map $\varphi$ for $\lambda=(4,4,2,1,1,1)$.}
\label{fig:special}
\end{figure}
\end{rmk}

\begin{lem}\label{lem:atmost1}
Let $\lambda$ and $\mu$ be partitions of $n$ such that $\lambda > \mu$. Then for $i=1,2,3$ we have
\begin{align*}
\dim_\bbF
\varphi_\lambda^{(q_i)}(\wbcalA_{n,\lambda}^{(q_i)}/\wbcalA_{n,\mu}^{(q_i)})
\le 1.
\end{align*}
\end{lem}

\begin{proof}
{}From the symmetry, it is enough to check the case $i=1$, hence we suppress the dependence. Write $\lambda=(\lambda_1,\ldots,\lambda_m)$, and take an element $p\in\calP_{n,\lambda}$ which gives us a non-zero element
\begin{align*}
(\varphi_\lambda p)(y)
=p(y_1,q_1y_1\ldots,q_1^{\lambda_1-1}y_1,\ldots,y_m,q_1y_m\ldots,q_1^{\lambda_m-1}y_m)\neq0.
\end{align*}
The total degree of $\varphi_\lambda p$ is $n(n-1)$, and the degree in each $y_i$ should satisfy
\begin{align*}
\deg_{y_i}(\varphi_\lambda p)(y) \le \sum_{k=1}^{\lambda_i}2(n-k)
= 2 n \lambda_i-(\lambda_i+1)\lambda_i,
\end{align*}
because of the expansion $p(x_1,\ldots,x_n)=\sum_{\nu \le 2\delta_n}c_\nu m_\nu(x_1,\ldots,x_n)$. Now the assumption $p\in\calP_{n,\lambda}=\bigcap_{\mu\not\le\lambda}\ker\varphi_\mu$ shows that $\varphi_\lambda p$ is zero if
\begin{align}\label{eq:zeros}
\begin{cases}
y_j=q_1^{\lambda_i-k}y_i&k=0,1,\ldots,\lambda_j-1,\\
y_j=q_1^{-k}y_i&k=0,1,\ldots,\lambda_j-1,\\
y_j=q_2q_1^{l-k}y_i&k=0,1,\ldots,\lambda_j-1,\ l=1,\ldots,\lambda_i-1,\\
y_j=q_2^{-1}q_1^{l-k}y_i&k=0,1,\ldots,\lambda_j-1,\ l=1,\ldots,\lambda_i-1,
\end{cases}
\end{align}
where $1\le i<j\le m$. Note that the list here counts the zeros with correct multiplicities. Hence $\varphi_\lambda p$ is divisible by
\begin{align*}
q_\lambda(y)\seteq\prod_{1\le i<j\le m}
\biggl\{
 \biggl[
  \prod_{k=0}^{\lambda_j-1}(y_j-q_1^{\lambda_i-k}y_i)(y_j-q_1^{-k}y_i)
 \biggl]\,
 \biggl[
 \prod_{k=0}^{\lambda_j-1}\prod_{l=1}^{\lambda_i-1}
  (y_j-q_2q_1^{l-k}y_i)(y_j-q_2^{-1}q_1^{l-k}y_i)
 \biggl]
\biggl\}.
\end{align*}
The total degree of $q_\lambda(y)$ is 
\begin{align*}
\deg q_\lambda(y)
=\sum_{1\le i<j\le m}(2\lambda_i+2\lambda_i(\lambda_j-1))
=\sum_{i<j}2\lambda_i\lambda_j=n^2-\sum_i\lambda_i^2
\end{align*}
and the degree in each $y_i$ is
\begin{align*}
\deg_{y_i}q_\lambda(y)=\sum_{j\neq i}2\lambda_i\lambda_j=2n\lambda_i-2\lambda_i^2.
\end{align*}
Therefore the quotient $r(y)\seteq(\varphi_\lambda p)(y)/q_\lambda(y)$ has the total degree 
\begin{align*}
\deg r(y)=\deg (\varphi_\lambda p)(y)-\deg q_\lambda(y)
=n(n-1)-(n^2-\sum_i\lambda_i^2)
=\sum_i\lambda_i(\lambda_i-1)
\end{align*}
and the degree in each variable $y_i$ satisfies the inequality
\begin{align*}
\deg_{y_i}r(y)
=\deg_{y_i}(\varphi_\lambda p)(y)-\deg_{y_i}q_\lambda(y)
\le[2n\lambda_i-\lambda_i(\lambda_i+1)]-[2n\lambda_i-2\lambda_i^2]
=\lambda_i(\lambda_i-1).
\end{align*}
This means that the $r(y)$ should be $\prod_{i=1}^m y_i^{\lambda_i(\lambda_i-1)}$ up to constant multiplication. 

Hence we have shown that the image $\varphi_\lambda \calP_{n,\lambda}$ is the one dimensional vector space spanned by 
\begin{align}
\label{eq:zeta}
\zeta_{\lambda}(y_1,\cdots,y_m)\seteq
\prod_{j=1}^m y_j^{\lambda_j(\lambda_j-1)}
q_\lambda(y).
\end{align}
\end{proof}

\begin{dfn}[Element $\ep_n$]\label{dfn:bottom}
For $n>1$ and $i=1,2,3$, define
\begin{align*}
\ep_n=\ep_n(x;q_i)=\ep_n(x_1,\ldots,x_n;q_i)\seteq
\prod_{1\le k<l\le n}\dfrac{(x_k-q_i x_l)(x_k-q_i^{-1}x_l)}{(x_k-x_l)^2}.
\end{align*}
We also set $\ep_1=1$. For a partition $\lambda=(\lambda_1,\ldots,\lambda_m)$ of $n$, we write $\ep_\lambda\seteq\ep_{\lambda_1}*\cdots*\ep_{\lambda_m}$ for simplicity.
\end{dfn}

\begin{prop}[Bottom of the Gordon filter]\label{prop:bottom}
$\wbcalA_{n,(1^n)}^{(q_i)}$ is a 1-dimensional subspace of $\wbcalA_n$ spanned by $\ep_n(x;q_i)$.
\end{prop}

\begin{proof}
We may assume $n\ge2$. Let $f\in \wtcalA_n$ and write it as $f=p \cdot \Delta_n^{-2}$. We have
\begin{align*}
f\in\ker\varphi^{(q_i)}_{(2,1^{n-2})}
&\iff p(y_1,q_i y_1,y_2,\ldots,y_{n-2})=0\\
&\iff (x_2-q_i x_1) \mid p\iff (x_2-q_i x_1)(x_2-q_i^{-1}x_1) \mid p\\
&\iff \prod_{1\le k<l\le n}(x_k-q_i x_l)(x_k-q_i^{-1}x_l) \mid p,
\end{align*}
because $p$ is a symmetric polynomial. Since the degree of $p$ is $n(n-1)$, $f$ must be a constant multiple of $\ep_n(x;q_i)$.
\end{proof}

Our next task is the determination of the filter $\wbcalA_{n,\nu}^{(q_i)}$ in which the $\ep_{\lambda'}(x;q_i)$ lies ($\lambda'$ is the transpose of $\lambda$). We need to study the image $\varphi_\mu^{(q_i)}\ep_{\lambda'}(x;q_i)$ in some systematic manner. For $\lambda'=(\lambda'_1,\ldots,\lambda'_l)$, we have 
\begin{align}
\label{eq:irrlevant}
\ep_{\lambda'}(x)=
\dfrac{1}{n!}\prod_{j=1}^l\lambda'_j!
\sum_{I} \ep_{\lambda',I}(x),
\quad
\ep_{\lambda',I}(x)\seteq
 \prod_{h=1}^l \ep_{\lambda'_h}(x_{I_h})
 \prod_{1\leq j<k\leq l}
 \prod_{\substack{\alpha \in I_j \\ \beta\in I_k}}\omega(x_\alpha,x_\beta),
\end{align}
where the running index $I=(I_1,\ldots,I_l)$ satisfies $I_1\sqcup \ldots\sqcup I_l=\{1,\ldots,n\}$, and we used the abbreviation $x_{I_j}\seteq\{x_\alpha\mid\alpha\in I_j\}$ ($1\leq j\leq l$). To state when $\varphi^{(q_i)}_\mu \ep_{\lambda',I}(x)$ vanishes, we introduce the following terminology.

\begin{dfn}[$\varphi_\mu^{(q_i)}$-irrelevant index]
Let $\lambda \vdash n$ and $I$ as above. Let $\mu=(\mu_1,\ldots,\mu_m)\vdash n$. We say that the index $I$ is {\it $\varphi_\mu^{(q_i)}$-irrelevant}, if there exists a pair $(h,\alpha,\beta)$ ($1\leq h \leq m$ and $0\leq \alpha<\beta\leq \mu_h-1$) for which we have $y_h q_i^\alpha \in \varphi_\mu^{(q_i)}(x_{I_k})$ and $y_h q_i^\beta \in \varphi_\mu^{(q_i)}(x_{I_j})$ for some $1\leq j<k\leq l$. We call an index $I$ {\it $\varphi_\mu^{(q_i)}$-relevant} if $I$ is {\it not} $\varphi_\mu^{(q_i)}$-irrelevant.
\end{dfn}

\begin{lem}
\label{lem:irrelevant}
If $I$ is $\varphi_\mu^{(q_i)}$-irrelevant, then we have
$\varphi_\mu^{(q_i)} \ep_{\lambda',I}(x)=0$.
\end{lem}

\begin{proof}
This immediately follows from the condition $\omega(q_ix,x)=0$.
\end{proof}

\begin{prop}\label{prop:atleast1}
$\ep_{\lambda'}(x;q_i)\in\wbcalA_{n,\lambda}^{(q_i)}$ and $\varphi_\lambda^{(q_i)}\ep_{\lambda'}(x;q_i)\neq 0 $.
\end{prop}

\begin{proof}
We set $i=1$ and suppress the dependence on $q_1$. We first show that 
$\varphi_\mu \ep_{\lambda'}=0$ for $\mu\not\le\lambda$. Denote $\mu=(\mu_1,\ldots,\mu_m)$ and $\lambda'=(\lambda_1',\ldots,\lambda_l')$. 
The set of running indices $I$ in (\ref{eq:irrlevant}) are classified 
into two;  (Case I) set of $\varphi_\mu$-irrelevant indices and 
(Case II) set of $\varphi_\mu$-relevant ones.
Lemma \ref{lem:irrelevant} means we can discard all the indices in Case I. 
Now we subdivide Case II 
$=\{I\mid \mbox{for any pair } (h,\alpha,\beta), \,
y_h q^\alpha \in \varphi_\mu(x_{I_j}) \mbox{ and } 
y_h q^\beta \in \varphi_\mu(x_{I_k})
\mbox{ means } j\leq k\}$.
\begin{align*}
\mbox{(Case IIa) }\qquad & \mbox{there exists a pair } (h,\alpha,\beta) \mbox{ such that }
y_h q^\alpha,y_h q^\beta \in \varphi_\mu(x_{I_j}) \mbox{ for some } j,\\
\mbox{(Case IIb) }\qquad & \mbox{for any pair } (h,\alpha,\beta),\, 
y_h q^\alpha \in \varphi_\mu(x_{I_j}) \mbox{ and } 
y_h q^\beta \in \varphi_\mu(x_{I_k})
\mbox{ means } j< k.
\end{align*}
In Case IIa we must have $\beta=\alpha+1$, since otherwise the index must be an element of Case I, which is a contradiction.
If $\beta=\alpha+1$, $\varphi_\mu \ep_{\lambda'}=0$ by the definition of 
$\epsilon_r$. Finally, one finds that Case IIb cannot occur. 
For Case IIb implies 
$\mu_1+\cdots+\mu_{h}\leq \lambda_1+\cdots+\lambda_{h}$
for all $h\geq 1$, which contradicts the assumption 
$\mu\not\le\lambda$.

Now we prove the second statement. Set $\lambda=(\lambda_1,\ldots,\lambda_\ell)$ ($\ell=\lambda'_1$). Then we find that there exists only one  $\varphi_\lambda$-relevant index and it is $I=(\{1,\lambda_1+1,\lambda_1+\lambda_2+1,\ldots,\lambda_1+\cdots+\lambda_{\ell-1}+1\}, \{2,\lambda_1+2,\lambda_1+\lambda_2+2,\ldots,\lambda_1+\cdots+\lambda_{\ell-1}+2\},\ldots,\{l,\ldots\})$. Then
\begin{align*}
(\varphi_\lambda \ep_{\lambda'})(y) =
\dfrac{1}{n!}\prod_{h=1}^l\lambda'_h!
 \prod_{i=1}^l \ep_{\lambda_i'}(y_1,\ldots,y_{\lambda_i'})
 \prod_{1\le j<k\le l}
 \prod_{\alpha=1}^{\lambda_j'}
 \prod_{\beta=1}^{\lambda_k'}\omega(y_\alpha q^{j-1},y_\beta q^{k-1}),
\end{align*}
Thus $(\varphi_\lambda \ep_{\lambda'})(y)\neq 0$ and the proof is completed.
\end{proof}

\begin{prop}[Dimension of $\wbcalA_n$]\label{prop:wbcalA_n:poincare}
For a partition $\lambda\vdash n$, $(\ep_\mu(x;q_i))_{\mu\ge\lambda'}$ form a basis of $\wbcalA_{n,\lambda}^{(q_i)}$. In particular, $\wbcalA_n$ has a basis $(\ep_\lambda(x;q_i))_{\lambda\,\vdash\, n}$ and $\dim_\bbF\wbcalA_n$ equals to the number of partitions of $n$.
\end{prop}
\begin{proof}
Lemma \ref{lem:atmost1} and Proposition \ref{prop:atleast1} gives us
$\dim_\bbF
\varphi_\lambda^{(q_i)}(\wbcalA_{\lambda}^{(q_i)}/\wbcalA_{\mu}^{(q_i)})
= 1.$
\end{proof}

\subsection{Commutativity of $\wbcalA$}
\label{subsec:comm}

\begin{prop}[Commutativity of $\wbcalA$]\label{prop:commute}
$\wbcalA$ is a commutative algebra.
\end{prop}
\begin{proof}
It is enough to show that $\{\ep_n(x;q_i)\mid n\in\bbZ\}$ is a commutative family. We can suppose $i=1$ without loss of generality. 

First we will show that $\ep_n*\ep_1=\ep_1*\ep_n$. Consider the specialization 
\begin{align*}
& \varphi_{(2,1^{n-1})}^{(q_1)}[\ep_n,\ep_1]\\
&=\dfrac{1}{n+1}\biggl(\ep_n(y_1,\ldots,y_n)\omega(y_1,q_1y_1)
\prod_{i=2}^{n}\omega(y_i,q_1y_1)
-\ep_n(q_1y_1,y_2,\ldots,y_n)\omega(y_1,q_1y_1)
\prod_{i=2}^n\omega(y_i,q_1y_1)\biggl)\\
&=\dfrac{1}{n+1}\ep_{n-1}(y_2,\ldots,y_n)\omega(y_i,q_1y_1)
\biggl(\bigg[\prod_{i=2}^n\ep_2(y_1,y_i)\omega(y_i,q_1y_1)\bigg]
-\bigg[\prod_{i=2}^n\ep_2(q_1y_1,y_i)\omega(y_1,y_i)\bigg]\biggl),
\end{align*}
where  $[\ep_n,\ep_1]\seteq\ep_n*\ep_1-\ep_1*\ep_n$. 
Noting that 
\begin{align*}
\ep_2(y_1,y_i)\omega(y_i,q y_1)=\ep_2(q y_1,y_i)\omega(y_1,y_i)
=\dfrac{(y_1-q_1^{-2}y_i)(y_1-q_1y_i)(y_1-q_2y_i)(y_1-q_3y_i)}
  {(y_1-y_i)^2(y_1-q_1^{-1}y_i)^2},
\end{align*}
we have $\varphi_{(2,1^{n-1})}[\ep_n,\ep_1]=0$. 
Lemma \ref{prop:bottom} means that 
$[\ep_n,\ep_1](x)=c \ep_{n+1}(x)$
with $c\in\bbF$. 
Using the specialization with respect to $q_2$, we have
\begin{align*}
&\varphi_{(n+1)}^{(q_2)}[\ep_n,\ep_1]=0,\qquad 
\varphi_{(n+1)}^{(q_2)} \ep_{n+1}
=\prod_{1\le i<j\le n+1}
\dfrac{(q_2^i-q_1q_2^j)(q_2^i-q_1^{-1}q_2^j)}{(q_2^i-q_2^j)^2}
\neq0,
\end{align*}
which indicates that $c=0$. Hence we have $\ep_n*\ep_1=\ep_1*\ep_n$.

Next we examine the commutator $[\ep_m,\ep_n]$ in general. Consider 
the specialization 
\begin{align*}
&\varphi^{(q_1)}_{(2,1^{m+n-2})}[\ep_m,\ep_n]\\
&=\dfrac{1}{(m+n-2)!}\sum_{\sigma\in\frakS_{m+n-2}}
\sigma\biggl(\ep_m(y_1,y_2,\ldots,y_m)\ep_n(q_1 y_1,y_{m+1},\ldots,y_{m+n-1})\\
&\hskip 3em
\bigg[\prod_{\alpha=2}^{m}\omega(y_\alpha,q_1y_1)\bigg]
\bigg[\prod_{\beta=n+1}^{m+n}\omega(y_1,y_\beta)\bigg]
\bigg[\prod_{\substack{2\le\alpha\le m\\m+1\le \beta\le m+n}}
\omega(y_\alpha,y_\beta)\bigg]\biggl)\\
&\hskip 3em
-\dfrac{1}{(m+n-2)!}
\sum_{\sigma\in\frakS_{m+n-2}}
\sigma\biggl(\ep_m(q_1 y_1,y_2,\ldots,y_m)\ep_n(y_1,y_{m+1},\ldots,y_{m+n-1})\\
&\hskip 3em
\bigg[\prod_{\alpha=2}^{m}\omega(y_1,y_\alpha)\bigg]
\bigg[\prod_{\beta=m+1}^{m+n}\omega(q_1 y_1,y_\beta)\bigg]
\bigg[\prod_{\substack{2\le\alpha\le m\\m+1\le \beta\le m+n}}
\omega(y_\beta,y_\alpha)\bigg]\biggl)\\
&=\dfrac{1}{(m+n-2)!}
\sum_{\sigma\in\frakS_{m+n-2}}
\sigma\biggl(
\ep_{m-1}(y_2,\ldots,y_m)\ep_{n-1}(y_{m+1},\ldots,y_{m+n-1})
\\
&
\bigg\{
\bigg[\prod_{\alpha=2}^{m}\ep_2(y_1,y_\alpha)\bigg]
\bigg[\prod_{\beta=m+1}^{m+n-1}\ep_2(q_1y_1,y_\beta)\bigg]
\bigg[\prod_{\alpha=2}^{m}\omega(y_\alpha,q_1y_1)\bigg]
\bigg[\prod_{\beta=m+1}^{m+n-1}\omega(y_1,y_\beta)\bigg]
\\
&\hskip 3em
\cdot
\bigg[\prod_{\substack{2\le\alpha\le m\\m+1\le \beta\le m+n}}
\omega(y_\alpha,y_\beta)\bigg]-\\
&
\bigg[\prod_{\alpha=2}^{m}\ep_2(q_1y_1,y_\alpha)\bigg]
\bigg[\prod_{\beta=m+1}^{m+n-1}\ep_2(y_1,y_\beta)\bigg]
\bigg[\prod_{\alpha=2}^{m}\omega(y_1,y_\alpha)\bigg]
\bigg[\prod_{\beta=m+1}^{m+n-1}\omega(y_\beta,q_1y_1)\bigg]
\\
&\hskip 3em
\cdot
\bigg[\prod_{\substack{2\le\alpha\le m\\m+1\le \beta\le m+n}}
 \omega(y_\beta,y_\alpha)\bigg]
\bigg\}\biggl),
\end{align*}
where $\sigma\in\frakS_{m+n-2}$ acts on $(y_2,\ldots,y_{m+n-1})$. 
The equation
\begin{align*}
&\prod_{\alpha=2}^{m}\ep_2(y_1,y_\alpha)
\prod_{\beta=m+1}^{m+n-1}\ep_2(q_1y_1,y_\beta)
\prod_{\alpha=2}^{m}\omega(y_1,y_\beta)
\prod_{\beta=m+1}^{m+n-1}\omega(y_\alpha,q_1y_1)\\
=
&\prod_{\alpha=2}^{m}\ep_2(q_1y_1,y_\alpha)
\prod_{\beta=m+1}^{m+n-1}\ep_2(y_1,y_\beta)
\prod_{\alpha=2}^{m}\omega(y_1,y_\alpha)
\prod_{\beta=m+1}^{m+n-1}\omega(y_\beta,q_1y_1)\\
=
&\prod_{i=2}^{m+n-1}
\dfrac{(y_1-q_1y_i)(y_1-q_2y_i)(y_1-q_3y_i)(y_1-q_1^{-2}y_i)}
{(y_1-y_i)^2(y_1-q_1^{-1}y_i)^2},
\end{align*}
shows that the image $\varphi^{(q_1)}_{(2,1^{m+n-2})}[\ep_m,\ep_n]$ is proportional to $[\ep_{m-1},\ep_{n-1}](y_2,\ldots,y_{m+n-1})$. By the induction hypothesis $[\ep_{m-1},\ep_{n-1}]=0$, we have $\varphi^{(q_1)}_{(2,1^{m+n-2})}[\ep_m,\ep_n]=0$, which implies that $[\ep_m,\ep_n]$ is a constant multiple of $\ep_{m+n}$. One can show that $\varphi_{(m+n)}^{(q_2)}[\ep_m,\ep_n]=0$, and $\varphi_{(m+n)}^{(q_2)}\ep_{m+n}\neq 0$ by a straightforward calculation. Hence we conclude that $\ep_m*\ep_n=\ep_n*\ep_m$.
\end{proof}

\subsection{Derivations $\partial^{(\infty,k)},\partial^{(0,k)} $ and the final step of the proof}
For $f\in \wbcalA_n$ 
we need to study the image $\partial^{(\infty,k)}f$ and 
 $\partial^{(0,k)}f$, including their existence. 
As for the existence, Proposition \ref{prop:hatA_2delta} assures.
To express the image of  $\partial^{(\infty,k)}$ and $\partial^{(0,k)} $
in a concise manner, we extend the definition of 
the star product as follows.
Let $f_1\in \whcalA_k$, $g_1\in \whcalA_l$, 
$f_2\in \whcalA_m$ and  $g_2\in \whcalA_n$, 
we set $(f_1\otimes g_1)*(f_2\otimes g_2):=(f_1*f_2)\otimes (g_1*g_2)\in 
\whcalA_{k+m}\otimes \whcalA_{l+n}$, and extend this by linearity.
In this extended case, the associativity of the star product 
holds as in Lemma \ref{lem:associative}.

\begin{lem}
\label{lem:derivation}
(1)
$\partial^{(0,k)}(\ep_n(x;q_i))=\partial^{(\infty,k)}(\ep_n(x;q_i))
=\ep_{n-k}(x_1,\ldots,x_{n-k})\otimes\ep_{k}(x_{n-k+1},\ldots,x_n)
\in \wbcalA_{n-k}\otimes \wbcalA_{k}$ for any $n,k\in\bbN$.
\\
(2)
For $f\in\whcalA_m$, $g\in\whcalA_n$ and $a=\infty$ or $0$, 
$\partial^{(a,k)}(f*g)\in \whcalA_{n-k}\otimes \whcalA_{k}$ exists 
and 
\begin{align*}
\partial^{(a,k)}(f*g)
=\sum_{i=0}^k \binom{k}{i}
(\partial^{(a,i)} f)*(\partial^{(a,k-i)} g).
\end{align*}
\end{lem}

\begin{proof}
(1) follows from a direct computation.
\\
(2)
We use tentative notation $y=(y_1,\ldots,y_{m+n})\seteq(x_1,\ldots,x_{m+n-k},\xi x_{m+n-k+1},\ldots,\xi x_{m+n})$ and $J\seteq\{m+n-k+1,\ldots,m+n\}$. Using the alternative expression (\ref{eq:star_rewrite}) yields
\begin{align*}
&\partial^{(\infty,k)}(f*g)(x)
=\dfrac{(m+n)!}{(m+n-k)!}\,
 \lim_{\xi\to\infty}\dfrac{m!n!}{(m+n)!}
  \sum_{\substack{I\subseteq\{1,\ldots,n+m\}\\ \# I =m}} f(y_I)\,g(y_{I^c})
   \prod_{\substack{\alpha\in I\\ \beta\in I^c}}\omega(y_\alpha,y_\beta)
\\
&=
 \dfrac{m!n!}{(m+n-k)!}\,
 \sum_{i=0}^k 
  \sum_{\substack{I\subseteq\{1,\ldots,n+m\}\\ \# I=m,\  \# (I\cap J)=i}}
  \lim_{\xi\to\infty}
   \biggl(
    f(y_I)\, g(y_{I^c})
    \prod_{\substack{\alpha\in I\\ \beta\in I^c}} \omega(y_\alpha,y_\beta)
   \biggr)
\\
&=
 \dfrac{(m-i)!(n-(k-i))!}{(m+n-k)!}\,
 \sum_{i=0}^k 
  \sum_{\substack{I\subseteq\{1,\ldots,n+m\}\\ \# I=m,\  \# (I\cap J)=i}}
  \biggl(
   \dfrac{m!}{(m-i)!}\big[\lim_{\xi\to\infty} f(y_I)\big]
   \dfrac{n!}{(n-(k-i))!}\big[\lim_{\xi\to\infty} g(y_{I^c})\big]
\\
&\hskip 20em
  \bigg[ \lim_{\xi\to\infty} 
   \prod_{\substack{\alpha\in I\\ \beta\in I^c}} \omega(y_\alpha,y_\beta)\bigg]
  \biggr)
\\
&=
 \binom{m+n-k}{m-i}^{-1}
 \sum_{i=0}^k 
  \sum_{\substack{I\subseteq\{1,\ldots,n+m\}\\ \# I=m,\  \# (I\cap J)=i}}
  \biggl(
   (\partial^{(\infty,i)} f) \, (\partial^{(\infty,k-i)}g)
    \prod_{\substack{\alpha\in I \cap J \\
                      \beta \in I^c \cap J}}
   \omega(x_\alpha,x_\beta)
   \prod_{\substack{\gamma\in I \cap J^c \\
                    \delta\in I^c \cap J^c}}
   \omega(x_\gamma,x_\delta)
  \biggr)
\\
&=
 \binom{m+n-k}{m-i}^{-1}
 \sum_{i=0}^k 
  \binom{k}{i} \binom{m+n-k}{m-i}
  \biggl(
   (\partial^{(\infty,i)} f)*(\partial^{(\infty,k-i)}g)
  \biggr)
\\
&=
 \sum_{i=0}^k 
  \binom{k}{i} 
  \biggl(
   (\partial^{(\infty,i)} f)*(\partial^{(\infty,k-i)}g)
  \biggr).
\end{align*}
Thus the limit $\partial^{(\infty,k)}(f*g)$ exists and the Leibniz rule are proved.  The case $\partial^{(0,k)}$ is the same.
\end{proof}

\begin{prop}
\label{prop:calA=wbcalA}
We have $\calA=\wbcalA$.
\end{prop}

\begin{proof}
$\calA_n\subset\wbcalA_n$ follows from the definitions, so we will prove the other direction $\wbcalA_n\subset\calA_n$. Since the family $\{\ep_\lambda(x;q_i)\mid \lambda\vdash n\}$ spans $\wbcalA_n$, it is enough to show that $\ep_\lambda(x;q_i)\in\calA$, that is, $\ep_\lambda$ satisfies the condition (ii) of Definition \ref{dfn:calA}. Set $\lambda=(\lambda_1,\lambda_2,\ldots,\lambda_l)$. Lemma \ref{lem:derivation} (2) yields
\begin{align*}
\partial^{(a,k)}(\ep_\lambda)
&=\sum \binom{m}{i_1,i_2,\ldots,i_l} 
\partial^{(a,i_1)}(\ep_{\lambda_1})*
\partial^{(a,i_2)}(\ep_{\lambda_2})*\cdots*
\partial^{(a,i_l)}(\ep_{\lambda_l})
\end{align*}
for $a=\infty$ and $0$. Then the conclusion follows from Lemma \ref{lem:derivation} (1).
\end{proof}

Hence we have proved Theorem \ref{thm:1}.

\section{Free field, Wronski relation, and 
proof of Theorems \ref{thm:2}\ and \ref{thm:3}}
\label{sect:mcd}

In this section we investigate the free field realization of the Macdonald difference operators,
which we need for our proof of Theorems \ref{thm:2}\ and \ref{thm:3}.
One more necessary ingredient is the {\it Wronski relation},
which gives a certain identity among operators acting on $\calF$.
 In what follows, we frequently use the notation for $q$-integers and the $q$-shifted factorial:
\begin{align*}
[n]_q\seteq \dfrac{1-q^n}{1-q},\quad
[n]_q!\seteq\prod_{k=1}^n[k]_q,\quad
(a;q)_k\seteq \prod_{i=0}^{k-1}(1-a q^i).
\end{align*}
We also use the infinite product 
\begin{align*}
(x;q)_\infty\seteq \prod_{i\ge 0}(1-x q^i),
\end{align*}
which is regarded as a formal power series in $x$ over $\bbF$.

\subsection{Preliminaries}
\label{subsec:pre}
We briefly recall some basic facts on the Macdonald symmetric polynomials/functions and on Macdonald's reproduction kernel.
Recall that we have denoted by $\Lambda_n$ the 
ring of symmetric polynomials in $x=(x_1,\ldots,x_n)$ over $\bbZ$, and 
by $\Lambda_n^k$ the space of homogeneous symmetric polynomials of degree $k$. 
The ring of symmetric functions $\Lambda$ is defined as the inverse limit of the $\Lambda_n$ in the category of graded rings. We set $\Lambda_\bbF=\Lambda\otimes_\bbZ \bbF$.

The Macdonald difference operator $D_n^r$ \cite[VI \S3]{M:1995:book}
acting on $\Lambda_{n,\bbF}$ is defined by
\begin{align}\label{eq:mdo}
D_n^r=D_n^r(q,t)\seteq t^{r(r-1)/2}
\sum_{\substack{I\subset\{1,2,\ldots,n\}\\ \# I=r}}
\prod_{\substack{i\in I \\ j\notin I}}\dfrac{t x_i-x_j}{x_i-x_j}
\prod_{k\in I}T_{q,x_k},
\end{align}
where $T_{q,x_i}$ denotes the $q$-difference operator 
$T_{q,x_i}f(x_1,\ldots,x_n)=f(x_1,\ldots,q x_i,\ldots,x_n)$.

\begin{prop}[Macdonald]
For $\lambda\vdash n$ such that $\ell(\lambda)\leq n$, $P_\lambda(x;q,t)\in \Lambda_{n,\bbF}$ is uniquely characterized by 
\begin{align}
\tag{c}
&P_\lambda(x;q,t)=m_\lambda+\sum_{\mu<\lambda} 
 c_{\lambda\mu}^{m\rightarrow P}m_\mu
 \qquad (c_{\lambda\mu}^{m\rightarrow P}\in \bbF),\\
\tag{d}
&D_n^r(q,t)P_\lambda(x;q,t)
=e_r^{(n)}(t^n s^\lambda_1,t^n s^\lambda_2,\ldots,t^n s^\lambda_n)
P_\lambda(x;q,t),
\end{align}
where $e_r^{(n)}(s_1\ldots,s_n)$ denotes the elementary symmetric polynomial in $n$ variables $(s_1,\ldots,s_n)$, and $s^\lambda_i\seteq q^{\lambda_i}t^{-i}$.
\end{prop}

Denote by $e_r(s)\in \Lambda_{\bbF,s}$ the $r$-th elementary symmetric function  in the {\it infinite} set of variables $s=(s_1,s_2,\ldots)$.

\begin{lem}
Let $\lambda\vdash n$ such that $\ell(\lambda)\leq n$.
Define the infinite sequence
$s^\lambda=(s_1^{\lambda},s_2^{\lambda},\ldots)$ by 
$s^\lambda_i \seteq t^{-i}q^{\lambda_i}$ and set the finite sequence 
$t^n s^\lambda\seteq(t^{n-1}q^{\lambda_1},\ldots,q^{\lambda_n})$. 
Then we have
\begin{align*}
e_r(s^\lambda)= \sum_{l=0}^r
\dfrac{t^{-n r-\binom{r-l+1}{2}}}{(t^{-1};t^{-1})_{r-l}} 
 e_l^{(n)}(t^n s^\lambda)\in \bbF[[t^{-1}]],
\end{align*}
where all the rational factors in the RHS such as $1/(1-t^{-k})$ ($k\geq 1$)
must be understood as the power series $\sum_{l\geq 0}t^{-kl}$.
\end{lem}

\begin{proof}
Introduce the generating functions for $e_r(s)$ and 
$e_r^{(n)}(s_1\ldots,s_n)$ as
\begin{align*}
E(s;u)=\prod_{i\ge 1} (1+s_i u)=\sum_{r\ge 0} e_r(s)u^r, 
\qquad E^{(n)}(s;u)=\prod_{i=1}^n (1+s_i u)=\sum_{r=0}^n e_r^{(n)}(s)u^r.
\end{align*}
Then we have
\begin{align*}
E(s^\lambda;u)
&=\prod_{i\ge 1}(1+s^\lambda_i u)
 =\prod_{i=1}^n (1+s^\lambda_i u) 
  \prod_{i\ge n+1} (1+t^{-i}u) 
=E^{(n)}(t^n s^\lambda ; t^{-n}u) \cdot 
(-t^{-1}t^{-n}u;t^{-1})_{\infty}\\
&= \sum_{l=0}^n e_l^{(n)}(t^n s^\lambda) \cdot 
   (t^{-n}u)^l \sum_{m \ge 0}  
   \frac{1}{(t^{-1};t^{-1})_m}(t^{-n-1}u)^m t^{-\binom{m}{2}}\\
&=
\sum_{r \ge 0} u^r \sum_{l=0}^r
\dfrac{t^{-n r-\binom{r-l+1}{2}}}{(t^{-1};t^{-1})_{r-l}} 
 e_l^{(n)}(t^n s^\lambda).
\end{align*}
\end{proof}

Thus we have proved 
\begin{prop}\label{prop:E_r}
Let $n\in \bbN_+$. Set a difference operator $E_r^{(n)}$ acting on $\Lambda_{n,\bbF}$ by 
\begin{align*}
E_r^{(n)}\seteq \sum_{l=0}^r
\dfrac{t^{-n r-\binom{r-l+1}{2}}}{(t^{-1};t^{-1})_{r-l}} D_n^l.
\end{align*}
For any partition $\lambda$ satisfying $\ell(\lambda)\leq n$, we have
\begin{align*}
E_r^{(n)} P_{\lambda}(x;q,t) = e_r(s^\lambda) P_{\lambda}(x;q,t),
\end{align*}
where $P_{\lambda}(x;q,t)\in \Lambda_{n,\bbF}$ denotes 
the Macdonald symmetric {\it polynomial}, 
$e_r(x)\in \Lambda_\bbF$ is the $r$-th elementary symmetric {\it function}
defined in (\ref{eq:E(y)}), 
and $s^\lambda=(t^{-1}q^{\lambda_1},t^{-2}q^{\lambda_2},\ldots)$. 
Hence the inductive limit $E_r:=\varprojlim_n E_r^{(n)}$ exists.
\end{prop}

\begin{rmk}
Note that if we renormalize the operator suitably as $ (t^{-1};t^{-1})_r E_r^{(n)}$,  then the modified eigenvalues $(t^{-1};t^{-1})_r e_r(s^\lambda)$ lie in $\bbF$.
\end{rmk}

Next we turn to the case of the symmetric function in $\Lambda_\bbF$. For a partition $\lambda$, denote by $m_i$ the number of parts $i$ in the partition $\lambda$. The scalar product on $\Lambda_\bbF$ is defined by \cite[VI.2.2]{M:1995:book}
\begin{align*}
\qtpr{p_\lambda,p_\mu}
\seteq
\delta_{\lambda,\mu}
\prod_{i \ge 1}i^{m_i}m_i!
\prod_{j \ge 1}\dfrac{1-q^{\lambda_j}}{1-t^{\lambda_j}},
\end{align*}

Set $b_\lambda=b_\lambda(q,t)\seteq\qtpr{P_\lambda(x;q,t),P_\lambda(x;q,t)}^{-1}$, and $Q_\lambda(x;q,t)\seteq b_\lambda(q,t)P_\lambda(x;q,t)$. Then $(Q_\lambda)$ forms a dual basis of $(P_\lambda)$. For two sets of independent indeterminates $x=(x_1,x_2,\ldots)$ and $y=(y_1,y_2,\ldots)$, the {\it reproduction kernel} is defined to be
\begin{align}\label{eq:kernel}
\Pi(x,y;q,t)\seteq \prod_{i,j\ge1}
 \dfrac{(t x_i y_j;q)_\infty}{(x_i y_j;q)_\infty}
\end{align}
Then we have \cite[VI.4.13]{M:1995:book}
\begin{align}\label{eq:tenkai}
\Pi (x,y;q,t)=\sum_\lambda P_\lambda(x;q,t)Q_\lambda(y;q,t).
\end{align}

\begin{lem} \label{lem:kernel}
Let $n\in\bbN_+$ and $y=(y_1,\cdots,y_n)$. Assume that an operator $\widehat{E}_r$ on the space of symmetric functions $\Lambda_\bbF$ in $x=(x_1,x_2,\ldots)$ satisfies $\widehat{E}_r \Pi(x,y)=E_{r,y}^{(n)} \Pi(x,y)$. Here the subscript $y$ indicates that the operator is acting on $y$. Then $\widehat{E}_r P_\lambda(x;q,t)=e_r(s^\lambda) P_\lambda(x;q,t)$ for any $\lambda\vdash n$ satisfying $\ell(\lambda)\leq n$.
\end{lem}

\begin{proof}
Note that $(Q_\lambda)_{\ell(\lambda)\leq n}$ is a basis of $\Lambda_{n,\bbF}$ and use the expansion (\ref{eq:tenkai}).
\end{proof}

The automorphism $\omega_{q,t}$ on $\Lambda_\bbF$ is defined by
\begin{align}\label{eq:omega_auto}
\omega_{q,t}(p_r)\seteq(-1)^{r-1}\dfrac{1-q^r}{1-t^r}p_r.
\end{align}
It is known that (\cite[(VI.2.15), (VI.5.1)]{M:1995:book})
\begin{align}\label{eq:dual}
\omega_{q,t}(g_n(x;q,t))=e_n(x),\quad
\omega_{q,t}(P_\lambda(x;q,t))=Q_{\lambda'}(x;t,q).
\end{align}

\subsection{Proofs of Propositions \ref{prop:calO}\ and \ref{prop:spectral_criterion}}
\label{subsec:prf_prop_calO_spec}

First we give proofs of Propositions \ref{prop:calO} and \ref{prop:spectral_criterion}. The former is easy, but the latter requires some preparation. 

\begin{proof}[{Proof of Proposition \ref{prop:calO}}]
Noting the symmetric property (\ref{eq:eta_ep}), we have
\begin{align*}
\calO&(f*g)=
\dfrac{1}{(m+n)!}
\bigg[
\sum_{\sigma\in\frakS_{m+n}}\sigma\biggl(
f(z_{1},\ldots,z_{m})
g(z_{m+1},\ldots,z_{m+n})\\
&\hskip 8em
\biggl(
\prod_{\substack{1\le i\le m\\m+1\le j\le m+n}}
\omega(z_{i},z_{j}) \biggr)
\biggl(\prod_{1\le i<j\le m+n}\omega(z_{i},z_{j})^{-1}\biggr)
\eta(z_{1})\cdots\eta(z_{m+n})\biggr)
\bigg]_1
\\
&=
\dfrac{1}{(m+n)!}
\bigg[
\sum_{\sigma\in\frakS_{m+n}}\sigma\biggl(
\dfrac{f(z_{1},\ldots,z_{m})}
{\displaystyle 
\prod_{1\le i<j\le m}\omega(z_{i},z_{j})}
\eta(z_{1})\cdots\eta(z_{m})
\\
&\hskip 12em
\dfrac{
g(z_{m+1},\ldots,z_{m+n})}
{\displaystyle 
\prod_{m+1\le k<l\le m+n}\omega(z_{k},z_{l})}
\eta(z_{m+1})\cdots\eta(z_{m+n})
\biggr)
\bigg]_1\\
&=
\dfrac{\#(\frakS_{m+n}/\frakS_m\times\frakS_n)}{(m+n)!}
\biggl[
\dfrac{f(z_1,\ldots,z_m)}{\displaystyle \prod_{1\le i<j\le m}\omega(z_i,z_j)}
\eta(z_{1})\cdots\eta(z_{m})
\biggr]_1
\biggl[
\dfrac{g(z_{1},\ldots,z_{n})}
{\displaystyle \prod_{1\le i<j\le m}\omega(z_i,z_j)}
\eta(z_{1})\cdots\eta(z_{n})
\biggr]_1
\\
&=
\dfrac{\#(\frakS_{m+n}/\frakS_m\times\frakS_n)}{(m+n)!/m!n!}
\calO(f)\calO(g)
=\calO(f)\calO(g).
\end{align*}
\end{proof}

To our proof of Proposition \ref{prop:spectral_criterion}, a key role is played by the following free field representation.

\begin{prop}[{\cite[Theorem 9.2]{S:2006}}]\label{prop:eigen:E_r}
Let $\ep_r(z;q)$ be the element in Definition \ref{dfn:bottom}. Let $E_r$ be the operator on $\Lambda_\bbF$ in Proposition \ref{prop:E_r}. Set
\begin{align}
\widehat{E}_r=\widehat{E}_r(q,t)
\seteq \dfrac{t^{-r(r+1)/2}}{(t^{-1};t^{-1})_{r}}
\dfrac{[r]_{t^{-1}}!}{r!}\calO(\ep_r(z;q)).
\end{align}
Then for any partition $\lambda$ we have 
\begin{align*}
\widehat{E}_r P_\lambda(x;q,t)=E_r P_\lambda(x;q,t).
\end{align*}
\end{prop}

We defer the proof of Proposition \ref{prop:eigen:E_r} to the next subsection, since it is a little lengthy, and we also need to continue similar developments further for later purpose. 

\begin{proof}[{Proof of Proposition \ref{prop:spectral_criterion}}]
Proposition \ref{prop:calO} means that $\calO(\cdot)$ is an algebra homomorphism. By construction, surjectivity is trivial. Changing the normalization, we set
\begin{align*}
f_{r}(z):=\dfrac{t^{-r(r+1)/2}}{(t^{-1};t^{-1})_{r}}
\dfrac{[r]_{t^{-1}}!}{r!}\ep_r(z;q),
\end{align*}
and denote $f_\lambda\seteq f_{\lambda_1}*f_{\lambda_2}*\cdots$. Expand $f\in \calA_n$ by the basis $(\ep_\lambda(z;q))$ as $f=\sum_{\mu\,\vdash\, n} c_\mu f_\mu(z)$. From Propositions \ref{prop:eigen:E_r}, \ref{prop:calO}, we find that $\calO(f)P_\lambda(x;q,t)=\sum_{\mu\,\vdash\, n} c_\mu e_\mu(s^\lambda)\cdot P_\lambda(x;q,t)$ for any $\lambda$. As for the injectivity, we must show that $\sum_{\mu\,\vdash\, n} c_\mu e_\mu(s^\lambda)=0$ for any $\lambda$ indicates $f=0$. From Lemma \ref{lem:zero=zero} below, we have $\sum_{\mu\,\vdash\, n} c_\mu e_\mu(x)=0$ in $\Lambda_\bbF$. Hence $c_\mu=0$, showing that $f=0$.
\end{proof}

\begin{lem}\label{lem:zero=zero}
Let $k\in \bbN$, and set $s^\lambda=(s_1^{\lambda},s_2^{\lambda},\ldots)
=(t^{-1}q^{\lambda_1},t^{-2}q^{\lambda_2},\ldots)$. If $f\in \Lambda^k_\bbF$ satisfies the condition that $f(s^\lambda)=0$ for all $\lambda$, then $f=0$.
\end{lem}

\begin{proof}
We prove this by induction. For $f\in\Lambda_\bbF^0$, the statement is trivial. Let $f \in \Lambda_\bbF^k$ with $k\in\bbN_+$. Expand $f$ in $x_1$ as $f=\sum_{i=0}^k x_1^{k-i} g_i(x_2,x_3,\ldots)$. Then we have $g_i(s^{\lambda}_2,s^\lambda_3,\ldots)\in \bbF[[t^{-1}]]$. Since the ring $\bbF[[t^{-1}]]$ is an integral domain, the polynomial remainder theorem applies. For a fixed  $(\lambda_2,\lambda_3,\ldots)$, consider the infinite sequence of evaluations 
$x_1=t^{-1}q^{\lambda_2},t^{-1}q^{\lambda_2+1},\ldots $. Then one finds that $ g_i(x_2,x_3,\ldots)=0$ for $i=0,1,\ldots,k-1$ by the induction hypothesis. Hence $f$ becomes a symmetric function which does not contain $x_1$. From the symmetry in $x_i$'s we conclude $f=0$.
\end{proof}

\subsection{Proof of Proposition \ref{prop:eigen:E_r}}
\label{subsec:ffr}
We recall the notation and a lemma in \cite{S:2006}.
\begin{dfn}
Set
\begin{align}
\label{eq:phi}
\phi(y) \seteq 
\exp\biggl(\sum_{n \ge 1} \dfrac{1-t^n}{1-q^n}\dfrac{a_{-n}}{n}y^n\biggr)
=\prod_{i\ge1}\dfrac{(t x_i y;q)_\infty}{(x_i y;q)_\infty}
=\sum_{n\ge0}g_n(x;q,t)y^n.
\end{align}
Here we have used the identification $p_n(x)=a_{-n}$.
\end{dfn}
%
\begin{lem}[{\cite[Lemma 9.3]{S:2006}}]
\label{lem:basic_tools}
We have
\begin{align*}
&
\phi(y_1)\cdots \phi(y_n) \cdot 1=\Pi(x_1,x_2,\ldots,y_1,\ldots,y_n;q,t), 
\\
&\eta(z)\phi(y)\cdot 1 
= \frac{1-y/z}{1-t y/z} : \phi(y)\eta(z) : \cdot\, 1,\quad
:\eta(t y)\phi(y): \cdot 1 = T_{q,y} \phi(y) \cdot 1,
\nonumber
\end{align*}
where $1$ denotes the highest vector in the Fock space and $\Pi$ is the reproduction kernel (\ref{eq:kernel}).
\end{lem}

\begin{proof}[{Proof of Proposition \ref{prop:eigen:E_r}}]
Our argument here is based on Lemma \ref{lem:kernel} and the 
reproduction kernel expressed as in the 
first relation of Lemma \ref{lem:basic_tools}. We examine the action of $\widehat{E}_r$ on the $y_i$'s. For simplicity of display, we use the abbreviation 
$\phi_y = \phi(y_1)\ldots\phi(y_n) $.
Note first that
\begin{align*}
\dfrac{[r]_{t^{-1}}!}{r!}\ep_r(z;q)
=\Sym\biggl(\prod_{1\le i<j\le r}\dfrac{1-z_j/z_i}{1-t^{-1}z_j/z_i}\biggr).
\end{align*}
Using (\ref{eq:eta_ope}), we have
\begin{align*}
\widehat{E}_r \phi_y\cdot 1
=&\dfrac{t^{-r(r+1)/2}}{(t^{-1};t^{-1})_{r}}
\biggl[
\biggl(
   \prod_{j=1}^r \prod_{i=1}^n \frac{1-y_i/z_j}{1-t y_i/z_j} \biggr)
   \biggl(
   \prod_{1 \le i<j\le r} \dfrac{1-z_j/z_i}{1-t^{-1}z_j/z_i}\biggr) \phi_y
   : \eta(z_1)\cdots \eta(z_r) : \cdot\, 1 
  \biggl]_1.
\end{align*}
The non-trivial residues occur at $z_r=t y_{k}$ ($1 \le k \le n$) and at $z_r=0$. Thus we can proceed as
\begin{align*}
\widehat{E}_r \phi_y\cdot 1
=& 
 \dfrac{t^{-r(r+1)/2}}{(t^{-1};t^{-1})_{r}}
 \sum_{k=1}^n 
  \biggl[
   (1-t^{-1})
   \biggl(\prod_{j=1}^{r-1} \prod_{i=1}^n \frac{1-y_i/z_j}{1-t y_i/z_j}\biggr)
   \biggl(\prod_{i \neq k} \dfrac{1-y_i/t y_k}{1-y_i/y_k} \biggr)
   \biggl(\prod_{1\leq i < j \leq r-1}  \dfrac{1-z_j/z_i}{1-t^{-1}z_j/z_i}
   \biggr)
\\
& \hskip 8em
   \times  
   \biggl(\prod_{1\leq i \leq r-1} \dfrac{1-t y_k/z_i}{1-y_k/z_i} \biggr)
   \cdot\phi_y : \eta(z_1)\cdots \eta(z_{r-1}) \eta(t y_k): \cdot\, 1 
  \biggl]_1 
\\
  &+\dfrac{t^{-r(r+1)/2}}{(t^{-1};t^{-1})_{r}}
    \biggl[ 
     t^{-n} 
     \biggl(\prod_{j=1}^{r-1}\prod_{i =1}^n \dfrac{1-y_i/z_j}{1-t y_i/z_j}
     \biggr)
     \biggl(\prod_{1\le i <j \le r-1} \dfrac{1-z_j/z_i}{1-t^{-1}z_j/z_i}
     \biggr)
\\
& \hskip 8em
     \times  
     \phi_y : \eta(z_1)\cdots \eta(z_{r-1}) : \cdot\, 1 
    \biggl]_1 
\\
=&
 \dfrac{t^{-r}(1-t^{-1})}{1-t^{-r}}\sum_{k=1}^n \prod_{i \neq k} 
   \dfrac{1-y_i/t y_k}{1-y_i/y_k}\,\phi(q y_k)\, 
 \widehat{E}_{r-1}\phi_{{\check{y_k}}} \cdot 1
   +\dfrac{t^{-n-r}}{1-t^{-r}}  \widehat{E}_{r-1}\phi_y \cdot 1,
\end{align*}
where $\phi_{{\check{y_k}}}\seteq\phi(y_1)\cdots\phi(y_{k-1})\phi(y_{k+1})\cdots \phi(y_n)$. By the induction hypothesis on $\widehat{E}_{r-1}$, we have
\begin{align}
\label{eq:tate:tmp3}
\begin{split}
\widehat{E}_r \phi_y\cdot 1
&=\bigg\{
   \dfrac{t^{-r}(1-t^{-1})}{1-t^{-r}} 
   \sum_{k=1}^n \prod_{i \neq k} \dfrac{1-y_i/t y_k}{1-y_i/y_k}T_{q,y_k}
   \biggl( 
    \sum_{l=0}^{r-1} 
    \dfrac{t^{-(n-1)(r-1)-\binom{r-l}{2}}}{(t^{-1};t^{-1})_{r-l-1}}
    D_{n-1,\check{y_k}}^{l}
   \biggr) \phi_y\\
&\hskip 2em
   +\dfrac{t^{-n-r}}{1-t^{-r}}  
   \biggl(
    \sum_{l=0}^{r-1} 
    \dfrac{t^{-n(r-1)-\binom{r-l}{2}}}{(t^{-1};t^{-1})_{r-l-1}} D_{n,y}^{l}
   \biggr) \phi_y
  \bigg\} \cdot 1,
\end{split}
\end{align}
where $D_{n-1,\check{y_k}}^{l}$ denotes the $l$-th Macdonald operator acting on the $n-1$ variables $y_1,\ldots,y_{k-1}$, $y_{k+1},\ldots,y_n$,
and $D_{n,y}^{l}$ denotes the $l$-th Macdonald operator acting on $y_1\ldots,y_n$. The first term of (\ref{eq:tate:tmp3}) can be computed as follows.
\begin{align*}
&  \sum_{k=1}^n \prod _{m \neq k} 
   \dfrac{t y_k - y_m}{y_k-y_m} T_{q,y_k} D_{n-1,\check{y_k}}^{l}
=  t^{\binom{l}{2}}\sum_{k=1}^n 
   \sum_{\substack{I \subset \{ 1,\ldots,n \} \\ k\notin I \\ |I|=l}} 
    \prod_{\substack{i \in I \\ j \not\in I \cup \{k\} }} 
     \dfrac{t y_i-y_j}{y_i-y_j}
    \prod_{m \neq k} \dfrac{t y_k-y_m}{y_k-y_m} 
    \prod_{i \in I \cup \{k\} }T_{q,y_i }\\
&\hskip 1em
=  t^{\binom{l}{2}}
    \sum_{k=1}^n 
    \sum_{\substack{J \subset \{ 1,\ldots,n \} \\ k \in J \\|J|=l+1}} 
     \prod_{\substack{i \in J \\ j \not\in J}}
      \dfrac{t y_i-y_j}{y_i-y_j} 
     \prod_{m \in J \setminus \{k\}}
      \dfrac{t y_k-y_m}{y_k-y_m} 
     \prod_{i \in J}T_{q,y_i}\\
&\hskip 1em
=  t^{\binom{l+1}{2}}t^{-l} 
    \sum_{\substack{J \subset \{1,\ldots,n\} \\ |J|=l+1}} 
     \biggl(
      \sum_{k \in J}
       \prod_{m \in J \setminus \{k\}} \dfrac{t y_k-y_m}{y_k-y_m} 
     \biggr)
     \prod_{\substack{i\in J \\ j \not\in J}} \dfrac{t y_i-y_j}{y_i-y_j} 
     \prod_{i \in J}T_{q,y_i} \\
&\hskip 1em
=  \dfrac{1-t^{-l-1}}{1-t^{-1}}\cdot D_{n,y}^{l+1},
\end{align*}
where $D_{n,y}^{l+1}$ denotes the $(l+1)$-st Macdonald operator acting on $y_1\ldots,y_n$. Here we have used the partial fraction expansion
\begin{align*}
\sum_{i=1}^n \prod_{j \neq i} \frac{t y_j-y_i}{y_j-y_i} =\frac{1-t^n}{1-t}.
\end{align*}

Hence we have
\begin{align*}
\widehat{E}_r \phi_y\cdot 1
&=\bigg\{
   \dfrac{t^{-r-n+1}(1-t^{-1})}{1-t^{-r}} 
   \sum_{l=0}^{r-1} 
    \dfrac{t^{-(n-1)(r-1)-\binom{r-l}{2}}}{(t^{-1};t^{-1})_{r-l-1}}
    \dfrac{1-t^{-l-1}}{1-t^{-1}}
    D_{n,y}^{l+1}\\
&\hskip 2em
   +\dfrac{t^{-n-r}}{1-t^{-r}}  
    \sum_{l=0}^{r-1} 
    \dfrac{t^{-n(r-1)-\binom{r-l}{2}}}{(t^{-1};t^{-1})_{r-l-1}}
     D_{n,y}^{l}
   \bigg\}\phi_y\cdot 1\\
&=
 \sum_{l=0}^{r} 
 \dfrac{t^{-n r-\binom{r-l+1}{2}}}{(t^{-1};t^{-1})_{r-l}}
  D_{n,y}^{l} \phi_y\cdot 1= E_r  \phi_y\cdot 1.
\end{align*}
\end{proof}

\subsection{Wronski relation of the Macdonald difference operators}\label{subsec:wronski}
Next, we proceed to the construction of the Wronski relation,
which is a key to the proof of Theorems \ref{thm:2}\ and \ref{thm:3}.

\begin{dfn}
Define the operator $\widehat{G}_n$ on the Fock space $\calF$ for $n\in\bbN$ by
\begin{align}
&\widehat{G}_n=\widehat{G}_n(q,t)\seteq
\dfrac{(-1)^n q^{\binom{n}{2}}}{(q;q)_n}
\dfrac{[n]_q!}{n!}\calO(\ep_n(z;t)). \label{eq:G_n}
\end{align}
\end{dfn}

\begin{prop}\label{prop:Wronski:op}
For $n\in\bbN$, we have the {\it Wronski relation} 
\begin{align}\label{eq:Wronski:op}
\sum_{k=0}^n(-1)^k(1-q^{k}t^{n-k})\widehat{E}_{n-k}\widehat{G}_{k}=0.
\end{align}
\end{prop}

\begin{proof}
For a while, we denote the dependence on the parameters by $q_1$ and $q_2$ (instead of $q^{-1}$ and $t$) for simplicity of display. 
We will show the next identity in $\calA$
\begin{align}\label{eq:sn}
S_n:=\sum_{k=0}^n\binom{n}{k}(q_1^{k}-q_2^{n-k})
\biggl(\dfrac{q_2-1}{q_1-1}\biggr)^k
\ep_{n-k}(x;q_1)*\ep_{k}(x;q_2)=0,
\end{align}
which is equivalent to (\ref{eq:Wronski:op}).
It is enough to verify that for any $\lambda=(\lambda_1,\ldots,\lambda_l)\vdash n$ one has
\begin{align}\label{eq:wronski:ind}
\varphi_\lambda^{(q_1)}(S_n)=0.
\end{align}
Our proof is an induction according to the dominance ordering of the partitions. We will fix the index $i=1$ and suppress the symbol `$(q_1)$'. 

First we show $\varphi_{(n)}(S_n)=0$. Skipping some trivial algebra, we have
\begin{align*}
S_n(1,q_1,\ldots,q_1^{n-1})
&=(q_1^{n-1}-q_2)\biggl(\dfrac{q_2-1}{q_1-1}\biggr)^{n-1}
  \ep_{n-1}(q_1,\ldots,q_1^{n-1};q_2)
  \prod_{i=1}^{n-1}\omega(1,q_1^i)
\\
&\hskip 1em
 +(q_1^n-1)\biggl(\dfrac{q_2-1}{q_1-1}\biggr)^{n}
   \ep_n(1,q_1,\ldots,q_1^{n-1};q_2)
\\
&=\biggl(\dfrac{q_2-1}{q_1-1}\biggr)^{n-1}
  \ep_{n-1}(q_1,\ldots,q_1^{n-2};q_2)
\\
&\hskip 1em 
  \cdot
  \biggl(
   (q_1^{n-1}-q_2) \prod_{i=1}^{n-1}\omega(1,q_1^i)
   +(q_1^n-1)\left(\dfrac{q_2-1}{q_1-1}\right)
    \prod_{i=1}^{n-1}\ep_2(1,q_1^i;q_2)
  \biggr)
=0.
\end{align*}

Next, fix a partition $\lambda$ and assume $\varphi_\nu(S_n)=0$ for any $\nu\ge\lambda$. Then the structure of the Gordon filtration indicates
\begin{align*}
S_n\in\bigcup_{\mu\,\vdash n,\ \mu <_1\lambda}\calA_\mu.
\end{align*}
Here the symbol $\mu <_1 \lambda$ denotes that $\mu<\lambda$ and that there is no partition $\nu\vdash n$ satisfying $\mu<\nu<\lambda$. Therefore we can write
\begin{align*}
S_n=\sum_{\mu\,\vdash n,\ \mu<_1\lambda}L_\mu/\Delta_n(x)^2,\quad
L_\mu\in\calP_\mu.
\end{align*}
For a partition $\mu$ satisfying $\mu<_1\lambda$, consider the specialization map $\varphi_\mu$. This map yields
\begin{align}\label{eq:wronski:unique}
\varphi_\mu(S_n)
=\sum_{\nu\,\vdash n,\ \nu<_1\lambda}\varphi_\mu(L_\nu/\Delta_n(x)^2)
=\varphi_\mu(L_\mu/\Delta_n(x)^2).
\end{align}
Using Lemma \ref{lem:atmost1}, we can write
\begin{align}\label{eq:wronski:proportional}
\varphi_\mu(L_\mu/\Delta_n(x)^2)=c\,\zeta_\mu/\varphi_\mu(\Delta_n(x)^2)
\end{align}
where $c\in\bbF$ and $\zeta_\mu$ is the polynomial (\ref{eq:zeta}).

We prepare the next definition and lemma.

\begin{dfn}[Snake evaluation]\label{dfn:snake}
For a partition $\lambda=(\lambda_1,\ldots,\lambda_l)$ of $n$, we define the specialization map $\psi_{\lambda} : \bbF(y_1,\ldots,y_l) \to \bbF$ by $y_i\mapsto q_1^{\sum_{k=i+1}^l (\lambda_k-1)} q_2^{i-1}$.
We call this specialization the {\it snake evaluation}.
\end{dfn}

\begin{lem}\label{lem:wronski}
For any $\mu=(\mu_1,\ldots,\mu_l)\vdash n$ we have the following.\\
(1) $\psi_{\mu}\circ\varphi_\mu(S_n)=0$.
\\
(2)
$\psi_{\mu}(\zeta_\mu)\neq0$.
\end{lem}

The proof of the lemma will be given afterwards. 
By the equation (\ref{eq:wronski:proportional}) and the lemma, we have $c=0$,
meaning that $\varphi_\mu(S_n)=0$ from (\ref{eq:wronski:unique}). 
Since $\mu$ was freely chosen from those partitions $\nu$ such that $\nu<_1\lambda$, 
we have completed the induction step.
\end{proof}

\begin{rmk}\label{rmk:snake}
We can give a diagrammatic explanation of the snake evaluation. First note that the composition $\psi_{\lambda}\circ\varphi_\lambda : \bbF(x_1,\ldots,x_n) \to\ \bbF$
can be written as
\begin{align*}
\begin{array}{l l l l l l l }
 \psi_{\lambda}\circ\varphi_\lambda(f(x_1,\ldots,x_n)) = 
&f(q_1^{(\lambda_{l}-1)+\cdots+(\lambda_2-1)},
&q_1^{(\lambda_{l}-1)+\cdots+\lambda_2},
&\ldots,&q_1^{(\lambda_{l}-1)+\cdots+(\lambda_1-1)},
\\
&\hskip1em
 q_1^{(\lambda_{l}-1)+\cdots+(\lambda_3-1)}q_2,
&q_1^{(\lambda_{l}-1)+\cdots+\lambda_3}q_2,
&\ldots,&q_1^{(\lambda_{l}-1)+\cdots+(\lambda_2-1)}q_2,
\\
&\hskip1em
&
&\ldots,
&
\\
&\hskip1em
 q_1^{\lambda_l-1}q_2^{l-2},
&q_1^{\lambda_l}q_2^{l-1},
&\ldots,
&q_1^{(\lambda_l-1)+(\lambda_{l-1}-1)}q_2^{l-1},
\\
&\hskip1em
 q_2^{l-1},
&q_1q_2^{l-1},
&\ldots,
&q_1^{\lambda_l-1}q_2^{l-1}).
\end{array}
\end{align*}
Consider the partition $\lambda=(4,4,2,1,1,1)$ (Figure \ref{fig:442111}). 
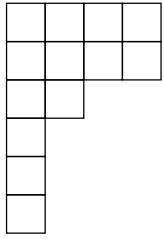
\begin{figure}[htbp]
\unitlength 0.1in
\begin{center}
\begin{picture}(  8.0000,  11.0000)(  4.0000,-15.0000)
\special{pn 8}%
\special{pa 400 400}%
\special{pa 600 400}%
\special{pa 600 600}%
\special{pa 400 600}%
\special{pa 400 400}%
\special{fp}%
\special{pn 8}%
\special{pa 600 400}%
\special{pa 800 400}%
\special{pa 800 600}%
\special{pa 600 600}%
\special{pa 600 400}%
\special{fp}%
\special{pn 8}%
\special{pa 800 400}%
\special{pa 1000 400}%
\special{pa 1000 600}%
\special{pa 800 600}%
\special{pa 800 400}%
\special{fp}%
\special{pn 8}%
\special{pa 1000 400}%
\special{pa 1200 400}%
\special{pa 1200 600}%
\special{pa 1000 600}%
\special{pa 1000 400}%
\special{fp}%
\special{pn 8}%
\special{pa 400 600}%
\special{pa 600 600}%
\special{pa 600 800}%
\special{pa 400 800}%
\special{pa 400 600}%
\special{fp}%
\special{pn 8}%
\special{pa 600 600}%
\special{pa 800 600}%
\special{pa 800 800}%
\special{pa 600 800}%
\special{pa 600 600}%
\special{fp}%
\special{pn 8}%
\special{pa 800 600}%
\special{pa 1000 600}%
\special{pa 1000 800}%
\special{pa 800 800}%
\special{pa 800 600}%
\special{fp}%
\special{pn 8}%
\special{pa 1000 600}%
\special{pa 1200 600}%
\special{pa 1200 800}%
\special{pa 1000 800}%
\special{pa 1000 600}%
\special{fp}%
\special{pn 8}%
\special{pa 400 800}%
\special{pa 600 800}%
\special{pa 600 1000}%
\special{pa 400 1000}%
\special{pa 400 800}%
\special{fp}%
\special{pn 8}%
\special{pa 600 800}%
\special{pa 800 800}%
\special{pa 800 1000}%
\special{pa 600 1000}%
\special{pa 600 800}%
\special{fp}%
\special{pn 8}%
\special{pa 400 1000}%
\special{pa 600 1000}%
\special{pa 600 1200}%
\special{pa 400 1200}%
\special{pa 400 1000}%
\special{fp}%
\special{pn 8}%
\special{pa 400 1200}%
\special{pa 600 1200}%
\special{pa 600 1400}%
\special{pa 400 1400}%
\special{pa 400 1200}%
\special{fp}%
\special{pn 8}%
\special{pa 400 1400}%
\special{pa 600 1400}%
\special{pa 600 1600}%
\special{pa 400 1600}%
\special{pa 400 1400}%
\special{fp}%
\end{picture}%
\end{center}
\caption{The Young diagram for $(4,4,2,1,1,1)$}
\label{fig:442111}
\end{figure}
\\
Then make a border strip (i.e. connected skew young diagram containing no $2\times 2$ block of squares \cite[I \S1, pp.5]{M:1995:book}) by shifting each row of the diagram of $\lambda$ (Figure \ref{fig:shift}). 
\begin{figure}[htbp]
\unitlength 0.1in
\begin{center}
\begin{picture}( 16.0000, 11.0000)(  2.0000,-13.0000)
\special{pn 8}%
\special{pa 200 1400}%
\special{pa 400 1400}%
\special{pa 400 1200}%
\special{pa 200 1200}%
\special{pa 200 1400}%
\special{fp}%
\special{pn 8}%
\special{pa 400 1200}%
\special{pa 200 1200}%
\special{pa 200 1000}%
\special{pa 400 1000}%
\special{pa 400 1200}%
\special{fp}%
\special{pn 8}%
\special{pa 200 1000}%
\special{pa 400 1000}%
\special{pa 400 800}%
\special{pa 200 800}%
\special{pa 200 1000}%
\special{fp}%
\special{pn 8}%
\special{pa 200 1000}%
\special{pa 400 1000}%
\special{pa 400 800}%
\special{pa 200 800}%
\special{pa 200 1000}%
\special{fp}%
\special{pn 8}%
\special{pa 400 800}%
\special{pa 600 800}%
\special{pa 600 600}%
\special{pa 400 600}%
\special{pa 400 800}%
\special{fp}%
\special{pn 8}%
\special{pa 200 800}%
\special{pa 400 800}%
\special{pa 400 600}%
\special{pa 200 600}%
\special{pa 200 800}%
\special{fp}%
\special{pn 8}%
\special{pa 400 600}%
\special{pa 600 600}%
\special{pa 600 400}%
\special{pa 400 400}%
\special{pa 400 600}%
\special{fp}%
\special{pn 8}%
\special{pa 600 600}%
\special{pa 800 600}%
\special{pa 800 400}%
\special{pa 600 400}%
\special{pa 600 600}%
\special{fp}%
\special{pn 8}%
\special{pa 800 600}%
\special{pa 1000 600}%
\special{pa 1000 400}%
\special{pa 800 400}%
\special{pa 800 600}%
\special{fp}%
\special{pn 8}%
\special{pa 1000 600}%
\special{pa 1200 600}%
\special{pa 1200 400}%
\special{pa 1000 400}%
\special{pa 1000 600}%
\special{fp}%
\special{pn 8}%
\special{pa 800 600}%
\special{pa 1000 600}%
\special{pa 1000 400}%
\special{pa 800 400}%
\special{pa 800 600}%
\special{fp}%
\special{pn 8}%
\special{pa 1000 400}%
\special{pa 1200 400}%
\special{pa 1200 200}%
\special{pa 1000 200}%
\special{pa 1000 400}%
\special{fp}%
\special{pn 8}%
\special{pa 1200 400}%
\special{pa 1400 400}%
\special{pa 1400 200}%
\special{pa 1200 200}%
\special{pa 1200 400}%
\special{fp}%
\special{pn 8}%
\special{pa 1000 400}%
\special{pa 1200 400}%
\special{pa 1200 200}%
\special{pa 1000 200}%
\special{pa 1000 400}%
\special{fp}%
\special{pn 8}%
\special{pa 1400 400}%
\special{pa 1600 400}%
\special{pa 1600 200}%
\special{pa 1400 200}%
\special{pa 1400 400}%
\special{fp}%
\special{pn 8}%
\special{pa 1600 400}%
\special{pa 1800 400}%
\special{pa 1800 200}%
\special{pa 1600 200}%
\special{pa 1600 400}%
\special{fp}%
\special{pn 8}%
\special{pa 200 400}%
\special{pa 400 400}%
\special{pa 400 200}%
\special{pa 200 200}%
\special{pa 200 400}%
\special{dt 0.03}%
\special{pn 8}%
\special{pa 400 400}%
\special{pa 600 400}%
\special{pa 600 200}%
\special{pa 400 200}%
\special{pa 400 400}%
\special{dt 0.03}%
\special{pn 8}%
\special{pa 600 400}%
\special{pa 800 400}%
\special{pa 800 200}%
\special{pa 600 200}%
\special{pa 600 400}%
\special{dt 0.03}%
\special{pn 8}%
\special{pa  800 400}%
\special{pa 1000 400}%
\special{pa 1000 200}%
\special{pa  800 200}%
\special{pa  800 400}%
\special{dt 0.03}%
\special{pn 8}%
\special{pa 200 600}%
\special{pa 400 600}%
\special{pa 400 400}%
\special{pa 200 400}%
\special{pa 200 600}%
\special{dt 0.03}%
\special{pn 8}%
\special{pa 800 600}%
\special{pa 1000 600}%
\special{pa 1000 400}%
\special{pa 800 400}%
\special{pa 800 600}%
\special{fp}%
\special{pn 8}%
\special{pa 800 600}%
\special{pa 1000 600}%
\special{pa 1000 400}%
\special{pa 800 400}%
\special{pa 800 600}%
\special{fp}%
\end{picture}%
\end{center}
\caption{The border strip made from the diagram for $(4,4,2,1,1,1)$}
\label{fig:shift}
\end{figure}
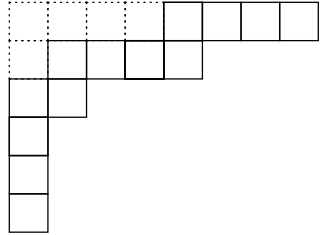
\\
Now we assign to each diagram a monomial $q_1^i q_2^j$ so that $i$ increases if one reads the boxes rightward and that $j$ increases if one reads downward (Figure \ref{fig:monom}).
\begin{figure}[htbp]
\begin{center}
\unitlength 0.1in
\begin{picture}(32.0, 18.0)(-3.0,-19.0)
\special{pn 8}%
\special{pa 200 200}%
\special{pa 500 200}%
\special{pa 500 500}%
\special{pa 200 500}%
\special{pa 200 200}%
\special{dt 0.03}%
\special{pn 8}%
\special{pa 200 500}%
\special{pa 500 500}%
\special{pa 500 800}%
\special{pa 200 800}%
\special{pa 200 500}%
\special{dt 0.03}%
\special{pn 8}%
\special{pa 200  800}%
\special{pa 500  800}%
\special{pa 500 1100}%
\special{pa 200 1100}%
\special{pa 200  800}%
\special{fp}%
\special{pn 8}%
\special{pa 200 1100}%
\special{pa 500 1100}%
\special{pa 500 1400}%
\special{pa 200 1400}%
\special{pa 200 1100}%
\special{fp}%
\special{pn 8}%
\special{pa 200 1400}%
\special{pa 500 1400}%
\special{pa 500 1700}%
\special{pa 200 1700}%
\special{pa 200 1400}%
\special{fp}%
\special{pn 8}%
\special{pa 200 1700}%
\special{pa 500 1700}%
\special{pa 500 2000}%
\special{pa 200 2000}%
\special{pa 200 1700}%
\special{fp}%
%
%
\special{pn 8}%
\special{pa 500  800}%
\special{pa 800  800}%
\special{pa 800 1100}%
\special{pa 500 1100}%
\special{pa 500  800}%
\special{fp}%
\special{pn 8}%
\special{pa 500 500}%
\special{pa 800 500}%
\special{pa 800 800}%
\special{pa 500 800}%
\special{pa 500 500}%
\special{fp}%
\special{pn 8}%
\special{pa  800 500}%
\special{pa 1100 500}%
\special{pa 1100 800}%
\special{pa  800 800}%
\special{pa  800 500}%
\special{fp}%
\special{pn 8}%
\special{pa 1100 500}%
\special{pa 1400 500}%
\special{pa 1400 800}%
\special{pa 1100 800}%
\special{pa 1100 500}%
\special{fp}%
\special{pn 8}%
\special{pa 1400 500}%
\special{pa 1700 500}%
\special{pa 1700 800}%
\special{pa 1400 800}%
\special{pa 1400 500}%
\special{fp}%
\special{pn 8}%
\special{pa 500 200}%
\special{pa 800 200}%
\special{pa 800 500}%
\special{pa 500 500}%
\special{pa 500 200}%
\special{dt 0.03}%
\special{pn 8}%
\special{pa  800 200}%
\special{pa 1100 200}%
\special{pa 1100 500}%
\special{pa  800 500}%
\special{pa  800 200}%
\special{dt 0.03}%
\special{pn 8}%
\special{pa 1100 200}%
\special{pa 1400 200}%
\special{pa 1400 500}%
\special{pa 1100 500}%
\special{pa 1100 200}%
\special{dt 0.03}%
\special{pn 8}%
\special{pa 1400 200}%
\special{pa 1700 200}%
\special{pa 1700 500}%
\special{pa 1400 500}%
\special{pa 1400 200}%
\special{fp}%
\special{pn 8}%
\special{pa 1700 200}%
\special{pa 2000 200}%
\special{pa 2000 500}%
\special{pa 1700 500}%
\special{pa 1700 200}%
\special{fp}%
\special{pn 8}%
\special{pa 2000 200}%
\special{pa 2300 200}%
\special{pa 2300 500}%
\special{pa 2000 500}%
\special{pa 2000 200}%
\special{fp}%
\special{pn 8}%
\special{pa 2000 200}%
\special{pa 2300 200}%
\special{pa 2300 500}%
\special{pa 2000 500}%
\special{pa 2000 200}%
\special{fp}%
\special{pn 8}%
\special{pa 2300 200}%
\special{pa 2600 200}%
\special{pa 2600 500}%
\special{pa 2300 500}%
\special{pa 2300 200}%
\special{fp}%
%
\put( 2.9,-19.2){\makebox(0,0)[lb]{$q_2^5$}}%
\put( 2.9,-16.2){\makebox(0,0)[lb]{$q_2^4$}}%
\put( 2.9,-13.2){\makebox(0,0)[lb]{$q_2^3$}}%
\put( 2.9,-10.2){\makebox(0,0)[lb]{$q_2^2$}}%
\put( 5.2,-10.2){\makebox(0,0)[lb]{$q_1q_2^2$}}%
\put( 5.2, -7.2){\makebox(0,0)[lb]{$q_1q_2^1$}}%
\put( 8.2, -7.2){\makebox(0,0)[lb]{$q_1^2q_2^1$}}%
\put(11.2, -7.2){\makebox(0,0)[lb]{$q_1^3q_2^1$}}%
\put(14.2, -7.2){\makebox(0,0)[lb]{$q_1^4q_2^1$}}%
\put(14.9, -4.2){\makebox(0,0)[lb]{$q_1^4$}}%
\put(17.9, -4.2){\makebox(0,0)[lb]{$q_1^5$}}%
\put(20.9, -4.2){\makebox(0,0)[lb]{$q_1^6$}}%
\put(23.9, -4.2){\makebox(0,0)[lb]{$q_1^7$}}%
\end{picture}%
\end{center}
\caption{The assignment of the monomials for $(4,4,2,1,1,1)$}
\label{fig:monom}
\end{figure}
\end{rmk}

\begin{proof}[{Proof of Lemma \ref{lem:wronski}}]
(1)
By the definition of the function $\ep_r$ and Remark \ref{rmk:snake}, we find that there are two non-vanishing terms occurring in $\psi_\mu\circ\varphi_\mu(S_n)$ as
\begin{align}\label{eq:wronski:tmp1}
\begin{split}
\psi_\mu\circ\varphi_\mu(S_n)
&=\binom{n}{l-1}(q_1^{n-l+1}-q_2^{l-1})
\left(\dfrac{q_2-1}{q_1-1}\right)^{n-l+1} 
(\ep_{l-1}(q_1)*\ep_{n-l+1}(q_2))(X,Y,z)\\
&+
\binom{n}{l}(q_1^{n-l}-q_2^{l})
\left(\dfrac{q_2-1}{q_1-1}\right)^{n-l} 
(\ep_{l}(q_1)*\ep_{n-l}(q_2))(X,Y,z),
\end{split}
\end{align}
where we used the notation
\begin{align*}
X
&\seteq
\{q_1^{(\mu_l-1)+\cdots+(\mu_{i+1}-1)}q_2^{i-1}\mid 1\le i\le l-1\},\\
Y
&\seteq
\{q_1^{(\mu_l-1)+\cdots+(\mu_{i+1}-1)+j}q_2^{i-1}
\mid 1\le i\le l-1,1\le j\le\mu_{i-1}\},\\
z&\seteq q_2^{l-1}.
\end{align*}
(See Remark \ref{rmk:snake2} as to the diagrammatic explanation.) As for the second term in (\ref{eq:wronski:tmp1}), a direct computation gives us 
\begin{align*}
&\binom{n}{l}
(q_1^{n-l}-q_2^{l})
\left(\dfrac{q_2-1}{q_1-1}\right)^{n-l} 
(\ep_{l}(q_1)*\ep_{n-l}(q_2))(X,Y,z)\\
&=(q_1^{n-l}-q_2^{l})
\left(\dfrac{q_2-1}{q_1-1}\right)^{n-l}
\ep_l(X,z;q_1)\,\ep_{n-l}(Y;q_2)
\prod_{x\in X\cup\{z\},\, y\in Y}\omega(x,y).
\end{align*}
A similar result can be derived for the first term in (\ref{eq:wronski:tmp1}). 
Thus we only need to show the equation
\begin{align}\label{eq:wronski:concl}
\begin{split}
(q_1^{n-l+1}-q_2^{l-1})\left(\dfrac{q_2-1}{q_1-1}\right)
&\prod_{y\in Y}\ep_2(y,z;q_2)\prod_{x\in X}\omega(x,z)
\\&+
(q_1^{n-l}-q_2^{l})\prod_{x\in X}\ep_2(x,z;q_1)
\prod_{y\in Y}\omega(z,y)
=0.
\end{split}
\end{align}
A straightforward calculation shows that this is correct.
\\
(2)
By the snake evaluation, $y_j$ in $\zeta_\mu$ is replaced by $q_2^{j-1}q_1^{(\mu_l-1)+\cdots+(\mu_{j+1}-1)}$. We only need to check that each factor in $\zeta_\mu$ 
does not become
zero after the snake evaluation. 
Checking is elementary, and we omit it. 
\end{proof}

\begin{rmk}\label{rmk:snake2}
The variables in $X$ are located at the leftmost boxes of the rows in the border diagram. The example $\lambda=(4,4,2,1,1)$ is given in Figure \ref{fig:colored}.

\begin{figure}[htbp]
\unitlength 0.1in
\begin{center}
\begin{picture}( 0.0000, 12.0000)(  18.0000,-26.0000)
\special{pn 8}%
\special{sh 0.900}%
\special{pa 600 2600}%
\special{pa 800 2600}%
\special{pa 800 2400}%
\special{pa 600 2400}%
\special{pa 600 2600}%
\special{fp}%
\special{pn 8}%
\special{sh 0.300}%
\special{pa 800 2400}%
\special{pa 600 2400}%
\special{pa 600 2200}%
\special{pa 800 2200}%
\special{pa 800 2400}%
\special{fp}%
\special{pn 8}%
\special{sh 0.300}%
\special{pa 600 2200}%
\special{pa 800 2200}%
\special{pa 800 2000}%
\special{pa 600 2000}%
\special{pa 600 2200}%
\special{fp}%
\special{pn 8}%
\special{pa 800 2000}%
\special{pa 1000 2000}%
\special{pa 1000 1800}%
\special{pa 800 1800}%
\special{pa 800 2000}%
\special{fp}%
\special{pn 8}%
\special{sh 0.300}%
\special{pa 600 2000}%
\special{pa 800 2000}%
\special{pa 800 1800}%
\special{pa 600 1800}%
\special{pa 600 2000}%
\special{fp}%
\special{pn 8}%
\special{sh 0.300}%
\special{pa 800 1800}%
\special{pa 1000 1800}%
\special{pa 1000 1600}%
\special{pa 800 1600}%
\special{pa 800 1800}%
\special{fp}%
\special{pn 8}%
\special{pa 1000 1800}%
\special{pa 1200 1800}%
\special{pa 1200 1600}%
\special{pa 1000 1600}%
\special{pa 1000 1800}%
\special{fp}%
\special{pn 8}%
\special{pa 1200 1800}%
\special{pa 1400 1800}%
\special{pa 1400 1600}%
\special{pa 1200 1600}%
\special{pa 1200 1800}%
\special{fp}%
\special{pn 8}%
\special{pa 1400 1800}%
\special{pa 1600 1800}%
\special{pa 1600 1600}%
\special{pa 1400 1600}%
\special{pa 1400 1800}%
\special{fp}%
\special{pn 8}%
\special{pa 1200 1800}%
\special{pa 1400 1800}%
\special{pa 1400 1600}%
\special{pa 1200 1600}%
\special{pa 1200 1800}%
\special{fp}%
\special{pn 8}%
\special{sh 0.300}%
\special{pa 1400 1600}%
\special{pa 1600 1600}%
\special{pa 1600 1400}%
\special{pa 1400 1400}%
\special{pa 1400 1600}%
\special{fp}%
\special{pn 8}%
\special{pa 1600 1600}%
\special{pa 1800 1600}%
\special{pa 1800 1400}%
\special{pa 1600 1400}%
\special{pa 1600 1600}%
\special{fp}%
\special{pn 8}%
\special{sh 0.300}%
\special{pa 1400 1600}%
\special{pa 1600 1600}%
\special{pa 1600 1400}%
\special{pa 1400 1400}%
\special{pa 1400 1600}%
\special{fp}%
\special{pn 8}%
\special{pa 1800 1600}%
\special{pa 2000 1600}%
\special{pa 2000 1400}%
\special{pa 1800 1400}%
\special{pa 1800 1600}%
\special{fp}%
\special{pn 8}%
\special{pa 2000 1600}%
\special{pa 2200 1600}%
\special{pa 2200 1400}%
\special{pa 2000 1400}%
\special{pa 2000 1600}%
\special{fp}%
\special{pn 8}%
\special{pa 600 1600}%
\special{pa 800 1600}%
\special{pa 800 1400}%
\special{pa 600 1400}%
\special{pa 600 1600}%
\special{dt 0.03}%
\special{pn 8}%
\special{pa  800 1600}%
\special{pa 1000 1600}%
\special{pa 1000 1400}%
\special{pa  800 1400}%
\special{pa  800 1600}%
\special{dt 0.03}%
\special{pn 8}%
\special{pa 1000 1600}%
\special{pa 1200 1600}%
\special{pa 1200 1400}%
\special{pa 1000 1400}%
\special{pa 1000 1600}%
\special{dt 0.03}%
\special{pn 8}%
\special{pa 1200 1600}%
\special{pa 1400 1600}%
\special{pa 1400 1400}%
\special{pa 1200 1400}%
\special{pa 1200 1600}%
\special{dt 0.03}%
\special{pn 8}%
\special{pa 600 1800}%
\special{pa 800 1800}%
\special{pa 800 1600}%
\special{pa 600 1600}%
\special{pa 600 1800}%
\special{dt 0.03}%
\special{pn 8}%
\special{pa 1200 1800}%
\special{pa 1400 1800}%
\special{pa 1400 1600}%
\special{pa 1200 1600}%
\special{pa 1200 1800}%
\special{fp}%
\special{pn 8}%
\special{sh 0.300}%
\special{pa 2600 1700}%
\special{pa 2800 1700}%
\special{pa 2800 1500}%
\special{pa 2600 1500}%
\special{pa 2600 1700}%
\special{fp}%
\put(29.0000,-16.5000){\makebox(0,0)[lb]{$:\,X$}}%
\special{pn 8}%
\special{pa 2600 2000}%
\special{pa 2800 2000}%
\special{pa 2800 1800}%
\special{pa 2600 1800}%
\special{pa 2600 2000}%
\special{fp}%
\put(29.0000,-19.5000){\makebox(0,0)[lb]{$:\,Y$}}%
\special{pn 8}%
\special{sh 0.900}%
\special{pa 2600 2300}%
\special{pa 2800 2300}%
\special{pa 2800 2100}%
\special{pa 2600 2100}%
\special{pa 2600 2300}%
\special{fp}%
\put(29.0000,-22.5000){\makebox(0,0)[lb]{$:\,z$}}%
\end{picture}%
\end{center}
\caption{Snake evaluation for $(4,4,2,1,1,1)$}
\label{fig:colored}
\end{figure}
\end{rmk}

Recall we have the following relation between the
elementary symmetric functions $e_n(x)$'s and the $g_n(x;q,t)$'s.

\begin{lem}\label{lem:qt_wronski}
For $n\in\bbN_+$ we have
\begin{align}\label{eq:Wronski:ply}
\sum_{k=0}^n(-1)^k(1-q^kt^{n-k})e_{n-k}(x)g_k(x;q,t)=0.
\end{align}
\end{lem}
\begin{proof}
Introduce generating functions as
\begin{align}
&G(y)\seteq\sum_{n\ge0}g_n(x;q,t)y^n=\exp\biggl(\sum_{n\ge1}\dfrac{1}{n}\dfrac{1-t^n}{1-q^n}p_n(x)y^n\biggr)=
\prod_{i\ge 1}
   \dfrac{\bigl(t x_i y;q\bigr)_\infty}{\bigl(x_i y;q\bigr)_\infty},
\label{eq:G(y)}\\
&E(y)\seteq\sum_{n\ge0}e_n(x;q,t)y^n
=\omega_{q,t}(G(y))
=\exp\biggl(\sum_{n\ge1}\dfrac{(-1)^{n-1}}{n}p_n(x)y^n\biggr)=
\prod_{i\ge 1}(1+x_i y)
.
\label{eq:E(y)}
\end{align}
Hence we have $E(y)G(-y)=E(t y)G(-q y)$, indicating the conclusion. 
\end{proof}

Remark that by using (\ref{eq:Wronski:ply}), we can uniquely 
express $g_r$ as a polynomial in $e_s$'s and vice versa.

\begin{prop}\label{prop:eigen:G_n}
We have
\begin{align*}
\widehat{G}_r(q,t)P_\lambda(x;q,t)=g_r(s^\lambda;q,t)P_\lambda(x;q,t).
\end{align*}
\end{prop}
\begin{proof}
Using the Wronski relation for the operators (\ref{eq:Wronski:op}),
we can uniquely express $\widehat{G}_r$ as a polynomial in
$\widehat{E}_s$'s, just in the same manner for $g_r$. 
Hence the eigenvalues of $\widehat{G}_r$ must be $g_r(s^\lambda)$.
\end{proof}

\begin{rmk}
We can introduce a difference operator $G_r^{(n)}$ acting on $\Lambda_{n,\bbF}$, 
and define the limit  $G_r:=\varprojlim_n G_r^{(n)}$
as was done for $E_r$ in Proposition \ref{prop:E_r}. 
Then we can show that the free field realization of $G_r$ 
is $\widehat{G}_r$. This will be treated in \S\S\,\ref{subsec:appendix_yoko}
\end{rmk}

\subsection{Triangularity and the proof of Theorems \ref{thm:2}\ and \ref{thm:3}}
\label{subsec:intersect}

\begin{lem}\label{lem:mcd:trng}
The Macdonald symmetric function $P_\lambda(x;q,t)$ has a triangular expansion with respect to the bases $(e_\lambda)$ and $(g_\lambda)$. Precisely, it can be written as
\begin{align*}
P_\lambda(x;q,t)=\sum_{\mu\ge\lambda}c_{\lambda\mu}^{g\to P}(q,t)g_\mu(x;q,t),
\quad
P_\lambda(x;q,t)=\sum_{\mu\ge\lambda'}c_{\lambda\mu}^{e\to P}(q,t)e_\mu(x),
\end{align*}
with $c_{\lambda\mu}^{g\to P}$, $c_{\lambda\mu}^{e\to P}\in\bbF$ and $c_{\lambda\lambda}^{g\to P}\neq0$, $c_{\lambda\lambda}^{e\to P}\neq0$.
\end{lem}

\begin{proof}
As to the first statement, see \cite[Lemma 8.1]{S:2006}. For the second one, we apply the automorphism $\omega_{q,t}$ (\ref{eq:omega_auto}) to the expansion $P_\lambda(x;q,t)=\sum_{\mu\ge\lambda}c_{\lambda\mu}^{g\to P}(q,t)g_\mu(x;q,t)$. By (\ref{eq:dual}), we obtain
\begin{align*}
Q_{\lambda'}(x;t,q)=\sum_{\mu\ge\lambda}c_{\lambda\mu}^{g\to P}(q,t)e_\mu(x).
\end{align*}
Rewriting $Q_{\lambda'}=b_{\lambda'}P_{\lambda'}$ and changing the notation $\mu$ and $\lambda$, we have the equation
\begin{align*}
P_{\lambda}(x;t,q)=
\sum_{\mu\ge\lambda'}b_\lambda(t,q)^{-1}c_{\lambda'\mu}^{g\to P}(q,t)e_\mu(x).
\end{align*}
Hence the second conclusion holds if one sets $c_{\lambda\mu}^{e\to P}(q,t) \seteq b_\lambda(q,t)^{-1}c_{\lambda'\mu}^{g\to P}(t,q)$.
\end{proof}

\begin{proof}[{Proof of Theorems \ref{thm:2} and \ref{thm:3}}]
First we will construct a non-zero element $f_\lambda$ in $\calA_{\lambda}^{(q^{-1})}\bigcap\calA_{{\lambda'}}^{(t)}$. By Lemma \ref{lem:mcd:trng}, we can write
\begin{align}\label{eq:tmp:1}
P_\lambda(x;q,t)
=\sum_{\mu\ge\lambda'}c_{\lambda\mu}^{e\to P}(q,t)e_\mu(x)
=\sum_{\mu\ge\lambda}c_{\lambda\mu}^{g\to P}(q,t)g_\mu(x;q,t).
\end{align}
Then we define the elements $f^{(q^{-1})}_\lambda(z),f^{(t)}_\lambda(z)\in \calA_n$ by
\begin{align}\label{eq:fo:trng}
\begin{split}
&f^{(q^{-1})}_\lambda(z)
\seteq\sum_{\mu\ge\lambda'}c_{\lambda\mu}^{e\to P}(q,t)\ep_{\mu}(z;q)
\prod_{i\geq 1} \dfrac{t^{-\mu_i(\mu_i+1)/2}}{(t^{-1};t^{-1})_{\mu_i}}
\dfrac{[\mu_i]_{t^{-1}}!}{\mu_i!},\\
&f^{(t)}_\lambda
\seteq
\sum_{\mu\ge\lambda}c_{\lambda\mu}^{g\to P}(q,t)\ep_{\mu}(z;t)
\prod_{i\geq 1}
\dfrac{(-1)^{\mu_i} q^{\mu_i(\mu_i-1)}}{(q;q)_{\mu_i}}
\dfrac{[\mu_i]_q!}{\mu_i!}.
\end{split}
\end{align}
By Proposition \ref{prop:wbcalA_n:poincare} (2) we find that
\begin{align}\label{eq:tmp:3}
f^{(q^{-1})}_\lambda\in\calA_\lambda^{(q^{-1})},\quad 
f^{(t)}_\lambda\in\calA_{\lambda'}^{(t)}.
\end{align}
On the other hand, Propositions \ref{prop:eigen:E_r}, \ref{prop:eigen:G_n} and (\ref{eq:tmp:1}) yield
\begin{align*}
\calO(f_\lambda^{(q^{-1})})P_\mu=\calO(f_\lambda^{(t)})P_\mu=P_\lambda(s^\lambda;q,t)P_\mu
\end{align*}
for any partition $\mu$. Since the family $(P_\mu)$ is a basis of $\calF\cong\Lambda_\bbF$, we should have $\calO(f_\lambda^{(q^{-1})})=\calO(f_\lambda^{(t)})$. Then by Proposition \ref{prop:spectral_criterion}, $f_\lambda^{(q^{-1})}=f_\lambda^{(t)}$ holds. Thus we can define 
\begin{align*}
f_\lambda\seteq f_\lambda^{(q^{-1})}=f_\lambda^{(t)}.
\end{align*}
Then from (\ref{eq:tmp:3}) we have $f_\lambda\in \calA_\lambda^{(q^{-1})}\cap\calA_{\lambda'}^{(t)}$.

Now we shall prove $\dim(\calA_\lambda^{(q^{-1})}\cap\calA_{\lambda'}^{(t)})=1$. Suppose that $f_\lambda'$ is an element of $\calA_\lambda^{(q^{-1})}\cap\calA_{\lambda'}^{(t)}$. By Proposition \ref{prop:wbcalA_n:poincare} we have
\begin{align*}
f_\lambda'
\in\bbF[\ep_\mu(z;q^{-1})\mid\mu\ge\lambda']\bigcap
\bbF[\ep_\mu(z;t)\mid\mu\ge\lambda].
\end{align*}
Then the triangular decomposition (\ref{eq:fo:trng}) yields
\begin{align*}
f_\lambda'
\in\bbF[f_\mu\mid\mu\le\lambda]\bigcap
\bbF[f_\mu\mid\mu\ge\lambda]
\end{align*}
Therefore $f_\lambda'$ must be proportional to $f_\lambda$.
\end{proof}

We also note the following remark.

\begin{cor}
The family $(f_\lambda)_{\lambda\,\vdash\,n}$ forms a basis of $\calA_n$.
\end{cor}

\subsection{Vertex operator $\xi(z)$ and the Ding-Iohara algebra $\calU(q,t)$}
\label{subsec:DI}
Before closing this section, we continue our study on the vertex operator $\eta(z)$ by applying a quantum group technique. 
Our goal in this subsection is to identify $\eta(z)$ 
as the level one representation of the Ding-Iohara algebra $\calU(q,t)$.
(See Appendix \ref{app:DI}.)

Set 
\begin{align}\label{eq:xi}
&\xi(z):=
\exp\Bigg(-\sum_{n>0} \dfrac{1-t^{-n}}{n}(t/q)^{n/2}a_{-n} z^{n}\Bigg)
\exp\Bigg(\sum_{n>0}  \dfrac{1-t^{n}}{n} (t/q)^{n/2} a_n z^{-n}\Bigg),
\end{align}
and consider its Fourier expansion $\xi(z)=\sum_{n\in \bbZ}\xi_n z^{-n}$.

\begin{lem}
We have
\begin{align}
&\xi(z)\xi(w)
 =\dfrac{(1-w/z)(1-q^{-1}t w/z)}{(1-q^{-1}w/z)(1-t w/z)}:\xi(z)\xi(w):,
\label{eq:xi_xi}\\
&\eta(z)\xi(w)
 =\dfrac{(1-q^{-1/2}t^{-1/2}w/z)(1-q^{1/2}t^{1/2} w/z)}
        {(1-q^{-1/2}t^{1/2} w/z)(1-q^{1/2}t^{-1/2}w/z)}
  :\eta(z)\xi(w):,
\label{eq:eta_xi}\\
&\xi(w)\eta(z)
 =\dfrac{(1-q^{-1/2}t^{-1/2}z/w)(1-q^{1/2}t^{1/2} z/w)}
        {(1-q^{-1/2}t^{1/2} z/w)(1-q^{1/2}t^{-1/2}z/w)}
  :\xi(w)\eta(z):,
\label{eq:xi_eta}\\
&\xi(z)\phi(y)\cdot 1
 =\dfrac{1-t^{3/2}q^{-1/2}y/z}
        {1-t^{1/2}q^{-1/2}y/z}
  :\phi(y)\xi(z):\cdot\, 1,
\nonumber
\\
&
 :\xi(y)\phi({q^{1/2}t^{-1/2}}y):\cdot\, 1
  =T_{q^{-1},y} \phi(q^{1/2}t^{-1/2}y)\cdot 1.
\nonumber
\end{align}

\end{lem}

The equation (\ref{eq:xi_xi}) gives the operator-valued symmetric Laurent series
\begin{align*}
&\dfrac{1}{\omega(z_1,z_2;q,t^{-1},q^{-1}t)}
\xi(z_1)\xi(z_2)\\
&=
\dfrac{q^{-1} t(1-z_2/z_1)^2 }
{(1-q^{-1} z_2/z_1)(1-t z_2/z_1)
}\dfrac{ (1-z_1/z_2)^2}
{(1-q^{-1} z_1/z_2)(1-t z_1/z_2)
}
:\xi(z_2)\xi(z_1):.
\end{align*}
Hence we can define the commutative ring 
$\calM'\seteq\{\calO(f;q^{-1},t^{-1}) \mid f\in \calA\}$
by the same argument we had for $\calM$.
It is clearly seen from the formula
\begin{align*}
\xi_0 \phi(y_1)\cdots \phi(y_n)\cdot 1
= (1-t) t^{n-1} D_{n,y}^1(q^{-1},t^{-1}) \phi(y_1)\cdots \phi(y_n)\cdot 1+
t^n  \phi(y_1)\cdots \phi(y_n)\cdot 1,
\end{align*}
that the vertex operator $\xi(z)$ provides us with the free field realization of the Macdonald operators of type $D_n^s(q^{-1},t^{-1})$. Recall the property of the Macdonald symmetric function $P_\lambda(x;q,t)=P_\lambda(x;q^{-1},t^{-1})$ \cite[(VI.4.13) (iv)]{M:1995:book}. This indicates the commutativity of the Macdonald difference operators $[D_n^r(q,t),D_n^s(q^{-1},t^{-1})]=0$.

\begin{prop}\label{prop:M_M'}
Let $f\in \calM$ and $g\in \calM'$.
Then $[f,g]=0$.
\end{prop}
\begin{proof}
{}From (\ref{eq:eta_xi}) and (\ref{eq:xi_eta}), we have
\begin{align*}
&[\eta(z),\xi(w)]=
\dfrac{(1-q)(1-t^{-1})}{1-q t^{-1}}
\Bigg(\delta((t/q)^{-1/2}z/w)\varphi^+((t/q)^{1/4}w)-
\delta((t/q)^{1/2}z/w)\varphi^-((t/q)^{-1/4}w)
 \Bigg),
\end{align*}
where $\delta(z)\seteq\sum_{n\in\bbZ}z^n$ is the formal delta function and
\begin{align}
&\varphi^+(z):=
\exp\Bigg(-\sum_{n>0} \dfrac{1-t^{n}}{n} (1-t^n q^{-n})(t/q)^{-n/4} a_n z^{-n}\Bigg)=
\sum_{n\in \bbN} \varphi^+_n z^{-n},
\label{eq:phi^+}
\\
&\varphi^-(z):=
\exp\Bigg(+\sum_{n>0} \dfrac{1-t^{-n}}{n} (1-t^n q^{-n})(t/q)^{-n/4} a_{-n} z^{n}\Bigg)=
\sum_{n\in \bbN} \varphi^-_{-n} z^{n}.
\label{eq:phi^-}
\end{align}
Then $[\eta_0,\xi_0]=\dfrac{(1-q)(1-t^{-1})}{1-q t^{-1}}
(\varphi^+_0-\varphi^-_0)=0$ holds since 
$\varphi^+_0=\varphi^-_0=1$. Thus $\xi_0$ is diagonalized  by the basis $(P_\lambda(x;q,t))$ since the eigenvalues are distinct.
Hence all the operators in $\calM'$ are simultaneously 
diagonalized by the same basis.
\end{proof}

Hereby we have obtained four vertex operators $\eta(z),\xi(z),\varphi^+(z)$ and
$\varphi^-(z)$. One finds that they satisfy the relations 
of the Drinfeld generators of the Ding-Iohara algebra.
Recall that the Ding-Iohara algebra
gives us a generalization of the quantum affine algebra
which has a function $g(z)$ depending on  {\it deformation parameters}
whenever we have $g(z)=g(z^{-1})^{-1}$. For our purpose 
we need to set 
\begin{align*}
g(z)\seteq
\dfrac{(1-q z)(1-t^{-1}z)(1-q^{-1} t z)}{(1-q^{-1}z)(1-t z)(1-q t^{-1}z)},
\end{align*}
indicating that the corresponding 
Ding-Iohara algebra has two parameters $q$ and $t$.
We denote this particular case of the Ding-Iohara algebra by $\calU(q,t)$.
For the sake of the readers' convenience,  
the definition and the basic properties of the Ding and Iohara
algebra, including the Drinfeld coproduct $\Delta$, the antipode $a$ and 
the counit $\varepsilon$, are summarized 
in the beginning of Appendix \ref{app:DI}.

One can construct two kinds of representations, 
one acting on the space of Laurent polynomials $V_x=\bbQ(q^{1/2},t^{1/2})[x,x^{-1}]$ and the other on the Fock space $\calF$. 
We denote the former by $\pi_x(\cdot)$ and the latter by $\rho(\cdot)$.
Then one can introduce the intertwining operator 
$\Phi^{\calF}_{V_x\otimes \calF}:V_x\otimes \calF\rightarrow \calF$
determined by the property $\Phi^{\calF}_{V_x\otimes \calF} \Delta(a)=
a\Phi^{\calF}_{V_x\otimes \calF}$ for any $a\in \calU$.
As for the detail, 
see Proposition \ref{prop:intertwiner}. 

It is tempting to study the 
meaning of  $\Phi^{\calF}_{V_x\otimes \calF}$ 
in terms of the Macdonald theory.
Set
\begin{align}
\label{eq:base_func}
\widetilde{\Phi}(y)\seteq
\exp
\Bigg(
 \sum_{n>0}\dfrac{1}{n}\dfrac{1-t^n}{1-q^{n}}t^{-n}a_{-n} y^{n}
\Bigg)
\exp\Bigg(
 -\sum_{n>0}\dfrac{1}{n}\dfrac{1-t^n}{1-q^{n}}q^{n}a_{n} y^{-n}
\Bigg).
\end{align}
From Proposition \ref{prop:intertwiner}, 
$\Phi(y):=\widetilde{\Phi}(q^{1/2}y)$ is nothing but the generating function 
of the intertwining operator  $\Phi^{\calF}_{V_x\otimes \calF}$.
The shift $q^{1/2}$ is introduced just for simplicity of display.
An easy calculation gives us the following.

\begin{prop}
\label{prop:recover_macd}
We have
\begin{align*}
&\eta(z)\widetilde{\Phi}(y_1)\cdots\widetilde{\Phi}(y_n)
=
t^{-n} \prod_{i=1}^n 
 \dfrac{(1-q t^{-1}z /y_i)(1-t z/y_i)}{(1-q z/y_i)(1-z/y_i)}
\widetilde{\Phi}(y_1)\cdots\widetilde{\Phi}(y_n)\eta(z)
\\
&\hskip 3em
 +(1-t^{-1})
 \sum_{i=1}^n t^{-i+1} 
 \prod_{j=1}^{i-1} \dfrac{(1-q t^{-1}z /y_j)(1-t z/y_j)}{(1-q z/y_j)(1-z/y_j)}
 \delta(y_i/z) T_{q,y_i} \widetilde{\Phi}(y_1)\cdots\widetilde{\Phi}(y_n),
\\
&\widetilde{\Phi}(y_1)\cdots\widetilde{\Phi}(y_n)
=\prod_{1\leq i<j\leq n} 
  \dfrac{(q t^{-1}y_j/y_i;q)_\infty}{(q y_j/y_i;q)_\infty }
 :\widetilde{\Phi}(y_1)\cdots\widetilde{\Phi}(y_n):,
\\
&:\widetilde{\Phi}(y_1)\cdots\widetilde{\Phi}(y_n):\cdot\, 1
=\phi(t^{-1}y_1)\cdots\pi(t^{-1}y_n)\cdot 1.
\end{align*}
\end{prop}

{}From this we recover 
\begin{align*}
\eta_0 \phi(y_1)\cdots \phi(y_n)\cdot 1=
t^{-n }\phi(y_1)\cdots \phi(y_n)\cdot 1+
(1-t^{-1})t^{-n+1} D_{n,y}^1 \phi(y_1)\cdots \phi(y_n)\cdot 1,
\end{align*}
which is the simplest equation of our construction in this section.

\subsection{Appendix: Remarks on the calculation in this section}
\label{subsec:appendix_yoko}

\subsubsection{Direct calculation of $G_r$}
We use the following proposition due to M. Noumi \cite{N:private}.

\begin{prop}[Noumi]\label{prop:Noumi}
Set
\begin{align}
H_n^l &\seteq 
 \sum_{\substack{\nu \in \mathbb{N}^n \\ |\nu|=l}} 
 \biggl(\prod_{1\le i < j\le n} \frac{q^{\nu_i}x_i - q^{\nu_j}x_j}{x_i-x_j}
 \biggr)
 \biggl(\prod_{i,j=1}^n \frac{(t x_i/x_j ; q)_{\nu_i}}{(q x_i/x_j;q)_{\nu_i}}
 \biggr)
 \prod_{i=1}^n T_{q,x_i}^{\nu_i}.\label{eq:Noumi}
\end{align}
Then the operator $H_n^l$'s act on $\Lambda_{n,\bbF}$ and $H_n^l\in \bbF[D_n^1,\ldots,D_n^n]$. Moreover we have
\begin{align*}
H_n^l P_\lambda(x;q,t)=g_l^{(n)}(t^n s^\lambda_1,\ldots,t^n s^\lambda_n;q,t)
P_\lambda(x;q,t)
\end{align*}
for any partition $\lambda$ such that $\ell(\lambda)\leq n$, 
where $g_l^{(n)}(s_1,\ldots,s_n;q,t)$ is defined by the generating function
\begin{align*}
\sum_{l\ge0}g_l^{(n)}(s_1,\ldots,s_n;q,t)y^l=
\exp\biggl(\sum_{m\ge1}\dfrac{1}{m}\dfrac{1-t^m}{1-q^m}\sum_{i=1}^n s_i^m y^m\biggr),
\end{align*}
and $s^\lambda_i\seteq t^{-i}q^{\lambda_i}$ for $i=1,\ldots,n$.
\end{prop}

\begin{prop}\label{prop:G_r}
Set
\begin{align*}
G_r^{(n)} \seteq 
 \dfrac{t^{-r n} q^{\binom{r}{2}}}{(-1)^r(q;q)_r}
 \sum_{l=0}^r (-1)^l q^{-\binom{l}{2}}q^{-l(r-l)}(q^{r-l+1};q)_l H_n^l.
\end{align*}
For any partition $\lambda$ such that $\ell(\lambda)\leq n$, 
we have
\begin{align*}
G_r^{(n)}P_\lambda(x;q,t)=g_r(s^\lambda) P_\lambda(x;q,t),
\end{align*}
where $P_\lambda(x;q,t)\in \Lambda_{n,\bbF}$ denotes the 
Macdonald symmetric {\it polynomial},
and $g_r(x)\in \Lambda_{\bbF}$ is the symmetric {\it function} 
defined in (\ref{eq:G(y)}), 
and $s^\lambda=(t^{-1}q^{\lambda_1},t^{-2}q^{\lambda_2},\ldots)$.
Hence the 
inductive limit $G_r\seteq\varprojlim_{n}G_r^{(n)}$ exists.
The operator $\widehat{G}_r$ given in (\ref{eq:G_n})
is a free field representation of $G_r$, namely we have
$\widehat{G}_n P_\lambda(x;q,t)=G_r P_\lambda(x;q,t)$ for any partition $\lambda$.
\end{prop}

\begin{proof}
As in the proof of Proposition \ref{prop:eigen:E_r}, we compute $\widehat{G}_r \phi_y \cdot 1$, where $\phi_y\seteq \phi(y_1)\cdots\phi(y_n)$. Noting that
\begin{align*}
\dfrac{[r]_{q}!}{r!}
\prod_{1\le i<j\le r}\ep_2(z_i,z_j;q)^{-1}
=\Sym\biggl(\prod_{1\le i<j\le r}\frac{1-z_j/z_i}{1-q z_j/z_i} \biggr),
\end{align*}
one finds that 
\begin{align*}
\widehat{G}_r \phi_y \cdot 1
&= \dfrac{(-1)^r q^{\binom{r}{2}}}{(q;q)_r}
   \biggl[
    \prod_{1 \le i < j \le r}\frac{1-z_j/z_i}{1-q z_j/z_i}
    : \eta(z_1) \cdots \eta(z_r) : \phi_y \cdot 1 
   \biggl]_1
\\
&= \dfrac{(-1)^r q^{\binom{r}{2}}}{(q;q)_r}
   \biggl[ 
    \biggl(\prod_{j=1}^r \prod_{i=1}^n \frac{1-y_i/z_j}{1-t y_i/z_j}\biggr)
    \biggl(\prod_{1 \le i < j \le r}\frac{1-z_j/z_i}{1-q z_j/z_i}\biggr)
    \phi_y : \eta(z_1) \cdots \eta(z_r) : \cdot\, 1
   \biggl]_1
\end{align*}
Extending the last integral, we introduce 
\begin{align*}
\overline{G}_{r}^{(\mu)}
\seteq
\biggl[ 
 \biggl(\prod_{j=1}^r \prod_{i=1}^n \frac{1-y_i/q^{\mu_i}z_j}{1-t y_i/z_j}
 \biggr)
 \biggl(\prod_{1 \le i < j \le r}\frac{1-z_j/z_i}{1-q z_j/z_i} \biggr)
 \phi_y : \eta(z_1) \cdots \eta(z_r) : \cdot\, 1
\biggl]_1
\end{align*}
with $\mu=(\mu_1,\ldots,\mu_n) \in \mathbb{N}^n$. We also generalize the operator $H_n^l$ to
\begin{align*}
H_n^{l,(\mu)} \seteq 
 \sum_{\substack{\nu \in \mathbb{N}^n \\ |\nu|=l}} 
 \biggl(\prod_{1\le i<j \le n}\frac{q^{\nu_i}y_i - q^{\nu_j}y_j}{y_i-y_j}
 \biggr) 
 \biggl(\prod_{i,j=1}^n 
  \frac{(t q^{\mu_j}y_i/y_j ; q)_{\nu_i}}{(q y_i/y_j;q)_{\nu_i}} 
 \biggr)
 \prod_{i=1}^n T_{q,y_i}^{\nu_i}.
\end{align*}
Here we denoted $|\nu|\seteq\nu_1+\cdots+\nu_n$ for $\nu\in\bbN^n$. 
Then the conclusion corresponds to the special case $\mu=(0,0,\ldots,0)$ of the next proposition.
\end{proof}

\begin{prop}\label{prop:yoko_para}
We have 
\begin{align*}
\overline{G}_{r}^{(\mu)} = 
 \dfrac{t^{-r n}}{q^{r|\mu|}}
 \sum_{l=0}^r (-1)^l q^{-\binom{l}{2}}q^{-l(r-l)}(q^{r-l+1};q)_l 
 H_n^{l,(\mu)} \phi_y \cdot 1.
\end{align*}
\end{prop}
\begin{proof}
Executing the integral, we have
\begin{align*}
\begin{split}
\overline{G}_r^{(\mu)}
=& \sum_{k=1}^n 
    \bigl(1-\dfrac{1}{q^{\mu_k}t}\bigr)
    \prod_{i \neq k}\dfrac{1-y_i/q^{\mu_i}t y_k}{1-y_i/y_k} T_{q,y_k}
    \overline{G}_{r-1}^{(\mu+\ve_k)}  
   +\dfrac{t^{-n}}{q^{|\mu|}} \overline{G}_{r-1}^{(\mu)}
\end{split}
\end{align*}
Here we defined $\ve_k\seteq(0,\ldots,0,\overset{k}{\check{1}},0,\ldots,0)$. The induction hypothesis on $r$ yields
\begin{align}\label{eq:yoko:tmp}
\begin{split}
\overline{G}_r^{(\mu)}
&=\biggl\{ 
   \frac{t^{-(r-1)n}}{q^{(r-1)(|\mu|+1)}} 
   \sum_{l=0}^{r-1}
   \bigg[
    (-1)^l q^{-\binom{l}{2}}q^{-l(r-l-1)}(q^{r-l};q)_l 
\\
&\hskip 10em
    \cdot\sum_{k=1}^n \bigl(1-\dfrac{1}{q^{\mu_k}t} \bigr)
    \biggl(\prod_{i \neq k} \dfrac{1-y_i/q^{\mu_i}t y_k}{1-y_i/y_k} \biggr)
     T_{q,y_k} H_l^{(\mu+\ve_k)} 
   \bigg]
\\
 &\hskip 2em
  +\dfrac{t^{-n}}{q^{|\mu|}} \dfrac{t^{-(r-1)n}}{q^{(r-1)|\mu|}} 
           \sum_{l=0}^{r-1} 
            (-1)^l q^{-\binom{l}{2}} q^{-l(r-l-1)}(q^{r-l};q)_l H_n^{l,(\mu) }
  \biggl\} \phi_y \cdot 1.
\end{split}
\end{align}
We now define 
\begin{align*}
h_{\nu}^{(\mu)} 
\seteq
\biggl(\prod_{i<j}\frac{q^{\nu_i}y_i - q^{\nu_j}y_j}{y_i-y_j}\biggr) 
\prod_{i,j=1}^n \frac{(t q^{\mu_j}y_i/y_j ; q)_{\nu_i}}{(q y_i/y_j;q)_{\nu_i}} 
\end{align*}
for $\nu = (\nu_1,\ldots,\nu_n) \in \mathbb{N}^n$. Then one obtains
\begin{align*}
H_n^{l,(\mu+\ve_k)} 
 = \sum_{|\nu|=l} h_{\nu}^{(\mu+\ve_k)} \prod_{i=1}^n T_{q,y_i}^{\nu_i}
 = \sum_{|\nu|=l} 
    \dfrac{1-q^{\mu_k + \nu_k}t}{1-q^{\mu_k}t} 
    \biggl(
     \prod_{i \neq k}\dfrac{q^{\mu_k + \nu_i} t y_i-y_k}{q^{\mu_k} t y_i-y_k} 
    \biggr)
    h_{\nu}^{(\mu)} \prod_{i=1}^n T_{q,y_i}^{\nu_i}.
\end{align*}
One can also find that
\begin{align*}
&h_{\nu + \ve_k}^{(\mu)} = 
\dfrac{1-q^{\mu_k + \nu_k}t}{1-q^{\nu_k + 1}}
\biggl(\prod_{i \neq k} 
 \dfrac{q^{\nu_k +1}y_k-q^{\nu_i}y_i}{q^{\nu_k}y_k-q^{\nu_i}y_i} 
 \dfrac{q^{\mu_i + \nu_k}t y_k-y_i}{q^{\nu_k +1}y_k-y_i}\biggr)
 h_{\nu}^{(\mu)},
\\
&T_{q,y_k}h_{\nu}^{(\mu)} = h_{\nu + \ve_k}^{(\mu)} \frac{1-q^{\nu_k +1}} {1-q^{\mu_k + \nu_k }t}
\biggl(\prod_{i \neq k} \frac{q^{\mu_k -1}t y_i-y_k}{y_i-q^{\mu_i}t y_k} 
\frac{y_k-q^{\nu_i}y_i}{y_k-q^{\mu_k + \nu_i -1}y_i}\biggr)
T_{q,y_k}.
\end{align*}
Thus the first part of (\ref{eq:yoko:tmp}) becomes
\begin{align}\label{eq:yoko_zenka_bubun}
\begin{split}
&\sum_{k=1}^n 
\frac{q^{\mu_k}t-1}{q^{\mu_k}t} 
\biggl(\prod_{i \neq k}\frac{1-y_i/q^{\mu_i}t y_k}{1-y_i/y_k} \biggr)
T_{q,y_k} H_n^{l,(\mu+\ve_k)}\\
=&\frac{t^{-n}}{q^{|\mu|}}
\sum_{|\nu|=l} \sum_{k=1}^n (q^{\nu_k +1}-1) 
\biggl(\prod_{i \neq k} \frac{y_k-q^{\nu_i}y_i}{y_k-y_i}\biggr)
h_{\nu + \ve_k}^{(\mu)} 
\biggl(\prod_{i \neq k} T_{q,y_i}^{\nu_i}\biggr)  T_{q,y_k}^{\nu_k +1} 
=
\dfrac{t^{-n}}{q^{|\mu|}}(q^{l+1}-1)H_n^{l+1,(\mu)}.
\end{split}
\end{align}
At the last equality we used the relation
\begin{align*}
\sum_{k=1}^n (q^{\nu_k}-1)
\biggl(\prod_{i \neq k} \frac{y_k-q^{\nu_i}y_i}{y_k-y_i}\biggr)  
= q^{|\nu|}-1,
\end{align*}
which is obtained by specializing $z=0$ in the next partial fraction expansion.
\begin{align*}
\prod_{i=1}^n \frac{z-q^{\nu_i}y_i}{z-y_i} 
= \sum_{k=1}^n \frac{(1-q^{\nu_k})y_k}{z-y_k} 
 \biggl(\prod_{i \neq k}\frac{y_k-q^{\nu_i}y_i}{y_k-y_i}\biggr) +1
\end{align*}

Using (\ref{eq:yoko_zenka_bubun}), we can compute (\ref{eq:yoko:tmp}) as 
\begin{align*}
\overline{G}_r^{(\mu)} =&
\biggl\{ \frac{t^{-(r-1)n}}{q^{(r-1)(|\mu|+1)}} \sum_{l=0}^{r-1} (-1)^l q^{-\binom{l}{2}}q^{-l(r-l-1)}(q^{r-l};q)_l  
\frac{t^{-n}}{q^{|\mu|}}(q^{l+1}-1)H_n^{l+1,(\mu)} \\
&+ \frac{t^{-r n}}{q^{r|\mu|}} 
\sum_{l=0}^{r-1} (-1)^l q^{-\binom{l}{2}}q^{-l(r-l-1)}(q^{r-l};q)_l H_n^{l,(\mu)} \biggl\} \phi_y \cdot 1.
\end{align*}
A direct calculation shows that the coefficient of $H_n^{l+1,(\mu)}\phi_y\cdot 1$ sums up to the desired one.
\end{proof}

\subsubsection{Remark on the eigenvalues}
We make a comment on a property of the eigenvalues $e_r(s^\lambda)$ and $g_r(s^\lambda;q,t)$ of the operators $\widehat{E}_r$ and $\widehat{G}_r$. 

\begin{lem}
For a partition $\lambda \vdash n$, define the {\it spectral parameter} $s^\lambda(q,t)$ by
\begin{align*}
s^\lambda(q,t)  \seteq q^\lambda t^{-\delta}
=(q^{\lambda_1} t^{-1},q^{\lambda_2} t^{-2},\ldots).
\end{align*}
Then
\begin{align*}
e_r\bigl(s^\lambda(t^{-1},q^{-1})\bigr)
 =(-1)^r g_r\bigl(s^{\lambda'}(q,t);q,t\bigr).
\end{align*}
\end{lem}

\begin{proof}
Using the generating functions 
(\ref{eq:G(y)}) and (\ref{eq:E(y)}), we can compute easily as
\begin{align*}
\sum_{r\ge 0}g_r\bigl(s^{\lambda'}\bigl(q,t);q,t) u^r 
& = \prod_{i\ge 1}
     \dfrac{\bigl(q^{\lambda'_i} t^{-i+1} u;q\bigr)_\infty}
           {\bigl(q^{\lambda_i'} t^{-i}   u;q\bigr)_\infty}    
= \bigl(q^{\lambda'_1} u;q\bigr)_\infty 
    \prod_{i\ge 1}
     \dfrac{\bigl(q^{\lambda'_{i+1}} t^{-i} u;q\bigr)_\infty}
           {\bigl(q^{\lambda_i'}     t^{-i} u;q\bigr)_\infty}  
    \\
& = \bigl(q^{\lambda'_1} u;q\bigr)_\infty 
    \prod_{i\ge 1} 
     \bigl(q^{\lambda'_{i+1}} t^{-i} u;q\bigr)_{\lambda_i'-\lambda_{i+1}'}.
\end{align*}
In the last line, the $i$-th term is $1$ unless $\lambda'_i>\lambda'_{i+1}$. If we denote the parts of $\lambda$ and $\lambda'$ as
\begin{align*}
\lambda
 =\bigl(\lambda_1,\ldots,\lambda_l\bigr) 
 =\bigl(m_1^{n_1},\ldots,m_k^{n_k}\bigr),      \quad
\lambda'
 =\bigl(\lambda'_1,\ldots,\lambda'_{l'}\bigr),
\end{align*}
then the condition $\lambda'_i>\lambda'_{i+1}$ implies that $i=n_j$ for some $j$ and that $\lambda'_{n_j}-\lambda'_{n_j+1}=m_j$. Thus
\begin{align*}
\sum_{r\ge 0}g_r\bigl(s^{\lambda'}(q,t);q,t\bigr) u^r 
& = \bigl(q^{\lambda'_1} u;q\bigr)_\infty 
    \bigl(t^{-n_1} u;q\bigr)_{m_j}
    \prod_{j=2}^k 
     \bigl(q^{m_{\lambda_1}+\cdots+m_{\lambda_{j-1}}} t^{-n_j} u;q\bigr)_{m_j}.
\end{align*}
Careful consideration yields
\begin{align*}
\sum_{r\ge 0}g_r\bigl(s^{\lambda'}(q,t);q,t\bigr) u^r 
& = \bigl(q^{\lambda'_1} u;q\bigr)_\infty 
    \prod_{i=1}^l (1-q^i t^{-\lambda_i} u)   
  = \prod_{i\ge 1} \biggl(1+(t^{-1})^{\lambda_i} (q^{-1})^{-i} (-u)\biggr) \\
& = \sum_{r\ge 0}e_r\bigl(s^\lambda(t^{-1},q^{-1})\bigr)(-u)^r.
\end{align*}

\end{proof}

\section{Elliptic algebras $\calA(p)$ and $\calM(p)$}
\label{sect:elliptic}

The goal of this section is to introduce an elliptic counterpart $\calM(p)$ of the commutative algebra of operators $\calM$, by using the algebra $\calA(p)$ and the elliptic analogue $\calU(q,t,p q^{-1}t )$ of the Ding-Iohara algebra. 

In \S\S\,\ref{subsec:ell_A} we summarize the properties of the elliptic counterpart $\calA(p)=\calA(q_1,q_2,q_3,p)$ of the algebra $\calA(q_1,q_2,q_3)$ over $\bbC\bbP^1$, recalling the original definition given in \cite{FO:1997}.

As for the free field side, we introduce a quasi-bialgebra $\calU(q,t,p q^{-1}t )$ in the sense of Drinfeld \cite{QHA} and Babelon-Bernard-Billey \cite{BBB}. In its construction, summarized in the appendix, we twist the Drinfeld coproduct $\Delta$ (see Proposition \ref{prop:Hopf-alg}) into a $p$-depending one $\Delta_{p q^{-1}t}$ (see Proposition \ref{prop:QBA}) by applying the twistor $F(p q^{-1}t)$ (see Definition \ref{dfn:F}).
Then we define the dressed Drinfeld generators $\eta(z,p q^{-1}t), \xi(z,p q^{-1}t), 
\varphi^\pm(z,p q^{-1}t)$ 
as in Definition \ref{dfn:dressed_currents}, which satisfy the relations in Proposition \ref{prop:dressed_copro}. Note that here we have used the deformation parameter $p q^{-1}t$ 
instead of $p$, to have the same elliptic nome $p$ for $\calA(p)$
appearing in our formulas.

\subsection{Commuting elliptic algebra $\calA(p)$}
\label{subsec:ell_A}
{}From now on, we regard  $q_1=q^{-1}$, $q_2=t$ and $p$ as complex parameters satisfying the conditions
$|q|<1, |t^{-1}|<1,|p|<1,|p q^{-1}t|<1$ and $q^it^j p^k\neq1\ \mbox{for any\ } 
(i,j,k)\in\bbZ^3\setminus\{(0,0,0)\}$. Set $q_3\seteq 1/q_1q_2$.
We use the notation for the theta function
\begin{align*}
\Theta_p(x)
\seteq(p;p)_\infty(x;p)_\infty(p/x;p)_\infty,
\end{align*}
written in the {\it multiplicative} notation.
It enjoys the quasi periodicity
$\Theta_p(e^{2\pi \sqrt{-1}}x)=\Theta_p(x)$ and
$\Theta_p(p x)=-x^{-1}\Theta_p(x)$.

\begin{dfn}[Space $\calA(p)$]
\label{dfn:ell_calA}
For $n\in\bbN$, the vector space $\calA_{n}(p)=\calA_{n}(q_1,q_2,q_3,p)$ is defined by the following conditions (i), (ii) and (iii).
\\
(i)
$\calA_0(p) \seteq \bbC$. For $n\geq 1$, $f(x_1,\ldots,x_n)\in \calA_n(p)$ is a 
symmetric meromorphic function over $\bbC^\times$ 
satisfying the double periodicity 
$f(x_1,\ldots,e^{2\pi \sqrt{-1}} x_i,\ldots,x_n)=f(x_1,\ldots,p x_i,\ldots,x_n)=
f(x_1,\ldots,x_n)$ for any $i$. \\
(ii)
The poles of $f\in \calA_n(p)$ are located only on the diagonal $\{(x_1,\ldots,x_n) \mid \exists (i,j),  i\neq j ,x_i=x_j\}$ modulo $p$-shifts, and the orders of the poles are at most two.
\\
(iii) For $n\geq 3$, $f\in \calA_n(p)$ satisfies the wheel conditions
\begin{align*}
f( x_1,q_1 x_1,q_1 q_2 x_1,x_4,\ldots)=0,\qquad
f( x_1,q_2 x_1,q_1 q_2 x_1,x_4,\ldots)=0.
\end{align*}

Then set the graded vector space 
$\calA(p)=\calA(q_1,q_2,q_3,p)\seteq\bigoplus_{n\geq 0}\calA_n(p)$.
\end{dfn}

Introduce the star product $*$ on $\calA(p)$ by the equation (\ref{eq:star}), 
while replacing the structure function $\omega$ by
\begin{align}
\label{eq:omega_ell}
\omega(x,y;q_1,q_2,q_3,p) \seteq 
\dfrac{\Theta_p(q_1 y/x)\Theta_p(q_2 y/x)\Theta_p(q_3 y/x)}{\Theta_p(y/x)^3}.
\end{align}

\begin{prop}\label{prop:calA(p)}
$(\calA(p),*)$ is a commutative algebra having a basis $(\ep_\lambda(x;q_i,p))$, where
\begin{align*}
&
\ep_\lambda \seteq \ep_{\lambda_1} * \cdots * \ep_{\lambda_l}
\quad (\lambda=(\lambda_1,\ldots,\lambda_l)),
\\
&
\ep_n(x_1,\ldots,x_n;q_i,p) \seteq 
 \prod_{1\le j<k\le n}
 \dfrac{\Theta_p(q_i x_j/ x_k)\Theta_p(x_j/q_i x_k)}{\Theta_p(x_j/x_k)^2}.
\end{align*}
\end{prop}

\begin{proof}
One can immediately find that the space $\calA(p)$ is closed by the star product $*$. 
Hence we can directly proceed to the argument of the Gordon filtration
for the study of the dimension of $\calA_n(p)$. 
Since the Gordon filtration is defined just by using the zero conditions, 
analogues of Definition \ref{dfn:gordon_filtr},
Lemma \ref{lem:atmost1}, Propositions \ref{prop:atleast1} and
\ref{prop:wbcalA_n:poincare} for $\calA(p)$ hold. 
Hence  $(\ep_\lambda(x;q_i,p))$ forms a basis of $\calA(p)$.
The commutativity of $\calA(p)$ can be proved in the same way 
as in Proposition \ref{prop:commute} for $\calA$.  
\end{proof}

Using the same technique for the proof of 
Proposition \ref{prop:Wronski:op}, 
one can show the following elliptic analogue of 
the Wronski relation (\ref{eq:sn}).

\begin{prop}[Elliptic Wronski relation]\label{prop:ell_wronski}
For $n\in\bbN_+$, we have
\begin{align*}
\sum_{k=0}^n
 \binom{n}{k}\Theta_p(q_2^{n-k}/q_1^k)
 \left(-\dfrac{\Theta_p(q_2)}{\Theta_p(1/q_1)}\right)^k
 \ep_{n-k}(x_1,\ldots,x_{n-k};q_1,p) * \ep_{k}(x_{n-k+1},\ldots,x_n;q_2,p)
=0.
\end{align*}
\end{prop}

\subsection{Elliptic current $\eta(z;p q^{-1}t)$}
Let $\eta(z;p q^{-1}t)$ be the vertex operator given in (\ref{eq:eta_p}),
with $p$ replaced by $p q^{-1}t$.
Explicitly, it is written as
\begin{align}\label{eq:vertex_operator_p}
\eta(z;p q^{-1}t)=
\exp\Bigg(\sum_{n>0} 
\dfrac{1-t^{-n}}{n}\dfrac{1-p^n q^{-n}t^n}{1-p^n}a_{-n} z^{n}\Bigg)
\exp\Bigg(-\sum_{n>0} \dfrac{1-t^{n}}{n}  a_n z^{-n}\Bigg).
\end{align}
\begin{lem}
We have
\begin{align*}
\begin{split}
&\eta(z;p q^{-1}t)\eta(w;p q^{-1}t)\\
&=(1-w/z)
\dfrac{(q t^{-1}w/z;p)_\infty (p t w/z;p)_\infty (p q^{-1}w/z;p)_\infty}
      {(q w/z;p)_\infty (t^{-1}w/z;p)_\infty(p q^{-1} t w/z;p)_\infty}
:\eta(z;p q^{-1}t)\eta(w;p q^{-1}t):.
\end{split}
\end{align*}
\end{lem}

{}From this lemma, we have the operator product 
\begin{align*}
&\dfrac{1}{\omega(z_1,z_2;q^{-1},t,q t^{-1},p)}
 \eta(z_1;p q^{-1}t)\eta(z_2;p q^{-1}t)\\
&=
\dfrac
 {q t^{-1} (z_2/z_1;p)_\infty^2 (p z_2/z_1;p)_\infty}
 {(q z_2/z_1;p)_\infty (t^{-1}z_2/z_1;p)_\infty (p q^{-1}t z_2/z_1;p)_\infty}
\\
&\hskip 6em
\times
\dfrac
 {(z_1/z_2;p)_\infty^2(p z_1/z_2;p)_\infty} 
 {(q z_1/z_2;p)_\infty(t^{-1}z_1/z_2;p)_\infty(p q^{-1}t z_1/z_2;p)_\infty}
:\eta(z_1;p q^{-1}t)\eta(z_2;p q^{-1}t):,
\end{align*}
which we will regard as a symmetric Laurent series 
defined on the annulus $M<|z_2/z_1|<M^{-1}$, where $M\seteq\max(|q|,|t^{-1}|,|p q^{-1} t|)<1$.

Introduce 
a mapping $\calO_p(\cdot)$ defined on $\calA_n(p)$ by
\begin{align*}
\calO_p(f)\seteq\biggl[
 \dfrac{f(z_1,\ldots,z_n)}
       {\prod_{1\le i<j\le n}\omega(z_i,z_j;q^{-1},t,q t^{-1};p)}
 \eta(z_1;p q^{-1}t)\cdots\eta(z_n;p q^{-1}t)\biggr]_1,
\end{align*}
and extend it by linearity.
Then the compatibility 
\begin{align*}
\calO_p(f*g)=\calO_p(f)\calO_p(g)
\end{align*}
holds for $f,g\in\calA(p)$.

Set $\calM(p)\seteq\{\calO_p(f) \mid f\in \calA(p)\}$.
The commutativity of the algebra $\calA(p)$ is inherited to $\calM(p)$. 

\begin{thm}\label{thm:calM(p)}
$\calM(p)$ is a family of commutative operators acting on the 
Fock space $\calF$.
\end{thm}

\begin{rmk}
Note that operators in $\calM(p)$ are acting on the {\it same} Fock space as in the case of the Macdonald operators $\calM$. One can check that $f\in\calM(p)$ and $g\in \calM$ do not commute in general.
By explicit calculation, one finds that the eigenvalues of the
operators in $\calM(p)$ in general can {\it not} be written as rational functions in parameters $q$, $t$ and $p$ in contrast to the Macdonald case.
See Conjecture \ref{con:eigen} below.
At present,
we do not know whether the 
mapping $\calO_p(\cdot)$ is injective or not. (As for the Macdonald case, see 
Proposition \ref{prop:spectral_criterion}.)
\end{rmk}

\subsection{Intertwining operator $\widetilde{\Phi}(y;p q^{-1}t)$
and the Ruijsenaars difference operator}
\label{subsec:ruij}

In this subsection, we consider an elliptic analogue of 
Proposition \ref{prop:recover_macd} in \S\S \ref{subsec:DI}.

Set 
\begin{align*}
\widetilde{\Phi}(y;p q^{-1}t)\seteq
\exp\Bigg(\sum_{n>0}
 \dfrac{1}{n}\dfrac{1-t^n}{1-q^{n}}\dfrac{1-p^n q^{-n}t^n}{1-p^n}
 t^{-n}a_{-n} y^{n}
\Bigg)
\exp\Bigg(-\sum_{n>0}\dfrac{1}{n}\dfrac{1-t^n}{1-q^{n}}q^{n}a_{n} y^{-n}
\Bigg).
\end{align*}
{}From (\ref{eq:Phi_boson2}), we have
$\Phi(y;p)=\widetilde{\Phi}(q^{1/2}y;p)$, where 
$\Phi(y;p)$ denotes the generating function of the intertwining operator
$\Phi^{\calF}_{V_x\otimes \calF}(p):V_x\otimes \calF\rightarrow \calF$
with respect to the elliptic algebra 
$\calU(q,t,p)$
described in Proposition \ref{prop:ITP_p}.

The elliptic analogue of Proposition \ref{prop:recover_macd}
is given as follows.
\begin{prop}
We have
\begin{align*}
&
\eta(z;p q^{-1}t) \widetilde{\Phi}(y_1;p q^{-1}t)
 \cdots \widetilde{\Phi}(y_n ;p q^{-1}t)\\
&
=t^{-n}\prod_{i=1}^n 
 \dfrac{\Theta_p(q t^{-1} z/y_i)}{\Theta_p(q z/y_i)}
 \dfrac{\Theta_p(t z/y_i)}{\Theta_p(z/y_i)}
 \widetilde{\Phi}(y_1 ;p q^{-1}t) \cdots 
 \widetilde{\Phi}(y_n ;p q^{-1}t)\eta(z;p q^{-1}t)\\
&+(1-t^{-1})
  \dfrac{(p/t;p)_\infty(p t/q;p)_\infty}{(p;p)_\infty(p/q;p)_\infty}\\
&\hskip 2em
 \times
 \sum_{i=1}^n t^{-i+1}
 \prod_{j=1}^{i-1}
 \dfrac{\Theta_p(q t^{-1} z/y_j)}{\Theta_p(q z/y_j)}
 \dfrac{\Theta_p(t z/y_j)}{\Theta_p(z/y_j)} \delta(y_i/z) T_{q,y_i} 
 \widetilde{\Phi}(y_1;p q^{-1}t)\cdots \widetilde{\Phi}(y_n;p q^{-1}t),
\\
&
 \widetilde{\Phi}(y_1;p q^{-1}t)\cdots \widetilde{\Phi}(y_n;p q^{-1}t)
\\
&=\prod_{1\leq i<j\leq n}
  \dfrac{\Gamma(q t^{-1}y_j/y_i;q,p)}{\Gamma(q y_j/y_i;q,p)}
  \prod_{1\leq k\neq l\leq n}
  \dfrac{(p t y_l/y_k,q,p)_\infty}{(p y_l/y_k,q,p)_\infty}
  :\widetilde{\Phi}(y_1;p q^{-1}t) \cdots \widetilde{\Phi}(y_n;p q^{-1}t):.
\end{align*}
Here we have used the notation for the double infinite product
$(x;q,p)_\infty \seteq \prod_{i,j\in \bbN}(1-q^i p^j x)$, and
the elliptic gamma function 
$\Gamma(x;q,p)\seteq(p q/x;q,p)_\infty/(x;q,p)_\infty$.
\end{prop}

Noting that $\Gamma(q z;q,p)=\Theta_p(z)\Gamma(z;q,p)$, we have
the following.

\begin{cor}\label{cor:Ruij}
Set 
\begin{align*}
\phi(y_1,\ldots ,y_n;p)\seteq
 \prod_{1\leq k\neq l\leq n}
  \dfrac{(p t y_l/y_k,q,p)_\infty}{(p y_l/y_k,q,p)_\infty}
:\widetilde{\Phi}(y_1;p q^{-1}t)\cdots \widetilde{\Phi}(y_n;p q^{-1}t):,
\end{align*}
for simplicity of display.
Then we have
\begin{align}\label{eq:Ruij}
\begin{split}
\bigg[\eta(z;p q^{-1}t)\bigg]_1 &\phi(y_1,\ldots ,y_n;p)
= 
\phi(y_1,\ldots ,y_n;p)
 \bigg[
  t^{-n}\prod_{i=1}^n 
  \dfrac{\Theta_p(q t^{-1} z/y_i)}{\Theta_p(q z/y_i)}
  \dfrac{\Theta_p(t z/y_i)}{\Theta_p(z/y_i)}
  \eta(z;p q^{-1}t)\bigg]_1\\
  &+(1-t^{-1})t^{-n+1}
    \dfrac{(p/t;p)_\infty(p t/q;p)_\infty}{(p;p)_\infty(p/q;p)_\infty}
    D_{n,y}^1(p)  \phi(y_1,\ldots ,y_n;p),
 \end{split}
\end{align}
where $[\cdot ]_1$ means the constant term in $z$,  $D_{n,y}^1(p)$ denotes
the Ruijsenaars difference operator \cite{Ruijsenaars:1987}
acting on the variable $y_i$'s:
\begin{align*}
D_n^1(p) \seteq \sum_{i=1}^n \prod_{j\neq i}
\dfrac{\Theta_p(t y_i/y_j)}{\Theta_p(y_i/y_j)}T_{q,y_i}.
\end{align*}
\end{cor}

\begin{rmk}
Contrary to the case of $p=0$, application of the operator in the first term of  RHS of (\ref{eq:Ruij}) to $1$, namely ${\displaystyle
\bigg[\prod_{i=1}^n 
\dfrac{\Theta_p(q t^{-1} z/y_i)}{\Theta_p(q z/y_i)}\dfrac{\Theta_p(t z/y_i)}{\Theta_p(z/y_i)}\eta(z;p q^{-1}t)\bigg]_1\cdot 1}$, 
remains quite nontrivial. This 
prevent us from interpreting $\phi(y_1,\ldots,y_n;p)\cdot 1$
as a reproduction kernel function for 
the Ruijsenaars difference operator $D_{n,y}^1(p)$.
\end{rmk}

\begin{rmk}
Langmann studied the elliptic Calogero-Sutherland model 
using the technique of the quantum field theory. See \cite{L1,L2} and 
references therein.
He used 
finite temperature correlation functions to 
introduce the elliptic parameter in his calculation. 
Corollary \ref{cor:Ruij} may be regarded as its difference version. 
However, precise identification between the finite temperature 
calculation and the quasi-Hopf treatment remains unclear.

\end{rmk}
\subsection{Conjecture about the eigenvalues of the operator $\big[\eta(z;p q^{-1}t)\big]_1$}
Here we briefly discuss the eigenvalues of 
the Ruijsenaars difference operator $D_{n,y}^1(p)$ and 
those of $\big[\eta(z;p q^{-1}t)\big]_1$. 
One can regard $D_{n,y}^1(p)$ as an operator acting on the
space of formal Laurent power series 
$y^{\lambda}\cdot \bbF[[y_2/y_1,y_3/y_2,\ldots,y_n/y_{n-1},p y_1/y_n,p]]$,
where $y^{\lambda}\seteq y_1^{\lambda_1}\cdots y_n^{\lambda_n}$ $(\lambda\in \bbC^n)$.
Namely, set
\begin{align*}
&f_\lambda(y,p)
 \seteq
 y^{\lambda}\sum_{i_1,\ldots,i_n\in \bbN}
 c_{i_1,\ldots,i_n}(\lambda) (y_2/y_1)^{i_1}\cdots (y_n/y_{n-1})^{i_{n-1}}
 (p y_1/y_{n})^{i_{n}},\\
&\varepsilon_\lambda(p,n)
 \seteq
 \sum_{k\in \bbN} p^k \varepsilon_{\lambda,k}(n) 
\end{align*}
where $c_{i_1,\ldots,i_n}(\lambda),
\varepsilon_{\lambda,k}(n) \in \bbF(q^{\lambda_1},\cdots ,q^{\lambda_n})$. 
Normalize $f_\lambda(y,p)$ as $f_\lambda(y,p)=y^{\lambda}+\cdots$ by letting
$c_{i,\ldots,i}=\delta_{i,0}$, and 
set $\varepsilon_{\lambda,0}(n)=
t^{n-1}q^{\lambda_1}+t^{n-2}q^{\lambda_2}+\cdots +q^{\lambda_n}$.
Then we can uniquely determine the coefficients 
$c_{i_1,\ldots,i_n}(\lambda)$ and $\varepsilon_{\lambda,k}(n)$
iteratively by perturbation in $p$ and 
by imposing the eigenvalue condition
 $D_{n,y}^1(p) f_\lambda(y,p)=\varepsilon_\lambda(p,n) f_\lambda(y,p)$.

Consider the case when $\lambda$ is a partition. 
Since the lowest order term in $p$ in this eigenvalue equation 
should be satisfied by the Macdonald polynomial, we must have
$f_\lambda(y,p)=P_\lambda(y;q,t)+O(p)$ and
$\varepsilon_\lambda(p,n)=\sum_{i=1}^n t^{n-i}q^{\lambda_i}+ O(p)$. 
The higher order corrections in $p$ can be calculated in principle. 
As an example, for the simplest case $\lambda=\emptyset$, we have
\begin{align*}
&f_\emptyset(y,p)
 =1+p q^{-1} t \dfrac{(1-t^{-1})(1-t^{-n})}{(1-q^{-1})(1-q^{-1}t^{-n+1})}
  \sum_{1\leq i\neq j\leq n}y_j/y_i+O(p^2),\\
&t^{-n+1}\varepsilon_\emptyset(p,n)
 =\dfrac{1-t^{-n}}{1-t^{-1}}
  +p q^{-1} t \dfrac{(1-q t^{-1})(1-t^{-n})(1-t^{-n+1})}{1-q^{-1}t^{-n+1}}
  +O(p^2).
\end{align*}
In general,
one may notice that the correction terms depend on $n$ in a nontrivial manner. 

A brute force calculation suggests the following. 

\begin{con}\label{con:eigen}
The eigenvalues of the operator $\big[\eta(z;p q^{-1}t)\big]_1$
on $\calF$ is given by the eigenvalue of
the Ruijsenaars difference operator $D_{n,y}^1(p)$ as
${\displaystyle
\lim_{n\rightarrow \infty}(1-t^{-1})
 \dfrac{(p/t;p)_\infty(p t/q;p)_\infty}{(p;p)_\infty(p/q;p)_\infty}t^{-n+1}
\varepsilon_\lambda(p,n)}$
with $\lambda$ being partitions.
\end{con}
 
\begin{rmk}
In view of (\ref{eq:Ruij}),
this conjecture is equivalent to saying that
the {\it unwanted } operator in (\ref{eq:Ruij}) satisfies
${\displaystyle 
\lim_{n\rightarrow \infty} t^{-n}
\bigg[\prod_{i=1}^n 
\dfrac{\Theta_p(q t^{-1} z/y_i)}{\Theta_p(q z/y_i)}
\dfrac{\Theta_p(t z/y_i)}{\Theta_p(z/y_i)}
\eta(z;p q^{-1}t)\bigg]_1\cdot 1=0}$.
\end{rmk}

Take the case $\lambda=\emptyset$ as an example. 
It is clear that $\big[\eta(z;p q^{-1}t)\big]_1\cdot 1=1$, showing that 
$1$ is an eigenvalue. This corresponds to the limit
\begin{align*}
&\lim_{n\rightarrow \infty}(1-t^{-1})
 \dfrac{(p/t;p)_\infty(p t/q;p)_\infty}{(p;p)_\infty(p/q;p)_\infty}
 \left(
{1-t^{-n}\over 1-t^{-1}}+p q^{-1} t 
{(1-qt^{-1})(1-t^{-n})(1-t^{-n+1})\over 1-q^{-1}t^{-n+1}}+O(p^2)\right)\\
&=
(1-t^{-1})
 \dfrac{(p/t;p)_\infty(p t/q;p)_\infty}{(p;p)_\infty(p/q;p)_\infty}
 \left(
{1\over 1-t^{-1}}+p q^{-1} t (1-qt^{-1})+O(p^2)\right)\\
&=1+O(p^2).
\end{align*}

\subsection{Elliptic bosons and the relation to the Okounkov-Pandharipande operator}\label{subsec:OP}
In this subsection, we claim that the
operator of Okounkov-Pandharipande \cite{OkounkovPandharipande:2004}
is derived from our elliptic deformation of the Macdonald operator 
$\big[\eta(z,p q^{-1}t)\big]_1$
when we set $q=e^\hbar, t=q^{\beta\hbar}$ and consider the limit $\hbar \rightarrow 0$ while $p$ being fixed.

To make our computation simple, we use the Heisenberg algebra 
with the commutation relation
\begin{align}\label{eq:ruij_boson}
[\lambda_m,\lambda_n]
=-\dfrac{1}{m}
\dfrac{(1-q^m)(1-t^{-m})(1-p^m q^{-m} t^{m})}{1-p^{m}}
\delta_{m+n,0}.
\end{align}
Relating $a_n$ with $\lambda_n$ in a suitable manner, 
we may write  the elliptic current $\eta(z,p q^{-1}t)$ as
\begin{align*}
\eta(z;p q^{-1}t)=:\exp\Biggl(\sum_{n\neq0}\lambda_n z^{-n}\Biggr):.
\end{align*}

Recall the conventional normalization for the 
boson $\overline{a}_n$ for the Virasoro algebra 
\begin{align*}
[\overline{a}_m,\overline{a}_n]=m\delta_{m+n,0}.
\end{align*}
Setting
\begin{align*}
\lambda_n\seteq
 \dfrac{1}{|n|}
 \sqrt{-\dfrac{(1-q^{|n|})(1-t^{-|n|})(1-p^{|n|}q^{-|n|}t^{|n|})}{1-p^{|n|}}}
 \cdot \overline{a}_n,
\end{align*}
we have
\begin{align*}
\lambda_n&=\bigg[
  \beta^{1/2}\hbar
  +\dfrac{n}{4}\dfrac{1+p^n}{1-p^n}(1-\beta)\beta^{1/2}\hbar^2\\
&\hskip 2em
  +\dfrac{n^2}{96}
   \left(4(2-3\beta+2\beta^2)\beta^{1/2}
     -3\dfrac{(1+p^n)^2}{(1-p^n)^2}(1-\beta)^2\beta^{1/2}
   \right)\hbar^3+O(\hbar^4)
 \bigg] \cdot \overline{a}_n.
\end{align*}
Then the operator $\big[\eta(z;p q^{-1}t)\big]_1$, which has appeared several times in our arguments, has the next expansion in $\hbar$.

\begin{prop}\label{prop:O-P}
We have
\begin{align*}
\big[\eta(z;p q^{-1}t)\big]_1
=&1+\beta\sum_{n\ge 1}\overline{a}_{-n}\overline{a}_n\hbar^2
  +\bigg[
    \beta(1-\beta)\sum_{n\ge 1}\dfrac{n}{2}\dfrac{1+p^n}{1-p^n}
    \overline{a}_{-n}\overline{a}_n
\\
&\hskip 3em
    +\dfrac{\beta^{3/2}}{2}
      \sum_{n,m\ge 1}\bigl(
       \overline{a}_{-n}\overline{a}_n\overline{a}_{n+m}   
       +\overline{a}_{-n-m}\overline{a}_n\overline{a}_m\bigr)
   \bigg]\hbar^3+O(\hbar^4).
\end{align*}
\end{prop}

With a suitable change of the notations for the parameters,
we notice that 
the third term coincides 
with the operator ${\sf M}(q,t_1,t_2)$ 
of Okounkov and Pandharipande in \cite{OkounkovPandharipande:2004}.
Hence this suggests that  the 
quasi-Hopf twisting of the Ding-Iohara algebra $\calU(q,t,p q^{-1} t)$
might be relevant to the quantum cohomology.

\appendix
\section{Ding-Iohara's quantum algebra}\label{app:DI}
This appendix deals with the Ding-Iohara quantum algebra \cite{DI:1997}. 
Restricting ourselves to a particular case, which we denote by $\calU(q,t)$, 
we study its representation theories and their connection with the results in the main text.

\subsection{Definition of $\calU(q,t)$}

Set
\begin{align*}
g(z)\seteq\dfrac{G^+(z)}{G^-(z)},\qquad
G^\pm(z)\seteq(1-q^{\pm1}z)(1-t^{\mp 1}z)(1-q^{\mp1}t^{\pm 1}z).
\end{align*}
Note that $g(z)$ satisfies the Ding-Iohara requirement $g(z)=g(z^{-1})^{-1}$.

\begin{dfn}[Ding-Iohara]\label{dfn:calU(q,t)}
Let $\calU=\calU(q,t)$ be a unital associative algebra generated by
the Drinfeld currents 
\begin{align*}
x^\pm(z)=\sum_{n\in \bbZ}x^\pm_n z^{-n},\qquad 
\psi^\pm(z)=\sum_{\pm n\in \bbN}\psi^\pm_n z^{-n},
\end{align*}
and the central element $\gamma^{\pm 1/2}$, satisfying the defining relations
\begin{align*}
&\psi^\pm(z) \psi^\pm(w)= \psi^\pm(w) \psi^\pm(z),\\ 
&\psi^+(z)\psi^-(w)=
\dfrac{g(\gamma^{+1} w/z)}{g(\gamma^{-1}w/z)}\psi^-(w)\psi^+(z),\\
&\psi^+(z)x^\pm(w)=g(\gamma^{\mp 1/2}w/z)^{\mp1} x^\pm(w)\psi^+(z),\\
&\psi^-(z)x^\pm(w)=g(\gamma^{\mp 1/2}z/w)^{\pm1} x^\pm(w)\psi^-(z),\\
&
[x^+(z),x^-(w)]=\dfrac{(1-q)(1-1/t)}{1-q/t}
\bigg( \delta(\gamma^{-1}z/w)\psi^+(\gamma^{1/2}w)-
\delta(\gamma z/w)\psi^-(\gamma^{-1/2}w) \bigg),\\
&x^\pm(z)x^\pm(w)=g(z/w)^{\pm1}x^\pm(w)x^\pm(z).
\end{align*}
\end{dfn}

\begin{prop}[Ding-Iohara]\label{prop:Hopf-alg}
The algebra $\calU$ has a Hopf algebra structure defined by the following
morphisms.\\
Coproduct $\Delta$:
\begin{align*}
&\Delta(\gamma^{\pm 1/2})=\gamma^{\pm 1/2} \otimes \gamma^{\pm 1/2},\\
&\Delta (x^+(z))=
x^+(z)\otimes 1+
\psi^-(\gamma_{(1)}^{1/2}z)\otimes x^+(\gamma_{(1)}z),\\
&\Delta (x^-(z))=
x^-(\gamma_{(2)}z)\otimes \psi^+(\gamma_{(2)}^{1/2}z)+1 \otimes x^-(z),\\
&\Delta (\psi^\pm(z))=
\psi^\pm (\gamma_{(2)}^{\pm 1/2}z)\otimes \psi^\pm (\gamma_{(1)}^{\mp 1/2}z),
\end{align*}
where $\gamma_{(1)}^{\pm 1/2}=\gamma^{\pm 1/2}\otimes 1$
and $\gamma_{(2)}^{\pm 1/2}=1\otimes \gamma^{\pm 1/2}$.\\
Counit $\varepsilon$:
\begin{align*}
\varepsilon(\gamma^{\pm 1/2})=1,\qquad \varepsilon(\psi^\pm(z))=1,
\qquad \varepsilon(x^\pm(z))=0.
\end{align*}
Antipode $a$:
\begin{align*}
&a(\gamma^{\pm 1/2})=\gamma^{\mp 1/2},\\
&a(x^+(z))=-\psi^-(\gamma^{-1/2}z)^{-1}x^+(\gamma^{-1} z),\\
&a(x^-(z))=-x^-(\gamma^{-1} z)\psi^+(\gamma^{-1/2}z)^{-1},\\
&a(\psi^{\pm}(z))=\psi^{\pm}(z)^{-1}.
\end{align*}
\end{prop}

\begin{rmk}
The $\Delta$ is usually called the Drinfeld coproduct.
Strictly speaking, this set of morphisms 
does not define a Hopf algebra, since it gives rise to 
infinite sums. 
\end{rmk}

For later purpose, we introduce the Heisenberg generator $b_n$ as follows.
\begin{lem}\label{lem:b_n}
Define $b_n$ ($n\in\bbZ\setminus\{0\}$) by
\begin{align*}
&\psi^+(z)=\psi^+_0 \exp\left(+\sum_{n>0} b_n \gamma^{n/2} z^{-n} \right) ,
&\psi^-(z)=\psi^-_0 \exp\left(-\sum_{n>0} b_{-n} \gamma^{n/2}z^{n} \right).
\end{align*}
Then we have
\begin{align*}
[b_m,b_n]=\dfrac{1}{m}(1-q^{-m})(1-t^m)(1-q^m t^{-m})(\gamma^m-\gamma^{-m})
\gamma^{-|m|}\delta_{m+n,0},
\end{align*}
and the coproduct for $b_n$ reads
\begin{align*}
\Delta(b_n)=b_n\otimes \gamma_{(2)}^{-|n|}+1\otimes b_n.
\end{align*}
\end{lem}

\subsection{Level zero and level one representations of $\calU(q,t)$}

We can construct two kinds of representations of $\calU$. One is realized on the space of Laurent polynomials $\bbQ(q^{1/2},t^{1/2})[x,x^{-1}]$, and the other is on the Fock space $\calF$. We denote the former by $\pi_x(\cdot)$ and the latter by $\rho(\cdot)$.

\begin{prop}\label{prop:DI_V_x}
We have a representation $\pi_x(\cdot)$ of $\calU(q,t)$ on $V_x\seteq\bbQ(q^{1/2},t^{1/2})[x,x^{-1}]$ by 
setting
\begin{align*}
&\pi_x(\gamma^{\pm 1/2})=1,\\
&\pi_x(\psi^+(z))
 =\dfrac{(1-q^{1/2}t^{-1}x/z)(1-q^{-1/2}t x/z)}{(1-q^{1/2}x/z)(1-q^{-1/2}x/z)}
 =1+\sum_{n\geq 1}(1-q/t)(1-t)\dfrac{1-q^n}{1-q}q^{-n/2}(x/z)^n,
\\
&\pi_x(\psi^-(z))
 =\dfrac{(1-q^{1/2}t^{-1}z/x)(1-q^{-1/2}t z/x)}{(1-q^{1/2}z/x)(1-q^{-1/2}z/x)}
 =1+\sum_{n\geq 1}(1-q/t)(1-t)\dfrac{1-q^{n}}{1-q}q^{-n/2}(z/x)^n,
\\
&\pi_x(x^\pm(z))=c^{\pm 1} (1-t^{\mp 1})\delta(q^{\mp1/2}x/z) T_{q^{\mp 1},x},
\end{align*}
where $c\in \bbQ(q^{1/2},t^{1/2})^{\times}$. 
Here we have used the $q$-shift operator 
$T_{q^{\pm 1},x} f(x)=f(q^{\pm 1} x)$. 
\end{prop}

\begin{prop}\label{prop:DI_Fock_rep}
Let $\frakh$ be the Heisenberg algebra generated by $a_n$ with the relations (\ref{eq:boson_macdonald}), and $\calF$ be the corresponding Fock space. Let $\eta(z),\xi(z),\varphi^+(z)$ and $\varphi^-(z)$ be the vertex operators given in (\ref{eq:vertex_operator}), (\ref{eq:xi}), (\ref{eq:phi^+}) and (\ref{eq:phi^-}).
Then we have a representation $\rho(\cdot)$ of $\calU(q,t)$ on $\calF$ by 
setting
$\rho(\gamma^{\pm 1/2})=(t/q)^{\pm 1/4}$,
$\rho(\psi^\pm(z))=\varphi^\pm(z)$, $\rho(x^+(z))=c\, \eta(z)$
and $\rho(x^-(z))=c^{-1} \xi(z)$,
where $c\in \bbQ(q^{1/2},t^{1/2})^{\times}$. 
\end{prop}

\begin{rmk}
We have identified the Heisenberg generators $a_n$ in (\ref{eq:boson_macdonald})
and $b_n$
introduced in Lemma \ref{lem:b_n} by
\begin{align*}
b_n=-\dfrac{1-t^n}{n}(1-t^n q^{-n})(t/q)^{-n/2}a_n,\quad
b_{-n}=-\dfrac{1-t^{-n}}{n}(1-t^n q^{-n})(t/q)^{-n/2}a_{-n}
\quad(n>0).
\end{align*}
\end{rmk}

When we need to show the dependence of $\pi_x(\cdot)$ on $c$, we denote it
by $\pi_{x,c}(\cdot)$, or more simply denote
the space by $V_{x,c}$. Similarly we use the notations $\rho_c(\cdot)$ and $\calF_c$ for the Fock representation.

For simplicity we call a representation of {\it level} $k$, if 
the central element $\gamma$ is realized by the constant $(t/q)^{k/2}$.
In this terminology, $\pi_x(\cdot)$ is of level zero, and  
 $\rho(\cdot)$ is a level one representation.

\subsection{Intertwining operator $\Phi^{\calF}_{V_x\otimes \calF}$}

In this subsection, 
we study the following intertwining operator 
$\Phi^{\calF}_{V_{x}\otimes \calF}:
V_{x,\alpha}\otimes \calF_\beta \rightarrow \calF_\gamma$
which is subject to the condition $\Phi^{\calF}_{V_{x}\otimes \calF} \Delta (a)=a\, \Phi^{\calF}_{V_{x}\otimes \calF}$
for any $a\in \calU(q,t)$.
We introduce the components $\Phi_n$ of $\Phi^{\calF}_{V_{x}\otimes \calF}$ by
$\Phi^{\calF}_{V_{x}\otimes \calF}( x^n\otimes v):=\Phi_n v$ ($v\in \calF$). 
Set the generating function as $\Phi(y)=\sum_{n\in \bbZ} \Phi_n y^{-n}$.

\begin{prop}
The intertwining property for $\Phi^{\calF}_{V_{x}\otimes \calF}:
V_{x,\alpha}\otimes \calF_\beta \rightarrow \calF_\gamma$ reads
\begin{align}
&\dfrac{(1-q^{1/2}t^{-1}y/z)(1-q^{-1/2}t y/z)}{(1-q^{1/2}y/z)(1-q^{-1/2}y/z)}
\Phi(y) \varphi^+((t/q)^{-1/4}z)
=\varphi^+(z) \Phi( y),\label{eq:Phi_1}\\
&\dfrac{(1-q^{1/2}t^{-1}z/y)(1-q^{-1/2}t z/y)}{(1-q^{1/2}z/y)(1-q^{-1/2}z/y)}
\Phi(y) \varphi^-((t/q)^{1/4}z)
=\varphi^-(z) \Phi( y),\label{eq:Phi_2}\\
&
\alpha (1-t^{-1})\delta(q^{1/2}y/z) \Phi(q y)+
\beta 
\dfrac{(1-q^{1/2}t^{-1}z/y)(1-q^{-1/2}t z/y)}{(1-q^{1/2}z/y)(1-q^{-1/2}z/y)}
\Phi(y) \eta(z)
=\gamma \,\eta(z)\Phi( y),\label{eq:Phi_3}\\
&
\alpha^{-1} (1-t)\delta(t^{-1/2}y/z) 
\Phi(q^{-1} y)\varphi^+((t/q)^{1/4}z)+
\beta^{-1} \Phi(y) \xi(z)
=\gamma^{-1}\xi(z)\Phi( y).\label{eq:Phi_4}
\end{align}
\end{prop}

\begin{proof}
Note that we have $\Phi^{\calF}_{V_{x}\otimes \calF} (\delta(x/y)\otimes v)=\Phi(y) v$ from the definition of the generating function. Writing $\Delta(a)=\sum_i a_i\otimes b_i$ for $a\in \calU(q,t)$, we have
\begin{align*}
\sum_i
\Phi^{\calF}_{V_{x}\otimes \calF}\Big(\pi_{x,\alpha}(a_i) \delta(x/y) \otimes 
\rho_{\beta}(b_i) v\Big)=\rho_{\gamma}(a)\Phi^{\calF}_{V_{x}\otimes \calF}
(\delta(x/y)\otimes v)=\rho_{\gamma}(a) \Phi(y)v
\end{align*}
{}from the intertwining property.
Setting $a=x^+(z)$, we can compute LHS as
\begin{align*}
&\Phi^{\calF}_{V_{x}\otimes \calF}\Big(
\alpha (1-t^{-1})\delta(q^{-1/2}x/z)T_{q^{-1},x}\delta(x/y) \otimes v \Big)
+\Phi^{\calF}_{V_{x}\otimes \calF}\Big(
\pi_{x,\alpha}(\psi^-(z))
\delta(x/y) \otimes \beta \eta(z)v \Big)\\
&=
\alpha (1-t^{-1})\delta(q^{1/2}y/z)
\Phi^{\calF}_{V_{x}\otimes \calF}\Big(
\delta(x/q y) \otimes v \Big)\\
&
\qquad\qquad\qquad\qquad +\beta
\dfrac{(1-q^{1/2}t^{-1}z/y)(1-q^{-1/2}t z/y)}{(1-q^{1/2}z/y)(1-q^{-1/2}z/y)}
\Phi^{\calF}_{V_{x}\otimes \calF}\Big(
\delta(x/y) \otimes \eta(z)v \Big),
\end{align*}
which indicates (\ref{eq:Phi_3}). The rest can be shown in the same manner.
\end{proof}

\begin{prop}\label{prop:intertwiner}
When $\alpha=\gamma$ and $\beta=t^{-1}\gamma$, 
the intertwining operator $\Phi^{\calF}_{V_{x}\otimes \calF}$ uniquely exists
up to normalization, and whose generating function $\Phi(y)$ is realized as
\begin{align}
&\Phi(y)=
\exp
\Bigg(
 \sum_{n>0}\dfrac{1}{n}\dfrac{1-t^n}{1-q^{n}}q^{n/2}t^{-n}a_{-n} y^{n}
\Bigg)
\exp\Bigg(
 -\sum_{n>0}\dfrac{1}{n}\dfrac{1-t^n}{1-q^{n}}q^{n/2}a_{n} y^{-n}
\Bigg).
\label{eq:Phi_boson}
\end{align}
\end{prop}

\begin{proof}
By using (\ref{eq:Phi_1}) and 
(\ref{eq:Phi_2}), one finds that $\Phi(y)$ must be proportional to 
the operator in (\ref{eq:Phi_boson}) if it exists. 
Examining the operator products, we can check 
that both (\ref{eq:Phi_3}) and (\ref{eq:Phi_4})
are fulfilled if and only if we have $\alpha=\gamma$ and $\beta=t^{-1}\gamma$.
\end{proof}


\subsection{Tensor representation $\pi_{x_1}\otimes \cdots \otimes \pi_{x_n}$
and the Macdonald difference operator $D_n^r$}
Let $n$ be a positive integer and consider the $n$-fold tensor space
$V_{x_1}\otimes \cdots \otimes V_{x_n}$.
Using the coproduct, we define the  representation of $\calU(q,t)$ via the composition $\pi_{x_1}\otimes \cdots \otimes \pi_{x_n} \Delta^{(n)}(\cdot)$,
where $\Delta^{(n)}$ is inductively defined to be $\Delta^{(2)}\seteq\Delta$ and $\Delta^{(n)}\seteq({\rm id}\otimes \cdots \otimes{\rm id}
\otimes \Delta)\circ \Delta^{(n-1)}$.
We set the parameter $c^{-1}$ on the $i$-th tensor component by $c^{-1}\seteq t^{n-i}$, namely we set 
$\pi_{x_i}(x^\pm(z)) \seteq
t^{\mp n\pm i} (1-t^{\mp 1})\delta(q^{\mp1/2}x/z) T_{q^{\mp 1},x_i}$.

\begin{prop}\label{prop:level0_tensor}
For any integer $1\leq r\leq n$, we have 
\begin{align*}
&\dfrac{t^{r(r-1)/2}}{(1-t)^r r!}
\Biggl[
 \dfrac{\epsilon_r(z;q)}
       {\prod_{1\leq i<j\leq r} \omega(z_i,z_j;q,t^{-1},q^{-1} t)}\,
 \pi_{x_1}\otimes \cdots \otimes \pi_{x_n} \Delta^{(n)}
 \big( x^-(z_1)\cdots x^-(z_r)\big)
\Biggr]_1\\
&=g D_n^r \,g^{-1},
\end{align*}
where $D_n^r$ denotes the Macdonald difference operator in (\ref{eq:mdo}),
the symbol $\left[\,\cdots\right]_1$ denotes the constant term in $z_i$'s, 
$\omega(x,y;q^{-1},t,q t^{-1})$ is the structure function (\ref{eq:omega}) 
for $\calA$, $\epsilon_r(z;q)\in \calA_r$ as in Definition \ref{dfn:bottom}, 
and $g$ is defined by
\begin{align*}
g=g(x_1,\ldots,x_n)\seteq
\prod_{1\leq i<j\leq n}
\dfrac{(q t^{-1}x_j/x_i;q)_\infty}{(q x_j/x_i;q)_\infty}.
\end{align*}
\end{prop}

\subsection{Quasi-Hopf twisting}

Now we turn to the elliptic deformation of $\calU(q,t)$.
As for the notation, we closely follow the one in \cite{JKOS2}.
Recall the Heisenberg generator $b_n$ introduced in Lemma \ref{lem:b_n}.

\begin{dfn}\label{dfn:dressed_currents}
Set $u^\pm(z;p)\in \calU(q,t)$ by
\begin{align*}
u^+(z;p)\seteq
\exp
\left(-\sum_{n>0}\dfrac{p^n \gamma^{-n}}{1-p^n \gamma^{-2n}}b_{-n}z^n\right),
\quad
u^-(z;p)\seteq
\exp\left(+\sum_{n>0}\dfrac{p^n}{1-p^n }b_{n}z^{-n}  \right),
\end{align*}
and set further 
\begin{align*}
&x^+(z;p) \seteq u^+(z;p) x^+(z),\\
&x^-(z;p) \seteq x^-(z) u^-(z;p) ,\\
&\psi^\pm(z;p) \seteq 
 u^+(\gamma^{\pm 1/2}z;p) \psi^\pm (z) u^-(\gamma^{\mp 1/2}z;p).
\end{align*}
\end{dfn}

One finds that the {\it dressed} Drinfeld currents $x^\pm(z;p), \psi^\pm(z;p)\in \calU(q,t)$ enjoy elliptic permutation relations. For $x^+(z)$ and $x^+(w)$, we have
\begin{align*}
&\Theta_{p\gamma^{-2}}(q^{-1}z/w)
 \Theta_{p\gamma^{-2}}(t z/w)
 \Theta_{p\gamma^{-2}}(q t^{-1}z/w) x^+(z)x^+(w)\\
&=-(z/w)^3 
 \Theta_{p\gamma^{-2}}(q^{-1}w/z)
 \Theta_{p\gamma^{-2}}(t w/z)
 \Theta_{p\gamma^{-2}}(q t^{-1}w/z) x^+(w)x^+(z).
\end{align*}
Recall the elliptic structure function $\omega$ in (\ref{eq:omega_ell}).
The above relation indicates that 
\begin{align*}
\omega(z_1,z_2;q^{-1},t,q t^{-1},p\gamma^{-2})^{-1} x^+(z_1)x^+(z_2)
\end{align*}
gives us a symmetric Laurent series in $z_1$ and $z_2$.
By the same argument as in the main text, 
we can construct, at least formally,  a family of commutative elements in $\calU(q,t)$,
depending on $p$.

\begin{dfn}\label{dfn:F}
We define the twistor $F(p)$ by
\begin{align}
F(p)\seteq
 \exp\left(
  \sum_{n>0}\dfrac{n p^n \gamma_{(2)}^{-n}}
  {(1-q^{-n})(1-t^n)(1-q^n t^{-n})(1-p^n  \gamma_{(2)}^{-2n})}
  b_n\otimes b_{-n}
 \right).
\end{align}
\end{dfn}

\begin{prop}\label{prop:QBA}
(1) The twistor $F(p)$ is invertible, and satisfies 
\begin{align*}
(\varepsilon \otimes {\rm id})F(p)=({\rm id}\otimes \varepsilon) F(p)=1.
\end{align*} 
Hence 
$(\calU(q,t),\Delta_p,\varepsilon, \Phi)$ is a quasi-bialgebra, where we set
\begin{align*}
&\Delta_p (a)\seteq F(p)\cdot \Delta(a)\cdot F(p)^{-1},\\
&\Phi \seteq 
 (F^{(23)}(p)({\rm id}\otimes \Delta)F(p))\cdot 
 (F^{(12)}(p)(\Delta \otimes {\rm id})F(p))^{-1}.
\end{align*}
We denote the resulting quasi-bialgebra by $\calU(q,t,p)$.\\
(2) $F(p)$ satisfies the {\it shifted cocycle condition }
\begin{align*}
F^{(23)}(p)({\rm id}\otimes \Delta)F(p)=
F^{(12)}(p\gamma_{(3)}^{-2})(\Delta \otimes {\rm id})F(p).
\end{align*}
\end{prop}

As for the notation and the definition of the quasi-bialgebra, we refer the readers to \cite{QHA}.
The element $\Phi$ here is called the Drinfeld associator. This should not 
be confused with the intertwining operator $\Phi^\calF_{V_x\otimes \calF}$
or whose generating function $\Phi(y)$. 

A simple calculation yields the following.
\begin{prop}\label{prop:dressed_copro}
We have
\begin{align*}
&\Delta_p (\gamma^{\pm 1/2})=\gamma^{\pm 1/2}\otimes\gamma^{\pm 1/2},\\
&\Delta_p (x^+(z;p))=
 x^+(z;p\gamma_{(2)}^{-2})\otimes 1+
 \psi^-(\gamma_{(1)}^{1/2}z;p\gamma_{(2)}^{-2})\otimes x^+(\gamma_{(1)}z;p),\\
&\Delta_p (x^-(z;p))=
 x^-(\gamma_{(2)}z;p\gamma_{(2)}^{-2})\otimes 
 \psi^+(\gamma_{(2)}^{1/2}z;p)+1 \otimes x^-(z;p),\\
&\Delta_p (\psi^\pm(z;p))=
 \psi^\pm (\gamma_{(2)}^{\pm 1/2}z;p \gamma_{(2)}^{-2})\otimes 
 \psi^\pm (\gamma_{(1)}^{\mp 1/2}z;p).
\end{align*}
\end{prop}

\subsection{Intertwining operator for the elliptic case}

First, we summarize the images of the dressed Drinfeld currents 
$x^\pm(z;p),\psi^\pm(z;p)$ under the level zero and level one 
representations $\pi_{x,c}(\cdot)$ and $\rho_c(\cdot)$.

\begin{lem}
The representation $\pi_{x,c}(\cdot)$ of $\calU(q,t)$ on $V_{x,c}=\bbQ(q^{1/2},t^{1/2})[x,x^{-1}]$ 
written in terms of the dressed Drinfeld currents reads
\begin{align*}
&\pi_{x,c}(\psi^+(z;p))=
\dfrac{\Theta_p(q^{1/2}t^{-1}(x/z)^{\pm 1})\Theta_p(q^{-1/2}t(x/z)^{\pm 1})}
      {\Theta_p(q^{1/2}(x/z)^{\pm 1})\Theta_p(q^{-1/2}(x/z)^{\pm 1})},\\
&\pi_{x,c}(x^\pm(z;p))=c^{\pm 1} (1-t^{\mp 1})
\dfrac{(p t^{\mp 1};p)_\infty(p q^{\mp 1}t^{\pm 1};p)_\infty}
{(p ;p)_\infty(p q^{\mp 1};p)_\infty}
\delta(q^{\mp1/2}x/z) T_{q^{\mp 1},x}.
\end{align*}
\end{lem}

Set
\begin{align}
\label{eq:eta_p}
&\eta(z;p)\seteq
 \exp\Bigg(\sum_{n>0} \dfrac{1-t^{-n}}{n}\dfrac{1-p^n}{1-p^n q^n t^{-n}}
 a_{-n} z^{n}\Bigg)
 \exp\Bigg(-\sum_{n>0} \dfrac{1-t^{n}}{n} a_n z^{-n}\Bigg),\\
\nonumber
&\xi(z;p)\seteq
 \exp\Bigg(-\sum_{n>0} \dfrac{1-t^{-n}}{n}(t/q)^{n/2}a_{-n} z^{n}\Bigg)\\
\nonumber
&\hskip 6em
 \exp\Bigg(\sum_{n>0} \dfrac{1-t^{n}}{n} 
 \dfrac{1-p^n q^n t^{-n}}{1-p^n}
 (t/q)^{n/2} a_n z^{-n}\Bigg),\\
\nonumber
&\varphi^+(z;p) \seteq
 \exp\Bigg(\sum_{n>0} \dfrac{1-t^{-n}}{n}
 (1-t^n q^{-n})\dfrac{p^n }{1-p^n q^n t^{-n}}
 (t/q)^{-3n/4}a_{-n} z^{n}\Bigg)\\
\nonumber
&\hskip 6em
 \times 
 \exp\Bigg(
   -\sum_{n>0} \dfrac{1-t^{n}}{n} (1-t^n q^{-n})
   \dfrac{1}{1-p^n}(t/q)^{-n/4} a_n z^{-n}
 \Bigg),\\
\nonumber
&\varphi^-(z;p)\seteq
 \exp\Bigg(
  \sum_{n>0} \dfrac{1-t^{-n}}{n}
  (1-t^n q^{-n})\dfrac{1}{1-p^n q^n t^{-n}}
  (t/q)^{-n/4}a_{-n} z^{n}
 \Bigg)\\
\nonumber
&\hskip 6em
 \times 
 \exp\Bigg(
  -\sum_{n>0} \dfrac{1-t^{n}}{n}(1-t^n q^{-n}) 
              \dfrac{p^n}{1-p^n}(t/q)^{-3n/4}
  a_n z^{-n}
 \Bigg).
\end{align}
Note that we have $\varphi^-(z;p) = \varphi^+(z p^{-1}q^{-1/2}t^{1/2};p)$.

\begin{lem}
The representation $\rho_c(\cdot)$ on the Fock space $\calF$ gives
the images of dressed Drinfeld currents as
\begin{align*}
\rho_c(x^+(z;p))=c\eta(z;p),\quad
\rho_c(x^-(z;p))=c^{-1}\xi(z;p),\quad
\rho_c(\psi^\pm(z;p))=\varphi^\pm(z;p).
\end{align*}
\end{lem}

Consider the intertwining operator $\Phi^{\calF}_{V_{x}\otimes \calF}(p): V_{x,\alpha}\otimes \calF_\beta \rightarrow \calF_\gamma$ with respect to the elliptic 
counterpart $\calU(q,t,p)$ of the Ding-Iohara algebra.
Namely, it
satisfies the condition $\Phi^{\calF}_{V_{x}\otimes \calF}(p) \Delta_p (a)=a\, \Phi^{\calF}_{V_{x}\otimes \calF}(p)$ for any $a\in \calU(q,t)$. Introduce the components $\Phi_{p,n}$ of $\Phi^{\calF}_{V_{x}\otimes \calF}(p)$ by $\Phi^{\calF}_{V_{x}\otimes \calF}(p)( x^n\otimes v)=\Phi_{p,n} v$ ($v\in \calF$), and set the generating function as $\Phi(y;p)\seteq\sum_{n\in \bbZ} \Phi_{p,n} y^{-n}$.

\begin{prop}\label{prop:ITP_p}
The intertwining property for $\Phi^{\calF}_{V_{x}\otimes \calF}(p):
V_{x,\alpha}\otimes \calF_\beta \rightarrow \calF_\gamma$ reads
\begin{align*}
&\dfrac
 {\Theta_{p q t^{-1}}(q^{1/2}t^{-1}y/z)\Theta_{p q t^{-1}}(q^{-1/2}t y/z)}
 {\Theta_{p q t^{-1}}(q^{1/2}y/z)\Theta_{p q t^{-1}}(q^{-1/2}y/z)}
 \Phi(y;p) \varphi^+((t/q)^{-1/4}z;p)
 =\varphi^+((t/q)^{-1/4}z,p) \Phi(y;p),\\
&\dfrac
 {\Theta_{p q t^{-1}}(q^{1/2}t^{-1}z/y)\Theta_{p q t^{-1}}(q^{-1/2}t z/y)}
 {\Theta_{p q t^{-1}}(q^{1/2}z/y)\Theta_{p q t^{-1}}(q^{-1/2}z/y)}
 \Phi(y;p) \varphi^-((t/q)^{1/4}z;p)
 =\varphi^-((t/q)^{1/4}z;p) \Phi(y;p),\\
&\alpha (1-t^{-1})
 \dfrac{(p q t^{-2};p q t^{-1})_\infty (p;p q t^{-1})_\infty}
       {(p q t^{-1};p q t^{-1})_\infty (p t^{-1};p q t^{-1})_\infty}
 \delta(q^{1/2}y/z) \Phi(q y;p)\\
&\hskip 2em 
+\beta 
 \dfrac{\Theta_{p q t^{-1}}(q^{1/2}t^{-1}z/y)\Theta_{p q t^{-1}}(q^{-1/2}t z/y)}       {\Theta_{p q t^{-1}}(q^{1/2}z/y)\Theta_{p q t^{-1}}(q^{-1/2}z/y)}
 \Phi(y;p) \eta(z;p)
=\gamma \,\eta(z;p)\Phi(y;p),\\
&
 \alpha^{-1} (1-t)
 \dfrac{(p q;p q t^{-1})_\infty (p q^2 t^{-2};p q t^{-1})_\infty}
       {(p q t^{-1};p q t^{-1})_\infty (p q^2 t^{-1};p q t^{-1})_\infty}
 \delta(t^{-1/2}y/z)\Phi(q^{-1} y;p)\varphi^+((t/q)^{1/4}z;p)\\
&\hskip 2em 
 +\beta^{-1} \Phi(y;p) \xi(z;p) 
=\gamma^{-1}\xi(z;p)\Phi(y;p).
\end{align*}
\end{prop}

\begin{prop}
When $\alpha=\gamma$ and $\beta=t^{-1}\gamma$, the intertwining operator $\Phi^{\calF}_{V_{x}\otimes \calF}(p)$ uniquely exists up to normalization, and whose generating function $\Phi(y;p)$ is realized as
\begin{align}
\Phi(y;p)=
\exp\Bigg(\sum_{n>0}
 \dfrac{1}{n}\dfrac{1-t^n}{1-q^{n}}\dfrac{1-p^n}{1-p^n q^n t^{-n}}
 q^{n/2}t^{-n}a_{-n} y^{n}
\Bigg)
\exp\Bigg(-\sum_{n>0}\dfrac{1}{n}\dfrac{1-t^n}{1-q^{n}}q^{n/2}a_{n} y^{-n}
\Bigg).
\label{eq:Phi_boson2}
\end{align}
\end{prop}

\subsection{Deformed $\calW_m$ algebra}

We end this appendix by pointing out the connection between 
the tensor representation of the Ding-Iohara algebra and 
the deformed $\calW_m$ algebra found in \cite{SKAO:1996,FF,AKOS:1996}.

Let $m\geq 2$ be a positive integer. Consider the $m$-fold tensor representation
$\rho_{y_1}\otimes\cdots \otimes\rho_{y_m}$ on 
$\calF_{y_1}\otimes \cdots \otimes \calF_{y_m}$.
Define $\Delta^{(m)}_p$ inductively by
$\Delta^{(2)}_p\seteq \Delta_p$  and 
$\Delta^{(m)}_p\seteq({\rm id}\otimes \cdots 
\otimes{\rm id}\otimes \Delta_p)\circ \Delta^{(m-1)}_p$. 
Since we have $\rho_{y_1}\otimes\cdots \otimes\rho_{y_m}\Delta^{(m)}_p(\gamma)=
\gamma_{(1)}\cdots \gamma_{(m)}=
(t/q)^{m/2}$, the level is $m$. 
By abuse of notation we write $\gamma=(t/q)^{1/2}$ for simplicity.

\begin{lem}
We have 
\begin{align*}
\rho_{y_1}\otimes\cdots \otimes\rho_{y_m} \Delta^{(m)}_p(x^+(z;p))
=\sum_{i=1}^m y_i \Lambda_i(z).
\end{align*}
Here the $\Lambda_i(z)$ is defined by the tensor
\begin{align*}
\Lambda_i(z)
&\seteq
 \varphi^-(\gamma^{1/2}z;p\gamma^{-2m+2})\otimes
 \varphi^-(\gamma^{3/2}z;p\gamma^{-2m+4})\otimes\cdots\\
&\hskip 3em 
 \cdots \otimes 
 \varphi^-(\gamma^{(2i-3)/2}z;p\gamma^{-2m+2i-2})
 \otimes \eta(\gamma^{i-1}z;p \gamma^{-2m+2i})\otimes 1\otimes \cdots\otimes 1,
\end{align*}
where $\eta(\gamma^{i-1}z;p \gamma^{-2m+2i})$
sits in the $i$-th tensor component.
\end{lem}

Set 
\begin{align}
\label{eq:f(z)}
&f(z)\seteq
 \exp\left(
  \sum_{n>0}
   \dfrac{1}{n}
   \dfrac{(1-q^n)(1-t^{-n})(1-p^n q^{(m-1)n}t^{-(m-1)n})} 
         {1-p^n q^{m n}t^{-m n}}z^n \right),
\\
\nonumber
&
\gamma_{i,j}(z)\seteq
\left\{
\begin{array}{ll}
1&(i=j),\\[2mm]
\dfrac{(1-q^{-1}z)(1-t z)}{(1-z)(1-q^{-1}t z)} &(i<j),\\[5mm]
\dfrac{(1-q z)(1-t^{-1}z)}{(1-z)(1-q t^{-1}z)} &(i>j).
\end{array}\right.
\end{align}

\begin{prop}\label{prop:W_m}
We have the operator product 
\begin{align*}
f(w/z) \Lambda_i(z)\Lambda_j(w)=\gamma_{i,j}(w/z) :\Lambda_i(z)\Lambda_j(w):
\end{align*}
for $1\leq i,j\leq m$.
\end{prop}

We note that
this operator product formula coincide with 
the one for the deformed $\calW_m$ algebra.
To be more precise we need to put $p=1$ in (\ref{eq:f(z)}) to recover the
original structure function $f(z)$ for the deformed $\calW_m$ algebra.

\begin{ack}
The authors thank Professor M.~Noumi
for valuable comments, discussion, and
allowing to include his private result in Proposition \ref{prop:Noumi}. 
K.~H. and S.~Y. would like to express 
their gratitude to him 
for inviting them in this subject.
\end{ack}

\end{document}